\documentclass[fleqn,10pt]{article}
\usepackage[utf8]{inputenc}

\usepackage{graphicx}
\usepackage{float} 
\usepackage{subfigure} 
\usepackage{amsmath,amssymb}
\usepackage{appendix}
\usepackage[hidelinks]{hyperref}
\usepackage{url}
\usepackage{comment,xcolor}
\usepackage[lined,boxed,linesnumbered,ruled]{algorithm2e}
\usepackage{bm}

\usepackage{setspace}
\usepackage{enumitem}
\setlist{nolistsep}

\usepackage{natbib}
\newcommand{\mbs}[1]{\mathbf{#1}}
\newcommand{\mbb}[1]{\mathbb{#1}}

\newtheorem{remark}{Remark}[section]

\title{Phase Field Modeling and Numerical Algorithm for Two-Phase Dielectric Fluid Flows}
\author{Jielin Yang$^1$, Ivan C.~Christov$^2$, Suchuan Dong$^1$\thanks{Author of correspondence. Emails: yang1659@purdue.edu (Yang), christov@purdue.edu (Christov), sdong@purdue.edu (Dong)} \\
$^1$Department of Mathematics,\\
$^2$School of Mechanical Engineering, \\
Purdue University, West Lafayette, IN47907, USA
}
\date{(May 15, 2023)}

\begin{document}

\maketitle

\begin{abstract}

We develop a method for modeling and simulating a class of two-phase flows consisting of two immiscible incompressible dielectric fluids and their interactions with imposed external electric fields in two and three dimensions. We first present a thermodynamically-consistent and reduction-consistent phase field model for two-phase dielectric fluids. The model honors the conservation laws and thermodynamic principles, and has the property that, if only one fluid component is present in the system, the two-phase formulation will exactly reduce to that of the corresponding single-phase system. In particular, this model accommodates an equilibrium  solution that is compatible with the zero-velocity requirement based on physics. This property provides a simpler method for simulating the equilibrium state of two-phase dielectric systems. We further present an efficient numerical algorithm, together with a spectral-element (for two dimensions) or a hybrid Fourier-spectral/spectral-element (for three dimensions) discretization in space, for simulating this class of problems. This algorithm computes different dynamic variables successively in an un-coupled fashion, and involves only coefficient matrices that are time-independent in the resultant linear algebraic systems upon discretization, even when the physical properties (e.g.~permittivity, density, viscosity) of the two dielectric fluids are different. This property is crucial and enables us to employ fast Fourier transforms for three-dimensional problems. Ample numerical simulations of two-phase dielectric flows under imposed voltage are presented to demonstrate the performance of the method herein and to compare the simulation results with theoretical models and experimental data.

\end{abstract}

\vspace{0.05cm}

\section{Introduction}



In the current work we focus on the modeling and simulation of a system of two immiscible incompressible dielectric fluids and their interaction with external electric fields. 
Dielectric fluids refer to fluids that are electrically non-conductive and can withstand high voltages without breakdown. They are traditionally used for cooling and insulating electrical equipment such as transformers and high-voltage cables. In recent years dielectric fluids have found widespread applications in electric vehicles, for cooling the electric motor, batteries, electric transmissions, and power electronics.

Using electric field to control fluid interface or droplets is a  widely-used technique for manipulating small amounts of liquids on surfaces.  Electrowetting-on-dielectric (EWOD)~\cite{mugele2005electrowetting}  is one of the most successful and  versatile approaches. The EWOD system typically involves conducting fluids or droplets on a dielectric substrate under an imposed voltage. The applications of EWOD range from ``lab-on-a-chip" devices~\cite{srinivasan2004integrated,cooney2006electrowetting}, to adjustable lenses~\cite{berge2000variable}, to  new types of electronic displays~\cite{hayes2003video}. While EWOD is versatile and effective in fluid manipulation, it requires the fluids to be conductive and an AC (alternating-current) electric potential.

Using dielectrophoresis (DEP)  underlies another class of  techniques for manipulating  fluids~\cite{batchelder1983dielectrophoretic,gascoyne2004dielectrophoresis}, and has an advantage over EWOD in some situations~\cite{edwards2018dielectrowetting}. Dielectrophoresis refers to the electromechanical force due to the polarization of a neutral material in non-uniform electric fields~\cite{edwards2018dielectrowetting}.
One can use the Korteweg-Helmholtz force density~\cite{landau2013electrodynamics} to explain the origination of liquid dielectrophoresis. Applying an electric field $\textbf{E}$ to the fluid results in   the Korteweg-Helmholtz force density
\begin{equation}\small
    \mbs f_{KT}=\rho_f\textbf{E}-\frac{1}{2}(\textbf{E}\cdot\mbs E)\nabla\epsilon+\nabla\left[\frac{1}{2}(\textbf{E}\cdot\mbs E)\rho\frac{\partial \epsilon}{\partial \rho}\right], \label{KT_force}
\end{equation}
where $\rho_f$ is the free electric charge density, $\epsilon$ is the permittivity, and $\rho$ is the fluid density. Based on equation~\eqref{KT_force}, when there is no free charge ($\rho_f=0$), as long as nonuniform polarisation of dipoles exists  within the liquid ($\nabla\epsilon\neq0$), the fluid will be influenced by the electric field; see \cite{geng2017dielectrowetting,pohl1951motion} for more detailed discussions of the DEP force theory.

The use of DEP to move bulk fluids can be traced to the work of Pellat in 1895~\cite{edwards2018dielectrowetting}. A dielectric siphon is described in~\cite{JonesPM1971} to pump fluids between two reservoirs. DEP is used to transport dielectric particles or droplets in~\cite{batchelder1983dielectrophoretic}  through a channel sandwiched by electrodes.
Transport of dielectric liquids at microscale and in microfluidic devices has been studied  with miniaturized electrodes in~\cite{JonesGWF2001,Jones2001,ChughK2008}.
In~\cite{brown2009voltage} the authors use DEP to spread a droplet onto coplanar electrodes to form a thin liquid film, estalishing the idea of interface localized liquid dielectrophoresis. In~\cite{mchale2011dielectrowetting} the effect of localized DEP on the wetting properties of solid-liquid interface has been investigated and the term dielectrowetting is introduced.
We refer to e.g.~\cite{xu2013dielectrophoretically,geng2017dielectrowetting,mchale2012developing} for a  review of this area and recent applications.

The two-phase system of dielectric fluids involves fluid interfaces, the associated surface tension, the contrast in  fluid properties (permittivity, density, and viscosity), contact lines and contact angles when a solid-wall boundary is  present, and the interaction with the imposed electric field. 
The  approach taken in the current work to handle the two phases belongs to the phase field framework.
Phase field (a.k.a. diffuse interface)~\cite{Rayleigh1892,Waals1893,AndersonMW1998,LowengrubT1998,Jacqmin1999,boyer2002theoretical,LiuS2003,ding2007diffuse} is one of the few techniques currently available for dealing with two-phase systems and fluid interfaces. It is particularly attractive because of its physics-based nature. With phase field the fluid interface is treated as a thin smooth transition layer (i.e.~diffuse). Besides the hydrodynamic variables, the system is characterized by an order parameter (or phase field function), which varies smoothly within the transition layer and is mostly uniform in the bulk phases. 
The evolution of the fluid phases is characterized by a free energy density function, which contain component terms that tend to promote the mixing of the two fluids and also component terms that tend to separate the fluids. The balance and interplay of these two tendencies determine
the dynamic profile of the fluid interface.
With this approach, the governing equations of the system can be derived rigorously based on the conservation laws and thermodynamic principles. We refer to e.g.~\cite{LowengrubT1998,KimL2005,abels2012thermodynamically,ShenYW2013,AkiDG2014,dong2014efficient,LiuSY2015,GongZW2017,Dong2018,RoudbariSBZ2018,Yue2020} (among others) for several thermodynamically consistent phase field models for two-phase and multiphase flows with various degrees of sophistication.

While  phase field  is successful for a range of two-phase and multiphase problems, investigations into this approach for modeling two-phase hydrodynamics coupled with the electric field effect are still quite limited. 
In~\cite{lin2012phase,yang20133d,yang2014phase} the authors employ the phase field method coupled with the Navier-Stokes equations to study the electrohydrodynamic (EHD) phenomenon, in particular the Taylor's leaky dielelctric model~\cite{saville1997electrohydrodynamics}. The authors of~\cite{tian2013numerical} investigate the electrohydrodynamic patterning based on the liquid dielectrophoresis. In \cite{wang2016numerical,xie2016two} the phase field method is used to study electrowetting and its applications.

The aforementioned studies on the coupled multiphase flow and electric field have a notable drawback.
These are  phenomenological models, and do not admit an energy law (or energy balance relation). In other words, these models are not  thermodynamically consistent.
To overcome this issue,
a phase field model is developed in~\cite{eck2009phase} for electrowetting (conductive fluids with free charges) based on the variational principles and the thermodynamics of irresversible processes near equilibrium. The model combines the 
multiphase flow, the electric field and the free charge system, 
and admits an energy balance relation. However, it only applies to cases when the two conductive fluids have the same density.

This model is extended in~\cite{CampilloGK2012} to take into account the density contrast and the transport of free ion species in the conductive fluids; see~\cite{Metzger2015,linga2018controlling,Metzger2019,LingaBM2019} for numerical algorithms developed based on this extended model.
Another diffuse interface model is proposed in~\cite{nochetto2014diffuse} for electrowetting on dielectric with different densities for the two fluids, which however appears not to be  Galilean invariant.
In~\cite{LingaBM2020} a thermodynamically consistent continuum model for single-phase electrohydrodynamic flows has been described. The model combines the Navier-Stokes equations and the Poisson-Nernst-Planck (PNP) equations, in which the fluid properties depend on the ion concentration fields.
We would also like to note the finite element method developed in~\cite{ZhaoR2021} employng a sharp-interface model for electrowetting on dielectric.


In the current work we look into
the dynamics of an isothermal system of two immiscible incompressible dielectric fluids and their interaction with
imposed external electric fields in two and three dimensions (2D/3D). The fluids considered here are  non-conductive and the system contains no free charges or ions. This setting is quite different from those studies reviewed in previous paragraphs related to electrowetting or electrohydrodynamics, where the fluids are electrolytic solutions and conductive and the transport of free ions is crucial to the system dynamics. Due to the liquid dielectrophoresis and the Korteweg-Helmholtz force, when an external voltage is applied, the interface between the dielectric fluids can experience large deformations, leading to the dielectrowetting phenomenon~\cite{edwards2018dielectrowetting}.

We first present a thermodynamically-consistent and reduction-consistent phase field model for two-phase dielectric fluid flows. Thermodynamic consistency refers to the property that the model honors the conservation laws and thermodynamic principles.
The current model is developed based on the conservations of mass and momentum and the second law of thermodynamics, in which the physical properties of the two fluids (permittivity, density, and viscosity) can be different. The model derivation process follows those of~\cite{abels2012thermodynamically,dong2014efficient,Dong2018}, with the quasi-static electromagnetic equations taken into account. 
Reduction consistency refers to the property that, when only one fluid component is present in the two-phase system (while the other fluid is absent),
the two-phase formulation will  exactly reduce to that of the corresponding single-phase system. We refer to~\cite{Dong2018} for discussions of  reduction consistency in general multiphase systems. The reduction consistency of a two-phase dielectric system places restrictions on the functional form of the mixture permittivity when expressed in terms of the phase field variable.
As discussed in~\cite{Dong2018}, reduction consistency reflects an inherent reduction relation within multiphase systems, and violation of  reduction consistency can lead to unphysical results from a model (e.g.~production of a fluid phase where it is absent). 
The phase field model here for dielectric fluids appears to have some connection to that of~\cite{CampilloGK2012} for conductive fluids and electrolytic solutions. We note that the reduction consistency issue was not considered in~\cite{CampilloGK2012} or
the related works of~\cite{Metzger2015,linga2018controlling,Metzger2019,LingaBM2019} for conductive fluids, and the model as given therein appears not reduction-consistent.

At equilibrium, the solution to the current model is compatible with the zero-velocity requirement based on physics. This property provides a  method for computing the equilibrium state (or steady state) of two-phase dielectric systems, which is of great practical interest and importance (e.g.~the equilibrium shapes of dielectric droplets under  imposed voltage), by solving a smaller reduced system of equations. This method is simpler and faster than integrating the full model in time  until the steady state is reached.

We then present a semi-implicit splitting type algorithm, together with a spectral-element spatial discretization for 2D and a hybrid Fourier-spectral/spectral-element discretization for 3D, for numerically solving the governing equations of this two-phase dielectric flow model. The computations for different dynamic variables (electric potential, phase field function, velocity, and pressure) are de-coupled with our method. For each dynamic variable, the resultant linear algebraic system upon discretization involves a constant and time-independent coefficient matrix, which can be pre-computed and saved for later use, despite the variable physical properties (permittivity, density, viscosity) of the two-phase mixture. This characteristic of the current algorithm is crucial, and it enables the use of Fourier spectral discretization and fast Fourier transform (FFT) in 3D simulations of two-phase dielectric flows with variable mixture properties. For 3D problems, with each dynamic variable, the computations of different Fourier modes are completely de-coupled with the current method. Thanks to these characteristics, the presented method is computationally very efficient.

These attractive properties of the current method are attained based on several strategies. The most important strategy, for producing a semi-discretized system having constant coefficients when variable material properties are present on the continuum level, is inspired by and built upon the algorithm from~\cite{dong2012time} (with modifications).
The main idea of~\cite{dong2012time} lies in a reformulation of the pressure/viscous terms in the momentum equation as follows,
\begin{equation*}\small
\begin{aligned}
&\frac{1}{\rho}\nabla p\approx \frac{1}{\rho_0}\nabla p+\left(\frac{1}{\rho}-\frac{1}{\rho_0}\right)\nabla p^*,
\qquad
\frac{\mu}{\rho}\nabla^2\textbf{u}\approx\nu_m\nabla^2\textbf{u}-\left(\frac{\mu}{\rho}-\nu_m\right)\nabla\times\nabla\times\textbf{u}^*,
\end{aligned}
\end{equation*} 
where $\rho$ and $\mu$ are  the variable density and variable viscosity of the mixture, $\rho_0$ and $\nu_m$ are two appropriate constants, $p$ and $\textbf{u}$ are the pressure and the divergence-free velocity, and $p^*$ and $\textbf{u}^*$ are explicit approximations of $p$ and $\mbs u$ with a prescribed order of accuracy. Such and similar reformulations  lead to a semi-discretized system of equations with constant coefficients, in spite of the variable fluid properties on the continuum level. 
This semi-discretized system with {\em constant} coefficients is critical to the success of subsequent Fourier spectral discretization in one or more directions in 3D space. This is because, if this system has variable coefficients, the FFT will induce convolutions with the coefficient functions in the frequency space, which will couple together all the Fourier modes of the unknown dynamic variables to be solved for.


The current algorithm is semi-implicit in nature, in which the nonlinear terms involved in the governing equations are treated explicitly and the linear terms are treated implicitly. As such the algorithm is only conditionally stable, in the sense that the time step size employed in the simulations cannot be large. On the other hand, this conditional stability in the algorithm enables very efficient computations within each time step, with constant pre-computable coefficient matrices and de-coupled computations for all dynamic variables. 
As opposed to the semi-implicit approach, one may also consider the development of energy-stable type schemes, which is not pursued here. Energy-stable schemes are discretizations designed to satisfy a discrete version of the energy law, irrespective of the time step size. The strength of energy-stable schemes lies in that they are  unconditionally stable and can allow the use of relatively large time step sizes in the simulations. We refer to e.g.~\cite{ShenY2010,Salgado2013,GuoLL2014,GrunK2014,ShenY2015,GuoLLW2017,YuY2017,RoudbariSBZ2018} (among  others) for several energy-stable schemes for two-phase problems.
The downside of the energy-stable algorithms lies in that their computational cost per time step can be very high. These schemes typically involve the solution of coupled nonlinear algebraic equations or coupled linear algebraic equations, and the linear algebraic systems resultant from these schemes involve time-dependent coefficient matrices, which require frequent re-computations (e.g.~at every time step). 


We present a number of numerical examples of two-phase dielectric flows under an imposed voltage in 2D and 3D to test the performance of the presented method. In particular, we compare the current simulation results with the theoretical models and the experimental data from the literature. The comparisons show that the phase field model and the numerical method developed herein can capture the physics of this class of flow problems  well.


The contributions of this paper lie in three aspects: (i) the reduction-consistent and thermodynamically-consistent phase field model for two-phase dielectric fluids, (ii) the simpler method for computing the equilibrium state of two-phase dielectric systems, and (iii) the efficient numerical algorithm for simulating two-phase dielectric flows.

The rest of this paper is organized as follows.
In Section~\ref{section2} we present the phase field model for two-phase dielectric flows and discuss the boundary/initial conditions and the equilibrium solution to this model.
In Section~\ref{section3} we present the numerical algorithm for solving the governing equations of this model, and discuss the spectral-element implementation for 2D problems and the hybrid Fourier-spectral and spectral-element implementation for 3D problems.
We employ several 2D and 3D two-phase dielectric flows to test the presented method in Section~\ref{sec:tests}, and in particular we compare the simulation results with theoretical models and the experimental data.
Section~\ref{sec:summary} concludes the presentation with some closing remarks. In Appendix A we outline the development of the current phase field model based on the conservation laws and thermodynamic principles.


\section{Phase Field Model for Two-Phase Dielectric Fluids}\label{section2}

Consider a domain $\Omega$ in two or three dimensions, and an isothermal system of two immiscible incompressible dielectric fluids in this domain.
The two fluids are assumed to be Newtonian, with constant densities $\rho_1$ and $\rho_2$, constant dynamic viscosities $\mu_1$ and $\mu_2$, and constant relative permittivity $\epsilon_1$ and $\epsilon_2$,
respectively.
We introduce a phase field variable $\phi$, which assumes the constant values $1$ and $-1$ in the bulk of the two fluids and has a smooth distribution in a thin layer surrounding the interface.

The material properties of the mixture are functions of the above parameters and the phase field variable $\phi$, with the mixture density $\rho=\rho(\rho_1,\rho_2,\phi)$, mixture viscosity $\mu=\mu(\mu_1,\mu_2,\phi)$, and mixture permittivity  $\epsilon=\epsilon(\epsilon_1,\epsilon_2,\phi)$. 
Specifically, we assume the following relations,
\begin{equation}\small
\left\{
\begin{aligned}
&\rho(\phi) =  \dfrac{\rho_1+\rho_2}{2}+\dfrac{\rho_1-\rho_2}{2}\phi, \qquad
\mu(\phi)=    \dfrac{\mu_1+\mu_2}{2}+\dfrac{\mu_1-\mu_2}{2}\phi, \\
&\epsilon(\phi)=  \dfrac{\epsilon_1+\epsilon_2}{2}+\dfrac{\epsilon_1-\epsilon_2}{2}\dfrac{\phi(3-\phi^2)}{2}.
\end{aligned}
\right.
\label{interpolation}
\end{equation} 
In the above relations, $\rho$ and $\mu$ are assumed to be linear with respect to $\phi$, which has been commonly used (see e.g.~\cite{ding2007diffuse,Dong2015}). However, for $\epsilon$ we employ a relation based on the Hermite interpolation. The benefit of Hermite interpolation is that $\dfrac{d\epsilon}{d\phi}=0$ in the bulk ($\phi=\pm1$), while a linear relation would result in a non-zero derivative.
The zero derivative of permittivity plays an important role in our modeling, which will become clearer in later discussions. The derivative of the mixture permittivity is,
$ 
\epsilon'(\phi)=\frac{\epsilon_2-\epsilon_1}{2}\frac{3(\phi^2-1)}{2}.
$ 

\subsection{Governing Equations}

The phase field model describing the motion of this system of fluids can be derived based on the conservation laws and thermodynamic principles. The development of this model has been discussed in detail in Appendix A. Here we only summarize the governing equations for this system.

Let $\mbs u$ denote the velocity, $\mathcal P$ denote the pressure, $\phi$ denote the phase field variable, $V$ denote the electric potential, and $\mbs E$ denote the electric field. 
Then the dynamics of this two-phase system is described by the following set of equations (see Appendix A for the derivation and specifically~\eqref{eq_124} for the general form of equations), 
\begin{align}\small
    &\dfrac{\partial \phi}{\partial t}+\textbf{u}\cdot\nabla\phi=\gamma_1\Delta\left(\lambda h(\phi)-\lambda \Delta\phi-\dfrac{\epsilon'}{2}\textbf{E}\cdot\mbs E\right),
    \label{eq_6}\\
    &\rho\left(\dfrac{\partial \mbs u}{\partial t}+\mbs u\cdot\nabla\mbs u\right)+\Tilde{\mbs J}\cdot\nabla\mbs u=-\nabla\cdot \left( \lambda\nabla \phi\otimes \nabla\phi \right)-\dfrac{\nabla\epsilon}{2}\textbf{E}\cdot\mbs E+\nabla\cdot\left[\mu\left(\nabla\mbs u +\nabla\mbs u^T\right)\right]-\nabla\mathcal P, \label{eq_7}
    \\
    &\nabla\cdot\textbf{u}=0, \label{eq_8}
    \\
    &\nabla\cdot(\epsilon\nabla V)=0, \label{eq_9} \\
    &
    \mbs E = \nabla V, \label{eq_a8}
\end{align}
where the flux term $\Tilde{\textbf{J}}$ is given by
\begin{equation}\small
    \Tilde{\textbf{J}}=-\gamma_1\dfrac{\rho_1-\rho_2}{2}\nabla\left(\lambda h(\phi)-\lambda\nabla^2\phi-\dfrac{\epsilon'}{2}\textbf{E}\cdot\mbs E\right).\label{J}
\end{equation}
In these equations $\gamma_1$ is the mobility coefficient, and $\lambda$ is the mixing energy density coefficient. 
$\rho$, $\mu$ and $\epsilon$ denote the density, dynamic viscosity, and permittivity of the mixture and are given in~\eqref{interpolation}.
$h(\phi)$ in equation~\eqref{eq_6} is defined by $\lambda h(\phi)=\dfrac{\partial F}{\partial \phi}$, where $F(\phi)$ is the interfacial mixing energy density function (with double well) given by, 
\begin{equation}\small
    F(\phi,\nabla\phi)=\dfrac{1}{2}\lambda|\nabla\phi|^2+\dfrac{\lambda}{4\eta^2}(\phi^2-1)^2.\label{double_well}
\end{equation}
The constant $\eta$ here is a length scale characterizing the interfacial thickness, and
$\lambda$ is related to the surface tension $\sigma$ of 
the two phases  by 
$ 
    \lambda=\dfrac{3}{2\sqrt{2}}\sigma \eta
$ 
\cite{YueFLS2004}.
So $h(\phi)$ is given by,
$ 
h(\phi)=\dfrac{1}{\eta^2}\phi(\phi^2-1).
$ 

With $\Tilde{\textbf{J}}$ given by~\eqref{J} and $\rho$ given in~\eqref{interpolation}, equation~\eqref{eq_6} is equivalent to,
\begin{equation}\small
    \dfrac{\partial \rho}{\partial t}+\textbf{u}\cdot\nabla\rho=-\nabla\cdot\Tilde{\textbf{J}}.
\end{equation}
Let $\mu_c$ denote a generalized chemical potential given by
\begin{equation}\small
\mu_c=\lambda h(\phi)-\lambda\nabla^2\phi-\dfrac{\epsilon'}{2}\textbf{E}\cdot\mbs E. \label{ecp}
\end{equation}
Then, $\Tilde{\textbf{J}}$ can be written as, $\Tilde{\textbf{J}}=-\gamma_1\frac{\rho_1-\rho_2}{2}\nabla\mu_c$.

\subsection{Reduction Consistency}
\label{rem_a0}

We require that the system consisting of~\eqref{eq_6}--\eqref{eq_a8} should be reduction-consistent~\cite{Dong2018,Dong2017}. In other words, if only one fluid component is present (while the other fluid is absent), the system of two-phase governing equations should exactly reduce to that of the corresponding single-phase system. This means that the system given by~\eqref{eq_6}--\eqref{eq_a8} should admit the following two solutions: 
\begin{itemize}
\item $(\mbs u, \mathcal{P}, V, \mbs E)$ and $\phi\equiv 1$: the first fluid is present, and the second fluid is absent. 

\item $(\mbs u, \mathcal{P}, V, \mbs E)$ and $\phi\equiv -1$: the second fluid is present, and the first fluid is absent.
\end{itemize}
It can be verified that these solutions are ensured if the following conditions on $\epsilon(\phi)$ are satisfied,
\begin{equation}\label{eq_a12}\small
\left.\frac{d\epsilon}{d\phi}\right|_{\phi=1}=0, \qquad
\left.\frac{d\epsilon}{d\phi}\right|_{\phi=-1}=0.
\end{equation}
The choice for $\epsilon(\phi)$ in~\eqref{interpolation} satisfies these conditions. Therefore the phase field model given by~\eqref{eq_6}--\eqref{eq_a8}, with the mixture properties given by~\eqref{interpolation}, is reduction consistent. 
It is noted that if one chooses a linear form for $\epsilon(\phi)$ (similar to $\rho(\phi)$ and $\mu(\phi)$ in~\eqref{interpolation}), then the system~\eqref{eq_6}--\eqref{eq_a8} will not be reduction consistent (when $\epsilon_1\neq\epsilon_2$).
We refer to~\cite{Dong2018} for more detailed discussions of the reduction consistency for multiphase systems.

From the physics perspective, 
the electric field influences the generalized chemical potential through the term $\dfrac{\epsilon'(\phi)}{2}\mbs E\cdot\mbs E$. Physically, the generalized chemical potential in the phase field equation should only have an effect on the interface (not in the bulk region), i.e.~$\mu_c$ should vanish in the bulk. This leads to the same conditions as given in~\eqref{eq_a12}.
Therefore, the Hermite interpolation relation for $\epsilon(\phi)$ in equation~\eqref{interpolation} is crucial to the current model.

\subsection{Energy Law}

The model given by equations~\eqref{eq_6}--\eqref{eq_a8} admits an energy law. Let $E(t)$
denote the total system energy, 
\begin{equation}\small
    E(t)=\int_{\Omega} \left(\dfrac{1}{2}\rho \textbf{u}\cdot \textbf{u}+F(\phi,\nabla \phi)+\dfrac{1}{2}\textbf{D}\cdot\textbf{E}\right)d\Omega+\int_{\partial \Omega_{s}}\Theta(\phi) dS. \label{eq_14}
\end{equation}
Here $\Omega$ and $\partial \Omega_{s}$ denote the flow domain and the solid domain boundary, respectively. $F(\phi,\nabla \phi)$ is the free energy density function defined in~\eqref{double_well}. The term $\dfrac{1}{2}\textbf{D}\cdot\textbf{E}$ represents the quasi-static electric energy of the system~\cite{landau2013electrodynamics}. $\Theta(\phi)$ is a wall energy density function, whose form is given later,  aiming to take into account the contact angle effect.

Taking the time derivative of~\eqref{eq_14} and using equations~\eqref{eq_6}--\eqref{double_well} and equation~\eqref{eq_112c} in Appendix A lead to the following energy balance equation,
\begin{equation}\small
\begin{split}
 \frac{dE}{dt}&=-\int_{\Omega}\gamma_1\left|
 \frac12(\rho_1-\rho_2)
 \nabla\left( \dfrac{\partial F}{\partial \phi}-\nabla\cdot\dfrac{\partial F}{\partial\nabla\phi}-\dfrac{\epsilon'}{2}\textbf{E}\cdot \mbs E\right)\right|^2 
 -\int_{\Omega}\mu\left\| \nabla\mbs u\right\|^2 \\
 &\quad -\int_{\partial \Omega}\left[\left(\dfrac{\partial F}{\partial \phi}-\nabla\cdot\dfrac{\partial F}{\partial\nabla\phi}-\dfrac{\epsilon'}{2}\textbf{E}\cdot\mbs E\right)\dfrac{\rho_1-\rho_2}{2}-\dfrac{1}{2}(\textbf{u}\cdot\textbf{u})\right](\tilde{\mbs J}\cdot\mbs n) \\
 &\quad +\int_{\partial \Omega}
 \left[\mu\nabla\textbf{u}\cdot\textbf{n}-F\textbf{n}-\dfrac{1}{2}(\textbf{u}\cdot\textbf{u})\textbf{n}\right]
 \cdot\textbf{u} 
 -\int_{\partial \Omega}
 \lambda\left(\mbs n\cdot\nabla\phi\right)
 \dfrac{\partial \phi}{\partial t}
 +\int_{\partial \Omega_{s}}\Theta'(\phi)\dfrac{\partial \phi}{\partial t} 
 -\int_{\partial \Omega}(\textbf{E}\times\textbf{H})\cdot\textbf{n}, 
 \end{split}
 \label{eq_a15}
\end{equation}
where $\partial\Omega$ denotes the boundary of $\Omega$ and 
$\partial\Omega_s\subset\partial\Omega$ is the solid portion of $\partial\Omega$.
The model ensures the dissipativeness of
the volume integral terms on the right hand side (RHS). Whether the boundary integral terms are dissipative depends on
the imposed boundary conditions, which can guide the choice for the appropriate forms of boundary conditions.
The term $\textbf{E}\times\textbf{H}$ is the  Poynting vector, representing the electromagnetic energy flux to the system~\cite{jackson1999classical}.

\subsection{Equilibrium Solution}
\label{sec:steady}

The incorporation of the electric field term into the  chemical potential and the choice of the $\epsilon(\phi)$ form in~\eqref{interpolation} (see also Remark~\ref{rem_a0}) play an important role in our model. It allows us to derive the energy inequality, thus leading to a thermodynamically consistent model. It also enables us to compute the equilibrium state (steady state) of the two-phase dielectric system by using essentially  the phase field equation only, instead of using the full system coupled with the Navier-Stokes equations. 
We note that in some other studies in the literature (e.g.~\cite{lin2012phase,yang20133d}), where the electric field term is absent from the chemical potential, this benefit does not exist and one needs to solve the full set of governing equations in time in order to find the steady state of the problem.

Specifically, the simpler method for computing the equilibrium solution to the current model is as follows.
At equilibrium ($\frac{\partial}{\partial t}=0$), 
the model represented by the equations~\eqref{eq_6}--\eqref{eq_a8} admits the following solution, 
\begin{equation}\small
\phi_s(\mbs x),\ P_s(\mbs x),\ V_s(\mbs x),\ \text{and}\ \mbs u_s(\mbs x)\equiv 0, 
\end{equation}
where 
\begin{subequations}\label{eq_15}
\begin{align}\small
&
\gamma_1\Delta\left[\lambda h(\phi_s)-\lambda \Delta\phi_s-\dfrac{\epsilon'(\phi_s)}{2}\textbf{E}\cdot\mbs E\right] = 0, \label{eq_15a}  \\
&
\nabla \mathcal{P}_s - \nabla\cdot(\lambda\phi_s\otimes\nabla\phi_s) 
-\frac{\nabla\epsilon}{2}(\mbs E\cdot\mbs E) = 0, \label{eq_15b} \\
&
\nabla\cdot\left(\epsilon\nabla V_s \right)=0,
\label{eq_15c} \\
&
\mbs E = \nabla V_s. \label{eq_15d}
\end{align}
\end{subequations}
It is important to note that the equilibrium state of the current model is compatible with the zero velocity ($\mbs u=0$) requirement based on physics. In the presence of an external electric field, the $(\mbs E\cdot\mbs E)$ term in~\eqref{eq_15a} will cause the equilibrium configuration of the fluid interface to deviate from that of the case with no electric field (e.g.~circular or spherical, due to the surface tension).

These characteristics of the current model suggest that we can employ an alternative system  to compute the steady-state solution to this model. We replace the equation~\eqref{eq_15a} by
the following Cahn-Hilliard equation,
\begin{equation}\label{eq_16}\small
\frac{\partial\phi_s}{\partial \tau}
= \gamma_1\Delta\left(\lambda h(\phi_s)-\lambda \Delta\phi_s-\dfrac{\epsilon'}{2}\textbf{E}\cdot\mbs E\right),
\end{equation}
where $\tau$ is a pseudo-time. 
We solve the system consisting of equations~\eqref{eq_16} and \eqref{eq_15c}--\eqref{eq_15d} by some time marching scheme until the steady state is reached.
This  in principle will produce the equilibrium solution to the original model consisting of~\eqref{eq_6}--\eqref{eq_a8}, with $\mbs u_s= 0$ and $\mathcal{P}_s$ computed by using~\eqref{eq_15b}.
The alternative system consisting of~\eqref{eq_16} and \eqref{eq_15c}--\eqref{eq_15d} is simpler and  faster to 
compute than the original full system consisting of equations~\eqref{eq_6}--\eqref{eq_a8}.
We will demonstrate the effectiveness of this simpler method for computing the equilibrium solution in Section~\ref{sec:tests} using numerical simulations, and also compare the results obtained using the simpler method and using the full model consisting of equations~\eqref{eq_6}--\eqref{eq_a8}.

\subsection{Normalization, Computational Domain, and Boundary/Initial Conditions}

In numerical simulations we employ the normalized non-dimensional form of the governing equations. The normalization discussed here is for the full model~\eqref{eq_6}--\eqref{eq_a8}. We employ a somewhat different normalization for the simpler steady-state model consisting of equations~\eqref{eq_15a}--\eqref{eq_15d}, which will be specified in a later section.
Let $L_0$ denote a characteristic length scale, $V_d$ a characteristic electric potential, $\mu_0$ a characteristic dynamic viscosity, and $\epsilon_0$ the vacuum permittivity with $\epsilon_0=8.85418781\times10^{-12}A^2\cdot s^4/(kg\cdot m^3)$. 
Table~\ref{tab_1} lists the normalization constants for different variables and parameters. For example, the normalized $\lambda$ is given by $\frac{\lambda}{\epsilon_0 V_d^2}$ according to this table.
The normalized governing equations have the same form as the original dimensional ones, and they are also given by the equations~\eqref{eq_6}--\eqref{eq_a8}. Henceforth, the variables and parameters appearing in the equations (and boundary/initial conditions) are understood to have been normalized appropriately, and we will not differentiate their dimensional and non-dimensional forms.

\begin{table}[tb]
\centering\small
\begin{tabular}{cc| cc}
     \hline
     variable & normalization constant & variable & normalization constant \\ \hline
     $x,y,z,\eta,d$ & $L_0$ & $V$ & $V_d$ \\
     $\epsilon,\epsilon_1,\epsilon_2$ & $\epsilon_0$ & $\mu,\mu_1,\mu_2$ & $\mu_0$ \\
    $\mbs u$ &  $u_0=\frac{\epsilon_0V_d^2}{L_0\mu_0}$ &
    $\mathcal{P}$, $P$, $p$ &  $\rho_0u_0^2$\\
    $\rho,\rho_1,\rho_2$ &   $\frac{\mu_0^2}{\epsilon_0V_d^2}$ &
    $\lambda$ & $\epsilon_0 V_d^2 $\\
    $\gamma_1$  & $L_0^2/\mu_0$  & $\mbs{E}$ & $V_d/L_0$ \\
    $\phi,\psi$ & $1$ & $t,\Delta t$ & $\frac{L_0^2\mu_0}{\epsilon_0V_d^2}$\\
     \hline
\end{tabular}
\caption{\small Normalization constants for variables and parameters. Choose $L_0,\epsilon_0,V_d,\mu_0$.}
\label{tab_1}
\end{table}

\begin{figure}[tb]\small
\centerline{
\includegraphics[height=1.5in]{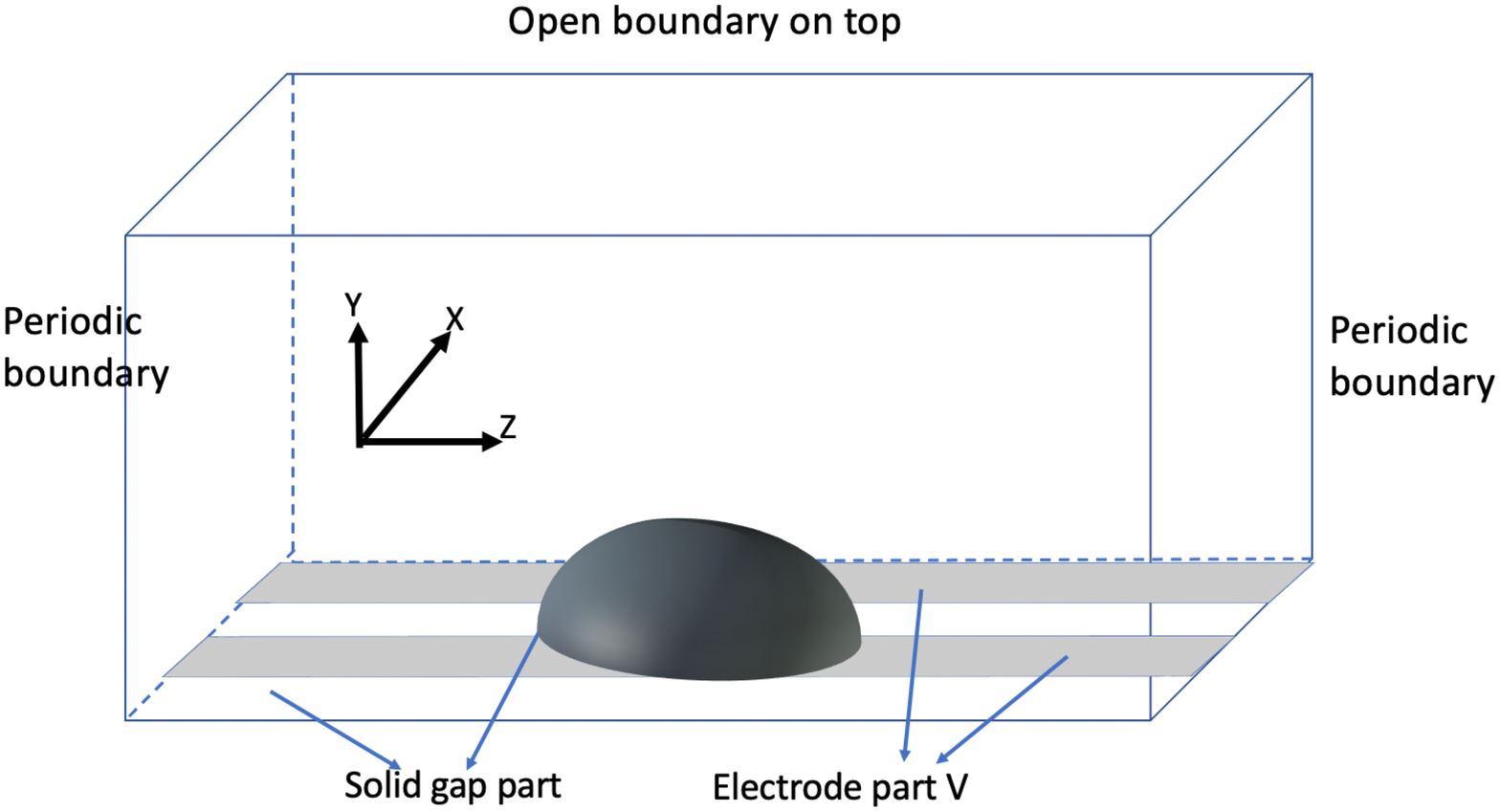}(a)\quad
\includegraphics[height=1.3in]{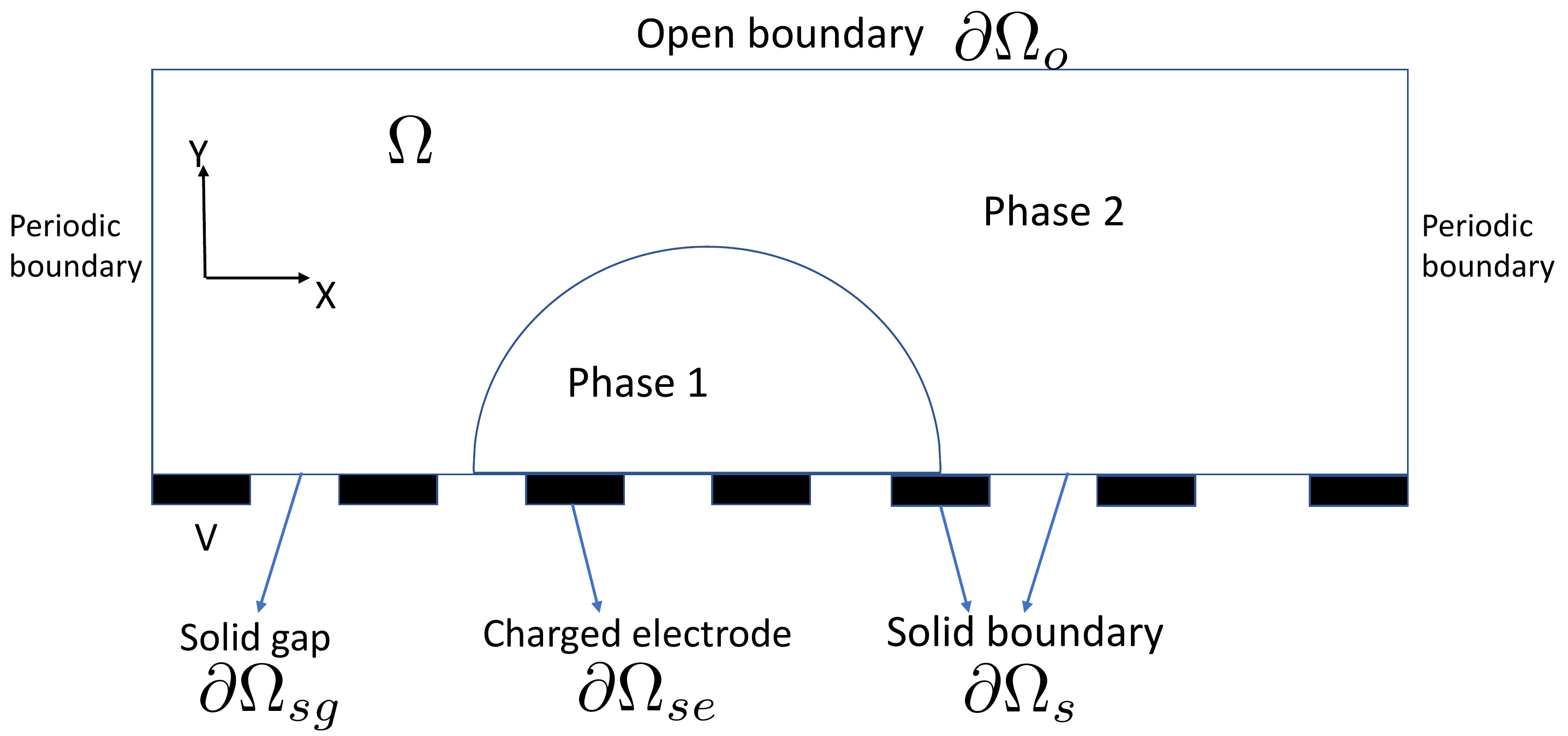}(b)
}
\caption{\small Typical flow domain and settings for (a) 3D and (b) 2D. On the bottom wall, the gray stripes denote the electrodes and the white stripes denote the gap between the electrodes.
} 
\label{general_setup_3d} 
\end{figure}

We perform two-phase dielectric flow simulations in both 2D and 3D. The flow domain and the settings considered here largely follow those of the experiments by McHale, Brown and collaborators~\cite{mchale2011dielectrowetting,brown2009voltage,brown2015dielectrophoresis}.
Especially  we assume that in 3D the domain has at least one homogeneous direction, so that Fourier expansions can be performed along that direction.
Figure~\ref{general_setup_3d} shows typical flow domains and configurations for 2D and 3D. We are interested in the deformation or motion of dielectric droplets on a solid wall. A regular array of parallel electrodes (gray stripes in plot (a), or the dark regions in plot(b)) are embedded on the bottom wall, stretching along the $z$ direction for 3D. The electrodes are separated by wall surfaces (white stripes in plot (a), or white regions in plot (b)). The top of the flow domain is open, while in the horizontal directions the flow is assumed to be periodic.

Let $\partial\Omega_o$ denote the open boundary on the top, $\partial\Omega_{se}$ (``se" standing for solid-electrode) denote the region of electrodes on the bottom wall, and $\partial\Omega_{sg}$ (``sg" standing for solid-gap ) denote the gap region between the electrodes on the wall. The  bottom wall constitutes the solid domain boundary $\partial\Omega_{s}$, with $\partial\Omega_{s}=\partial\Omega_{se}\cup \partial\Omega_{sg}$.

We employ periodic conditions for all the field variables in the horizontal directions and the following boundary conditions for the top and bottom sides of the domain:\\[5pt]
\noindent\underline{Top open boundary $\partial\Omega_o$:}
\begin{equation}\small
\text{(phase field eq.)}\ \left\{
\begin{aligned}
&\textbf{n}\cdot\nabla\left(\lambda h(\phi)-\lambda \Delta \phi-\dfrac{\epsilon'}{2}\textbf{E}\cdot\mbs E\right) = 0, \quad on\quad \partial \Omega_{o}\\
& \nabla \phi \cdot \textbf{n} = 0, \quad on\quad \partial \Omega_{o}.
\end{aligned}
\right.
\label{eq_a19}
\end{equation}
\begin{equation}\small
\text{(momentum eq.)} \ \
\dfrac{\partial \textbf{u}}{\partial n}=0, \quad P=0, \quad on\quad \partial \Omega_{o}.
\label{eq_a24}
\end{equation}
\begin{equation}\small
\text{(electric potential)}\ 
\dfrac{\partial V}{\partial n}=0, \quad on\quad \partial \Omega_{o}.
\label{eq_a25}
\end{equation}
\\
\noindent\underline{Bottom wall ($\partial\Omega_{se}\cup\partial\Omega_{sg}=\partial\Omega_{s}$):}
\begin{equation}\small
\text{(phase field eq.)}\ 
\left\{
\begin{aligned}
&\textbf{n}\cdot\nabla\left(\lambda h(\phi)-\lambda \Delta \phi-\dfrac{\epsilon'}{2}\textbf{E}\cdot\mbs E\right) = 0, \quad on\quad \partial \Omega_{s}
\\
& \lambda\nabla \phi \cdot \textbf{n} + \Theta'(\phi)= 0, \quad on\quad \partial \Omega_{s}
\label{eq_19}
\end{aligned}
\right.
\end{equation}
\begin{equation}\small
\text{(momentum eq.)}\ \
    \textbf{u} = \textbf{0},\quad on\quad \partial \Omega_{s},
    \label{eq_a22}
\end{equation}
\begin{equation}\label{eq_21}\small
\text{(electric potential)}\ 
\left\{
\begin{aligned}
&\dfrac{\partial V}{\partial n}=0, \quad on\quad \partial \Omega_{sg}\\
& V = \mathcal{V}_e, \quad on\quad \partial \Omega_{se}.
\end{aligned}
\right.
\end{equation}

In equation~\eqref{eq_19} $\Theta(\phi)$ denotes the wall energy density function, which accounts for the contact angle effect, given by
\begin{equation}\small
\Theta(\phi) = \gamma \cos(\theta_s)\dfrac{\phi(\phi^2-3)}{4}+\dfrac{1}{2}(\gamma_{s1}+\gamma_{s2}). \label{theta_def}
\end{equation}
where $\gamma,\gamma_{s1},\gamma_{s2}$ are interfacial tension between phase1-phase2, phase1-solid and phase2-solid, and $\theta_s$ is the static contact angle. This functional form is essentially a Hermite interpolation of interfacial tensions; see \cite{dong2012imposing} for more details.
In equation~\eqref{eq_21} $\mathcal{V}_e$ denotes the imposed
voltage on the electrodes. 
We will in general impose an alternate negative/positive voltage on adjacent electrodes as in the experiments  (see e.g.~\cite{mchale2011dielectrowetting}). 
In the gap region between the electrodes, we have employed a simple condition $\frac{\partial V}{\partial n}=0$. This essentially assumes that the electric field at the wall (gap region) has only a tangent component. Note that this condition is exact if the fluids and the wall have matching permitivities or when the fluid is homogeneous~\cite{engan1969excitation}.
In more general cases,
this boundary condition may not be exactly accurate. 
We adopt this boundary condition  because of its simplicity, and that the simulation results indicate that it can capture the flow physics reasonably well.
We note that the set of boundary conditions~\eqref{eq_a19}--\eqref{eq_21}
is reduction-consistent with $\epsilon(\phi)$ given by~\eqref{interpolation}
and $\Theta(\phi)$ given by~\eqref{theta_def}.

Finally we employ the following initial conditions,
\begin{align}\small
&
\mbs u(\mbs x,t=0) = \mbs u_0(\mbs x), \label{eq_a26}\\
&
\phi(\mbs x,t=0) = \phi_0(\mbs x), \label{eq_a27}
\end{align}
where $\mbs u_0$ and $\phi_0$ denote the initial 
distributions for the velocity and the phase field function.

\section{Numerical Algorithm}\label{section3}

\subsection{Algorithm Formulation}

The system consisting of the equations~\eqref{eq_6}--\eqref{eq_a8}, the boundary conditions~\eqref{eq_a19}--\eqref{eq_21} and the
periodic conditions along the horizontal directions,
and the initial conditions~\eqref{eq_a26}--\eqref{eq_a27}
constitute the initial/boundary value problem we need to solve for the velocity, pressure, phase field, and the 
electric potential.

For the purpose of numerical testing, we modify some of the equations and boundary conditions slightly by adding certain prescribed source terms. These source terms are useful for testing the convergence of the method using manufactured solutions, and they will be set to zero in actual flow simulations. Specifically, we re-write equations~\eqref{eq_6}, \eqref{eq_7} and~\eqref{eq_9} into,
\begin{align}\small
&\dfrac{\partial \phi}{\partial t}+\textbf{u}\cdot\nabla\phi=\gamma_1\Delta\left(\lambda h(\phi)-\lambda \Delta\phi-\dfrac{\epsilon'}{2}\textbf{E}\cdot\mbs E\right)
+ g(\mbs x,t),
    \label{eq_a27}\\
    &\dfrac{\partial \mbs u}{\partial t}+ \mbs N(\mbs u) +\frac{1}{\rho}\Tilde{\mbs J}\cdot\nabla\mbs u=-\frac{\lambda}{\rho}\nabla^2\phi\nabla\phi
    -\dfrac{\epsilon'}{2\rho}(\textbf{E}\cdot\mbs E)\nabla\phi
+\frac{\mu}{\rho}\nabla^2\mbs u
+\frac{1}{\rho}\nabla\mu\cdot\mathcal{D}(\mbs u)
-\frac{1}{\rho}\nabla P
    +\frac{1}{\rho}\mbs f(\mbs x,t)
    \label{eq_a28} \\
    &
    \nabla\cdot(\epsilon V) = f_V(\mbs x,t), \label{eq_28}
\end{align}
where $g$, $\mbs f$ and $f_V$ are prescribed source terms, and
\begin{equation}\small
P = \mathcal{P} + \frac{\lambda}{2}\nabla\phi\cdot\nabla\phi, \quad
\mathcal{D}(\mbs u) = \nabla\mbs u + \nabla\mbs u^T, \quad
\mbs N(\mbs u) = \mbs u\cdot\nabla\mbs u.
\end{equation}
The boundary conditions~\eqref{eq_a19}--\eqref{eq_a24} are modified as,
\begin{align}\small
    &
    \mbs n\cdot\nabla\left(\lambda h(\phi)-\lambda\nabla^2\phi -\frac12\epsilon'\mbs E^2 \right) = g_1(x,t), \quad
    \mbs n\cdot\nabla\phi = g_2(\mbs x,t), \quad
    \mbs x\in\partial\Omega_{o}; \label{eq_29}
    \\
    &
    \frac{\partial\mbs u}{\partial n}=\mbs f_1(\mbs x,t), \quad
    P = f_2(\mbs x,t), \quad 
    \mbs x\in\partial\Omega_{o}; \label{eq_30}
\end{align}
where $\mbs f_1$, $f_2$, $g_1$, $g_2$ are 
prescribed source terms. The boundary conditions~\eqref{eq_19}--\eqref{eq_a22} are modified as,
\begin{align}\small
    &
    \mbs n\cdot\nabla\left(\lambda h(\phi)-\lambda\nabla^2\phi -\frac12\epsilon'\mbs E^2 \right) = g_1(x,t), \quad
    \mbs n\cdot\nabla\phi + \frac{1}{\lambda}\Theta'(\phi)= g_3(\mbs x,t), \quad \mbs x\in\partial\Omega_{s}; \label{eq_32}
    \\ &
    \mbs u=\mbs w(\mbs x,t), \quad \mbs x\in\partial\Omega_{s}; \label{eq_33}
\end{align}
where $g_3$ and $\mbs w$ are  prescribed source terms.

We next present an algorithm for solving the system
consisting of equations~\eqref{eq_a27}--\eqref{eq_28}, \eqref{eq_8}, \eqref{eq_a8}, \eqref{eq_29}--\eqref{eq_33}, \eqref{eq_a25}, \eqref{eq_21}, together with the periodic conditions in the horizontal directions.
Let $n\geqslant 0$ denote the time step index, $\Delta t$ denote the time step size, and
$(\cdot)^n$ denote the variable $(\cdot)$ at time
step $n$.
Given $(\mbs u^n,P^n,\phi^n,V^n)$, we compute these quantities at step $(n+1)$ successively by the following procedure: \\[6pt]
\noindent\underline{Electric potential $V^{n+1}$ and electric field $\mbs E^{n+1}$:}
\begin{subequations}\label{eq_a35}
\begin{align}\small
&\nabla \cdot(\varepsilon_0\nabla V^{n+1})=f_V^{n+1}- \nabla \cdot\left[(\epsilon(\phi^{*,n+1})-\varepsilon_0)\nabla V^{*,n+1}\right],
\label{ep}\\
&\frac{\partial V^{n+1}}{\partial n}=0, \quad \text{on} \ \partial\Omega_{o} \cup \partial\Omega_{sg}, 
\label{epa}\\
&V^{n+1}=\mathcal{V}_e, \quad \text{on} \ \partial\Omega_{se}, 
\label{epbb}
\\
& \mbs E^{n+1} = \nabla V^{n+1}.
\label{eq_E}
\end{align}
\end{subequations}

\noindent\underline{Phase field $\phi^{n+1}$:}
\begin{subequations}\label{eq_a36}
\begin{align}\small
   &\frac{\gamma_0\phi^{n+1}-\hat{\phi}}{\Delta t}+\textbf{u}^{*,n+1}\cdot \nabla \phi^{*,n+1}
   =-\lambda \gamma_1 \nabla^2\left[\nabla^2 \phi^{n+1}-\frac{S}{\eta^2}(\phi^{n+1}-\phi^{*,n+1})\right. \notag\\
   &\qquad \qquad \qquad \qquad \qquad \qquad \qquad
   \left.-h(\phi^{*,n+1})+\frac{\epsilon'(\phi^{*,n+1})}{2\lambda} \left|\mbs E^{n+1}\right|^2\right]+g^{n+1}, 
\label{phi}\\
&\textbf{n}\cdot \nabla\left[\nabla^2\phi^{n+1}-\frac{S}{\eta^2}(\phi^{n+1}-\phi^{*,n+1})-h(\phi^{*,n+1})+\frac{\epsilon'(\phi^{*,n+1})}{2\lambda} \left|\mbs E^{n+1}\right|^2 \right]=g^{n+1}_1, \ \text{on} \ \partial\Omega_{o}\cup \partial\Omega_{s}, 
\label{phia}\\
&\textbf{n} \cdot \nabla \phi^{n+1}=g^{n+1}_2, \  \text{on} \ \partial\Omega_{o}, 
\label{phib}\\
&-\mbs n\cdot\nabla \phi^{n+1}-\frac{\Theta'(\phi^{*,n+1})}{\lambda}=g_3^{n+1},\  \text{on} \ \partial\Omega_{s}. 
\label{phic}
\end{align}
\end{subequations}

\noindent\underline{Pressure $P^{n+1}$:}
\begin{subequations}\label{eq_a37}
\begin{align}\small
&\frac{\gamma_0\Tilde{\textbf{u}}^{n+1}-\hat{\textbf{u}}}{\Delta t}+\dfrac{1}{\rho_0}\nabla P^{n+1}=-\textbf{N}(\textbf{u}^{*,n+1})+\left(\dfrac{1}{\rho_0}-\dfrac{1}{\rho^{n+1}}\right)\nabla P^{*,n+1}-\dfrac{\mu^{n+1}}{\rho^{n+1}}\nabla\times\nabla\times\textbf{u}^{*,n+1}\notag\\
   &\qquad \qquad \qquad \qquad \qquad \qquad +\dfrac{1}{\rho^{n+1}}\nabla\mu^{n+1}\cdot\mathcal{D}(\textbf{u}^{*,n+1})-\dfrac{\lambda}{\rho^{n+1}}\nabla^2\phi^{n+1}\nabla\phi^{n+1}+\dfrac{f^{n+1}}{\rho^{n+1}}\notag\\
   &\qquad \qquad \qquad \qquad \qquad \qquad- \dfrac{1}{\rho^{n+1}}\Tilde{\textbf{J}}^{n+1}\cdot\nabla\textbf{u}^{*,n+1}-\dfrac{\epsilon'(\phi^{n+1})}{2\rho^{n+1}} \left|\mbs E^{n+1}\right|^2 \nabla\phi^{n+1}, \label{u} \\
   &\nabla\cdot\Tilde{\textbf{u}}^{n+1}=0, 
   \label{ua}\\
   &\frac{\partial \Tilde{\textbf{u}}^{n+1}}{\partial n}=\mbs f_1^{n+1},\ \text{on} \ \partial\Omega_{o}, 
   \label{ub}\\
   &P^{n+1}=f_2^{n+1},\ \text{on} \ \partial\Omega_{o},
   \label{uc}\\
   &\Tilde{\textbf{u}}^{n+1}\cdot\textbf{n}=\textbf{w}^{n+1}\cdot\textbf{n},\ \text{on} \ \partial\Omega_{s}. 
   \label{ud}
\end{align}
\end{subequations}

\noindent\underline{Velocity $\textbf{u}^{n+1}$:}
\begin{subequations}\label{eq_a38}
\begin{align}\small
   &\frac{\gamma_0\textbf{u}^{n+1}-\hat{\textbf{u}}}{\Delta t}+\dfrac{1}{\rho_0}\nabla P^{n+1}-\nu_m\nabla^2\textbf{u}^{n+1}
   =-\textbf{N}(\textbf{u}^{*,n+1})+\left(\dfrac{1}{\rho_0}-\dfrac{1}{\rho^{n+1}}\right)\nabla P^{*,n+1}
   \notag \\
   &  \qquad \qquad 
   +\left(\nu_m-\dfrac{\mu^{n+1}}{\rho^{n+1}}\right)\nabla\times\nabla\times\textbf{u}^{*,n+1}
   +\dfrac{1}{\rho^{n+1}}\nabla\mu^{n+1}\cdot\mathcal{D}(\textbf{u}^{*,n+1})-\dfrac{\lambda}{\rho^{n+1}}\nabla^2\phi^{n+1}\nabla\phi^{n+1} \notag \\
   & \qquad \qquad
   +\dfrac{f^{n+1}}{\rho^{n+1}}
   - \dfrac{1}{\rho^{n+1}}\Tilde{\textbf{J}}^{n+1}\cdot\nabla\textbf{u}^{*,n+1}-\dfrac{\epsilon'(\phi^{n+1})}{2\rho^{n+1}} \left|\mbs E^{n+1}\right|^2 \nabla\phi^{n+1} \label{v} \\
   &\frac{\partial \textbf{u}^{n+1}}{\partial n}=\mbs f_1^{n+1},\ \text{on} \ \partial\Omega_{o},
   \label{va}\\
&\textbf{u}^{n+1}=\textbf{w}^{n+1},\ \text{on} \ \partial\Omega_{s}.
   \label{vb}
\end{align}
\end{subequations}
In the horizontal directions ($x$ in 2D, $x$ and $z$ in 3D) we impose periodic conditions for ($V^{n+1},\phi^{n+1},P^{n+1},\mbs u^{n+1}$). These periodic conditions are not explicitly included in
the above system of equations.

The meanings of those symbols involved in the above equations are as follows.
$\Tilde{\textbf{J}}^{n+1}$ in \eqref{u} and \eqref{v} is given by,
\begin{equation}\small
\Tilde{\textbf{J}}^{n+1}=\dfrac{1}{2}(\rho_1-\rho_2)\lambda\gamma_1\nabla\left[\ \nabla^2\phi^{n+1}-\dfrac{S}{\eta^2}(\phi^{n+1}-\phi^{*,n+1})-h(\phi^{*,n+1})+\dfrac{\epsilon'(\phi^{*,n+1})}{2\lambda} \left|\mbs E^{n+1}\right|^2 \right],
\end{equation}
in which $S$ is a stabilization parameter to be described later.
$\Tilde{\textbf{u}}^{n+1}$ is an auxiliary approximation of $\mbs u$ at time step $(n+1)$.  
Let $\xi$ denote a generic variable. Then in the above equations the expression $\dfrac{1}{\Delta t}D\xi^{n+1}=\dfrac{1}{\Delta t}(\gamma_0\xi^{n+1}-\hat{\xi})$ represents an approximation of $\left.\dfrac{\partial \xi}{\partial t}\right|^{n+1}$ by the $J$-th order backward differentiation formula (BDF), where $J=1$ or $2$, and
$\xi^{*,n+1}$ represents a $J$-th order explicit approximation of $\xi^{n+1}$. They are explicitly given by,
\begin{equation}\small
D\xi^{n+1}=
\left\{
\begin{array}{ll}
\xi^{n+1}-\xi^n, & if \quad J=1,\\
\frac32\xi^{n+1}-2\xi^n+\frac12\xi^{n-1}, & if \quad J=2;\\
\end{array}
\right.
\quad
\xi^{*,n+1}=
\left\{
\begin{array}{ll}
\xi^n, & if \quad J=1,\\
2\xi^n-\xi^{n-1}, & if \quad J=2.\\
\end{array}
\right.
\end{equation} 
Note that $\gamma_0=1$ for $J=1$, and $\frac32$ for $J=2$.

The $\varepsilon_0$ in equation \eqref{ep},
$\rho_0$ in~\eqref{u}, and $\nu_m$ in~\eqref{v} are constant algorithmic parameters. With those terms involving these constants in the above formulation, the algorithm leads to linear algebraic systems with constant and time-dependent coefficient matrices upon discretization, which makes the current method computationally highly efficient.
We choose $\varepsilon_0=\max(\epsilon_1,\epsilon_2)$, where $\epsilon_1,\epsilon_2$ are the permittivities of two dielectric fluids. 
In~\eqref{ep} we have treated the $\varepsilon_0\nabla V$ term on the left hand side (LHS) implicitly and the correction term $(\epsilon-\varepsilon_0)\nabla V$ in the RHS explicitly. The consistent approximation of these terms ensures the $J$-th order accuracy of the overall scheme.
We choose $\rho_0$ and $\nu_m$ following~\cite{dong2012time,dong2014outflow}, specifically with 
$\rho_0=\min(\rho_1,\rho_2)$ and $\nu_m\geq \frac{1}{2}\frac{\max(\mu_1,\mu_2)}{\min(\rho_1,\rho_2)}$, where $\rho_1$ and $\rho_2$ are the densities, and
$\mu_1$ and $\mu_2$ are the dynamic viscosities of two dielectric fluids, respectively.

The term $\frac{S}{\eta^2}(\phi^{n+1}-\phi^{*,n+1})$ in equation \eqref{phi} is a  stabilization term, 
where $S$ is a chosen constant satisfying $S\geq \eta^2 \sqrt{\frac{4\gamma_0}{\lambda \gamma_1\Delta t}}$. This allows us to reformulate the fourth order Cahn-Hilliard equation into two decoupled Helmholtz type equations (see \cite{dong2012imposing} for more details).
\begin{remark}
    In the above algorithm 
    we would need to compute the  initial distribution of the electric potential (and the electric field) to start the computation, i.e.~solving for $V(\mbs x,t=0)=V_0(\mbs x)$.
    We  use a fixed point iteration to compute the initial distribution,
    \begin{equation}\small
        \nabla \cdot\left(\varepsilon_0\nabla V^{(k+1)}\right)=f_V^{ini}- \nabla \cdot\left[(\epsilon(\phi_0)-\varepsilon_0)\nabla V^{(k)}\right],
        \label{eq_41}
    \end{equation}
    where $\phi_0$ is the initial phase field distribution given in~\eqref{eq_a27},  $f_V^{ini}=f_V(\mbs x,t=0)$,
    and the superscript in $V^{(k)}$ refers to the iteration index. 
    The initial distribution of the electric potential is obtained upon convergence of this iteration.
    
\end{remark}

\subsection{Implementation and Spatial Discretization}

We next discuss how to implement the algorithm represented by equations~\eqref{ep}--\eqref{vb} using high-order $C^0$ spectral elements in 2D and a hybrid Fourier spectral/spectral-element method in 3D. We first derive a weak form of the algorithm, which is suitable for  both 2D and 3D. Then we further  transform the weak form in 3D to a form specifically for the hybrid Fourier spectral/spectral-element discretization.

Given $(\textbf{u}^n,P^n,\phi^n,V^n)$, we wish to compute $(\textbf{u}^{n+1},P^{n+1},\phi^{n+1},V^{n+1})$.
We will first derive the weak forms about these variables in the continuous space by taking the $L^2$ inner product between an arbitrary test function and the equations about these variables. Then we  restrict these variables and the test functions to appropriate function spaces to attain the final weak forms.

Let $e(\mbs x)$ denote an arbitrary test function. Taking the $L^2$ inner product between $e$ and~\eqref{ep} and using the integration by parts, we attain the weak form for $V^{n+1}$,
\begin{equation}\small
    \int_{\Omega}\varepsilon_0 \nabla V^{n+1}\cdot\nabla e=-\int_{\Omega} (\epsilon(\phi^{*,n+1})-\varepsilon_0)\nabla V^{*,n+1}\cdot\nabla e-\int_{\Omega} f_V^{n+1}e,\quad \forall e. \label{eq_41}
\end{equation}
The weak form of the electric field
$\mbs E^{n+1}(\mbs x)=(E_x^{n+1}(\mbs x),E_y^{n+1}(\mbs x),E_z^{n+1}(\mbs x))$ is obtained by taking the $L^2$ inner product between $e(\mbs x)$ and equation~\eqref{eq_E},
\begin{equation}\small
   \int_{\Omega} \mbs E^{n+1} e=\int_{\Omega} \nabla V^{n+1} e, \quad\forall e. 
   \label{eq_43}
\end{equation}

Equation~\eqref{phi} can be written as (see~\cite{dong2012time} for details),
\begin{subequations}
\begin{align}\small
    &
    \nabla^2\psi^{n+1} - \left(\alpha+\dfrac{S}{\eta^2}\right)\psi^{n+1} = Q_1+\nabla^2 Q_2,
    \label{eq_44a} \\
    &
    \nabla^2\phi^{n+1}+ \alpha\phi^{n+1} = \psi^{n+1}, \label{eq_44b}
\end{align}
\end{subequations}
where $\alpha=-\frac{S}{2\eta^2}\left(1+\sqrt{1-\frac{4\gamma_0}{\lambda\gamma_1\Delta t}\frac{\eta^4}{S^2}}\right)$, $\psi^{n+1}$ is an auxiliary variable defined by~\eqref{eq_44b}, and
\begin{equation}\small
Q_1=\dfrac{1}{\lambda\gamma_1}\left(g^{n+1}-\textbf{u}^{*,n+1}\cdot \nabla \phi^{*,n+1}+\dfrac{\hat{\phi}}{\Delta t}\right),
\quad 
Q_2=h(\phi^{*,n+1})-\dfrac{S}{\eta^2}\phi^{*,n+1}-\dfrac{\epsilon'(\phi^{*,n+1})}{2\lambda}|\mbs E^{n+1}|^2.
\end{equation}
Let $\omega(\mbs x)$ denote an arbitrary test function. 
The weak forms for~\eqref{eq_44a} and~\eqref{eq_44b}
are attained by taking the $L^2$ inner product between
$\omega(\mbs x)$ and these equations,
\begin{align}\small
    \int_{\Omega} \nabla \psi^{n+1} &\cdot \nabla \omega 
    + \left(\alpha+\dfrac{S}{\eta^2}\right)\int_{\Omega} \psi^{n+1} \cdot \omega = \int_{\partial\Omega_{o}}\left[g_1^{n+1}+\left(\alpha+\dfrac{S}{\eta^2}\right)g_2^{n+1}\right]\omega \notag\\
    +&\int_{\partial\Omega_{s}} \left[g_1^{n+1}+\left(\alpha+\dfrac{S}{\eta^2}\right)\left(-g_3^{n+1}-\dfrac{\Theta'(\phi^{*,n+1})}{\lambda}\right)\right]\omega
    -\int_{\Omega} Q_1 \omega+\int_{\Omega} \nabla Q_2\cdot \nabla \omega, \quad \forall \omega; \label{eq_46}
\end{align}
\begin{align}\small
  &\int_{\Omega} \nabla \phi^{n+1}\cdot \nabla \omega -\alpha \int_{\Omega} \phi^{n+1}\omega= \int_{\partial \Omega_{o}} g_2^{n+1}\omega+\int_{\partial\Omega_{s}}\left[-g_3^{n+1}-\dfrac{\Theta'(\phi^{*,n+1})}{\lambda}\right]\omega 
  -\int_{\Omega}\psi^{n+1} \omega, \quad \forall \omega. \label{eq_47}
\end{align}


Let $q(\mbs x)$ denote an arbitrary test function that vanishes on  $\partial \Omega_{o}$. 
Taking the $L^2$ inner product between $\nabla q$ and
equation~\eqref{u} leads to the weak form about $P^{n+1}$,
\begin{align}\small
  &\int_{\Omega}\nabla P^{n+1}\cdot\nabla q=\rho_0\int_{\Omega} \left[\textbf{T}+\nabla\left(\dfrac{\mu^{n+1}}{\rho^{n+1}}\right)\times\boldsymbol{\omega}^{*,n+1}\right]\cdot\nabla q\notag\\
  &\qquad \qquad \qquad \qquad-\rho_0\int_{\partial \Omega}\dfrac{\mu^{n+1}}{\rho^{n+1}}\textbf{n}\times\boldsymbol{\omega}^{*,n+1}\cdot\nabla q-\dfrac{\rho_0\gamma_0}{\Delta t}\int_{\partial \Omega_s} \textbf{n}\cdot \mbs w^{n+1}q,\quad \forall q
  \ \text{with}\ q|_{\partial\Omega_o}=0,
  \label{eq_48}
\end{align}        
where $\bm\omega = \nabla\times\mbs u$ and
\begin{align}\small
\textbf{T}&=\dfrac{1}{\rho^{n+1}}\left[f^{n+1}-\lambda(\psi^{n+1}-\alpha\phi^{n+1})\nabla\phi^{n+1}-\dfrac{\epsilon'(\phi^{n+1})}{2}\left|\textbf{E}\right|^2\nabla\phi^{n+1}+\nabla\mu^{n+1}\cdot\mathcal{D}(\textbf{u}^{*,n+1})\right. \notag\\
&\qquad \qquad \left.-\Tilde{\textbf{J}}^{n+1}\cdot\nabla\textbf{u}^{*,n+1}\right]
+\dfrac{\hat{u}}{\Delta t}-\textbf{N}(\textbf{u}^{*,n+1})+\left(\dfrac{1}{\rho_0}-\dfrac{1}{\rho^{n+1}}\right)\nabla P^{*,n+1}.
\label{eq_49}
\end{align}  
When deriving the above equation, the following identity has been used,
\begin{equation*}\small
\begin{split}
    &\dfrac{\mu}{\rho}\nabla \times \boldsymbol{\omega}\cdot\nabla q=\nabla \cdot\left(\dfrac{\mu}{\rho}\boldsymbol{\omega}\times\nabla q\right)-\nabla\left(\dfrac{\mu}{\rho}\right)\times\boldsymbol{\omega}\cdot\nabla q.
\end{split}
\end{equation*}

For the weak form of equation~\eqref{v}, let $v(\mbs x)$ be an arbitrary test function that vanishes on $\partial \Omega_s$,
and we take the $L^2$ inner product between $v(\mbs x)$ and equation(\ref{v}) to get
\begin{align}\small
&\int_{\Omega}\nabla\textbf{u}^{n+1}\cdot \nabla v+\dfrac{\gamma_0}{\nu_m\Delta t}\int_{\Omega}\textbf{u}^{n+1} v=\dfrac{1}{\nu_m}\int_{\Omega}\left(\textbf{T}-\dfrac{1}{\rho_0}\nabla P^{n+1}\right)v-\dfrac{1}{\nu_m}\int_{\Omega}\left(\dfrac{\mu^{n+1}}{\rho^{n+1}}-\nu_m\right)\boldsymbol{\omega}^{*,n+1}\times\nabla v\notag\\
    &\qquad 
    +\dfrac{1}{\nu_m}\int_{\Omega}\nabla\left(\dfrac{\mu^{n+1}}{\rho^{n+1}}-\nu_m\right)\times\boldsymbol{\omega}^{*,n+1}v+\int_{\partial \Omega_{o}}\mbs f_1^{n+1}v
    -\dfrac{1}{\nu_m}\int_{\partial \Omega_{o}}\left(\dfrac{\mu^{n+1}}{\rho^{n+1}}-\nu_m\right)\textbf{n}\times\boldsymbol{\omega}^{*,n+1}v, \notag \\
    & \qquad\quad
    \forall v \ \text{with}\ 
    v|_{\partial\Omega_s}=0,
    \label{eq_52}
\end{align}
where the following identity has been used,
\begin{equation*}\small
v\left(\nu_m-\dfrac{\mu}{\rho}\right)\nabla\times\boldsymbol{\omega}  =\nabla\times\left[v\boldsymbol{\omega}\left(\nu_m-\dfrac{\mu}{\rho}\right)\right]-v\nabla\left(\nu_m-\dfrac{\mu}{\rho}\right)\times\boldsymbol{\omega}-\left(\nu_m-\dfrac{\mu}{\rho}\right)\nabla v\times\boldsymbol{\omega}.
\end{equation*}


\subsubsection{Two Dimensions}

For two-dimensional (2D, $\Omega\subset\mathbb{R}^2$) problems we employ  $C^0$ spectral elements for spatial discretizations. We partition the domain $\Omega$ using a spectral element mesh. Let $\Omega_h$ denote the discretized domain, $\Omega_h=\cup_{e=1}^{N_{e}}\Omega_{h}^e$, where $\Omega_h^e$ ($1\leqslant e\leqslant N_e$) denotes the element $e$ and $N_e$ is
the number of elements in the mesh.
Let $\partial\Omega_h$, $\partial\Omega_{oh}$, $\partial\Omega_{sh}$ denote the discretized versions of the domain boundary $\partial\Omega$, open boundary $\partial\Omega_o$, and solid boundary $\partial\Omega_s$.
Then $\partial\Omega_h = \partial\Omega_{oh}\cup\partial\Omega_{sh}
=\partial\Omega_{oh}\cup\partial\Omega_{seh}\cup\partial\Omega_{sgh}$,
where $\partial\Omega_{seh}$ and $\partial\Omega_{sgh}$ are the discretized versions of the solid-electrode boundary and the solid-gap boundary, respectively.
Let $\Pi_{K}(\Omega_h^e)$ denote the linear space of polynomials defined on $\Omega_h^e$ with their degrees characterized by $K$ ($K$ will be referred to as the element order hereafter).
Define
\begin{equation}\small
\left\{
\begin{split}
&
X_h = \{\ v\in H^1(\Omega_h)\ :\ v|_{\Omega_h^e}\in\Pi_K(\Omega_h^e),
\ 1\leqslant e\leqslant N_e \ \}, \\
&
X_{h0}^E = \{\  v\in X_h\ :\ 
v|_{\partial\Omega_{seh}}=0
\ \}, \\
&
X_{h0}^{P} = \{\ v\in X_h\ :\ 
v|_{\partial\Omega_{oh}}=0
\ \}, \\
&
X_{h0}^u = \{\  v\in X_h\ :\ 
v|_{\partial\Omega_{sh}}=0
\ \}.
\end{split}
\right.
\end{equation}
In what follows we use $(\cdot)_h$
to denote the discretized version of $(\cdot)$.

The 2D fully discretized equations consists of the following:\\[5pt]
\noindent\underline{For $V_h^{n+1}$:} 
find $V_h^{n+1}\in X_h$ such that
\begin{subequations}\label{eq_54}
\begin{align}\small
&
    \int_{\Omega_h}\varepsilon_0 \nabla V_h^{n+1}\cdot\nabla e_h=-\int_{\Omega_h} (\epsilon(\phi_h^{*,n+1})-\varepsilon_0)\nabla V_h^{*,n+1}\cdot\nabla e_h -\int_{\Omega_h} f_{Vh}^{n+1}e_h,\quad \forall e_h\in X_{h0}^E; \label{eq_54a}
    \\
    &
    V_h^{n+1} = \mathcal V_e, 
    \quad \text{on}\ \partial\Omega_{seh}.
    \label{eq_54b}
\end{align}
\end{subequations}
\noindent\underline{For $\mbs E_h^{n+1}$:}
find $\mbs E_h^{n+1}\in [X_h]^2$ such that
\begin{equation}\small
   \int_{\Omega_h} \mbs E_h^{n+1} e_h=\int_{\Omega_h} \nabla V_h^{n+1} e_h, \quad\forall e_h\in X_h. 
   \label{eq_55}
\end{equation}
\noindent\underline{For $\psi_h^{n+1}$:} find $\psi_h^{n+1}\in X_h$ such that
\begin{align}\small
    \int_{\Omega_h} &\nabla \psi_h^{n+1} \cdot \nabla \omega_h 
    + \left(\alpha+\dfrac{S}{\eta^2}\right)\int_{\Omega_h} \psi_h^{n+1} \cdot \omega_h = \int_{\partial\Omega_{oh}}\left[g_{1h}^{n+1}+\left(\alpha+\dfrac{S}{\eta^2}\right)g_{2h}^{n+1}\right]\omega_h \notag\\ 
    &+\int_{\partial\Omega_{sh}} \left[g_{1h}^{n+1}+\left(\alpha+\dfrac{S}{\eta^2}\right)\left(-g_{3h}^{n+1}-\dfrac{\Theta'(\phi_h^{*,n+1})}{\lambda}\right)\right]\omega_h
    -\int_{\Omega_h} Q_{1h} \omega_h+\int_{\Omega_h} \nabla Q_{2h}\cdot \nabla \omega_h, \notag \\
    & \ \forall \omega_h\in X_h. \label{eq_56}
\end{align}
\noindent\underline{For $\phi_h^{n+1}$:} find $\phi_h^{n+1}\in X_h$ such that
\begin{align}\small
  \int_{\Omega_h} \nabla \phi_h^{n+1}\cdot \nabla \omega_h -\alpha \int_{\Omega_h} \phi_h^{n+1}\omega_h=& \int_{\partial \Omega_{oh}} g_{2h}^{n+1}\omega_h+\int_{\partial\Omega_{sh}}\left[-g_{3h}^{n+1}-\dfrac{\Theta'(\phi_h^{*,n+1})}{\lambda}\right]\omega_h\notag\\
  &-\int_{\Omega_h}\psi_h^{n+1} \omega_h, \quad \forall \omega_h\in X_h. \label{eq_57}
\end{align}
\noindent\underline{For $P_h^{n+1}$:} find $P_h^{n+1}\in X_h$ such that
\begin{subequations}
\begin{align}\small
  &\int_{\Omega_h}\nabla P_h^{n+1}\cdot\nabla q_h=\rho_0\int_{\Omega_h} \left[\textbf{T}_h+\nabla\left(\dfrac{\mu_h^{n+1}}{\rho_h^{n+1}}\right)\times\boldsymbol{\omega}_h^{*,n+1}\right]\cdot\nabla q_h\notag\\
  &\qquad \qquad \qquad \qquad-\rho_0\int_{\partial \Omega_h}\dfrac{\mu_h^{n+1}}{\rho_h^{n+1}}\textbf{n}\times\boldsymbol{\omega}_h^{*,n+1}\cdot\nabla q_h-\dfrac{\rho_0\gamma_0}{\Delta t}\int_{\partial \Omega_{sh}} \textbf{n}\cdot \mbs w_h^{n+1}q_h,\quad \forall q_h\in X_{h0}^P.
  \label{eq_58a} \\
  &
  P_h^{n+1}=f_{2h}^{n+1}, \ \text{on} \ \partial\Omega_{oh}. \label{eq_58b}
\end{align}\end{subequations}
\noindent\underline{For $\mbs u_h^{n+1}$:} find $\mbs u_h^{n+1}\in [X_h]^2$ such that
\begin{subequations}
\begin{align}\small
&\int_{\Omega_h}\nabla v_h\cdot\nabla\textbf{u}_h^{n+1} +\dfrac{\gamma_0}{\nu_m\Delta t}\int_{\Omega_h}\textbf{u}_h^{n+1} v_h
    =\dfrac{1}{\nu_m}\int_{\Omega_h}\left(\textbf{T}_h-\dfrac{1}{\rho_0}\nabla P_h^{n+1}\right)v_h \notag \\
    &\qquad
    -\dfrac{1}{\nu_m}\int_{\Omega_h}\left(\dfrac{\mu_h^{n+1}}{\rho_h^{n+1}}-\nu_m\right)\boldsymbol{\omega}_h^{*,n+1}\times\nabla v_h
     +\dfrac{1}{\nu_m}\int_{\Omega_h}\nabla\left(\dfrac{\mu_h^{n+1}}{\rho_h^{n+1}}\right)\times\boldsymbol{\omega}_h^{*,n+1}v_h+\int_{\partial \Omega_{oh}}\mbs f_{1h}^{n+1}v_h\notag\\
    &\qquad-\dfrac{1}{\nu_m}\int_{\partial \Omega_{oh}}\left(\dfrac{\mu_h^{n+1}}{\rho_h^{n+1}}-\nu_m\right)\textbf{n}\times\boldsymbol{\omega}_h^{*,n+1}v_h, \quad \forall v_h\in X_{h0}^u; \label{eq_59a} \\
    &
\textbf{u}_h^{n+1}=\textbf{w}_h^{n+1},\ 
    \text{on}\ \partial\Omega_{sh}.
    \label{eq_59b}
\end{align}\end{subequations}

Therefore, given $(\mbs u^n, P^n, \phi^n, V^n)$, one can compute $V^{n+1}$, $\mbs E^{n+1}$, $\psi^{n+1}$, $\phi^{n+1}$, $P^{n+1}$ and $\mbs u^{n+1}$ by solving equations~\eqref{eq_54a}--\eqref{eq_59b} successively in an uncoupled fashion. The solution procedure is summarized in Algorithm~\ref{alg_1}.

\begin{algorithm}[tb]\small
  \DontPrintSemicolon
  \SetKwInOut{Input}{input}\SetKwInOut{Output}{output}

  \Input{$V^n$, $\phi^n$, $P^n$, $\mbs u^n$.}
  \Output{ $V^{n+1}$, $\mbs E^{n+1}$, $\psi^{n+1}$, $\phi^{n+1}$, $P^{n+1}$, $\mbs u^{n+1}$. }
  \BlankLine  \BlankLine
  solve equations~\eqref{eq_54a}--\eqref{eq_54b}
  for $V^{n+1}$\;
  solve equation~\eqref{eq_55} for $\mbs E^{n+1}$\;
  solve equation~\eqref{eq_56} for $\psi^{n+1}$\;
  solve equation~\eqref{eq_57} for $\phi^{n+1}$\;
  solve equations~\eqref{eq_58a}--\eqref{eq_58b}
  for $P^{n+1}$\;
  solve equations~\eqref{eq_59a}--\eqref{eq_59b}
  for $\mbs u^{n+1}$\;
  
  \caption{Solution Procedure within a Time Step for 2D Dielectric Flows}
  \label{alg_1}
\end{algorithm}

\subsubsection{Three Dimensions}

For three dimensions (3D, $\Omega\subset\mathbb{R}^3$) we concentrate on problems with one homogeneous direction in this work, so that Fourier expansions can be employed along that direction,
as stated previously. Let us assume that the homogeneous direction is along the $z$ axis, and we employ a hybrid spectral-element/Fourier spectral discretization to solve the problem, with spectral element discretization in the $xy$ plane and Fourier spectral discretization along the $z$ direction.

We take the  domain along the $z$ direction as $z\in [0,L_z]$, and assume that the domain and all the dynamic variables are periodic at $z=0$ and $z=L_z$,
where $L_z$ is the dimension of the computational domain in $z$.
Then the following relations hold,
\begin{equation}\small
\left\{
\begin{split}
&
\Omega=\Omega_{2D}\otimes[0,L_z],\quad \partial \Omega=\partial \Omega_{2D}\otimes[0,L_z], \quad
    \partial\Omega_s = \partial\Omega_{s}^{2D}\otimes[0,L_z],
    \quad
    \partial\Omega_o = \partial\Omega_{o}^{2D}\otimes[0,L_z],
    \\
    &
    \partial\Omega_{se} = \partial\Omega_{se}^{2D}\otimes[0,L_z],
    \quad
    \partial\Omega_{sg} = \partial\Omega_{sg}^{2D}\otimes[0,L_z].
\end{split}
\right.
\end{equation}
In the above relations $\Omega$ is the 3D domain, and
$\Omega_{2D}$ is the computational domain in the $xy$ plane (i.e.~projection of $\Omega$ onto the $xy$ plane).
Similarly, $\partial\Omega_{2D}$, $\partial\Omega_s^{2D}$, $\partial\Omega_o^{2D}$, $\partial\Omega_{se}^{2D}$ and $\partial\Omega_{sg}^{2D}$ are projections onto the $xy$ plane of the 3D boundaries
$\partial\Omega$, $\partial\Omega_s$, $\partial\Omega_o$, $\partial\Omega_{se}$ and $\partial\Omega_{sg}$, respectively.
In addition, we have the following relations,
\begin{equation}\small
\mbs n = (\mbs n_{2D},0), \quad
\mbs n_o = (\mbs n_o^{2D},0), \quad
\mbs n_s = (\mbs n_s^{2D},0).
\end{equation}
Here $\mbs n$, $\mbs n_o$ and $\mbs n_s$
denote the outward-pointing unit vectors normal to $\partial\Omega$, $\partial\Omega_o$ and $\partial\Omega_s$,
respectively. $\mbs n_{2D}$, $\mbs n_o^{2D}$ and $\mbs n_s^{2D}$ are the outward-pointing unit vectors normal to $\partial\Omega_{2D}$,
$\partial\Omega_o^{2D}$ and $\partial\Omega_s^{2D}$,
respectively.

Let $N_z$ denote the number of Fourier grid points in $z$. We introduce the Fourier basis functions,
\begin{equation}\small
    \Phi_k(z)=e^{i\beta_kz},\quad \beta_k=\frac{2\pi k}{L_z},\quad -\frac{N_z}{2}\leq k \leq \frac{N_z}{2}-1.
\end{equation}
Then, for a generic function $f(x,y,z)$ we have the Fourier expansion relation, 
\begin{equation}\small
f(x,y,z)=\sum_{k=-N_z/2}^{N_z/2-1}\hat{f}_k(x,y)\Phi_k(z), \quad
\int_0^{L_z}f(x,y,z)\bar\Phi_k(z)dz = L_z\hat{f}_k(x,y),
\end{equation}
where $\bar\Phi_k$ is the complex conjugate of $\Phi_k$, and $\hat{f}_k(x,y)$ denotes the $k$-th Fourier mode of $f(x,y,z)$.

We define the basis and test functions in 3D by,
for $-\frac{N_z}{2}\leq k \leq \frac{N_z}{2}-1$,
\begin{equation}\small
\left\{
\begin{split}
&
Q_k(x,y,z)=l(x,y)\Phi_k(z), \quad
\text{(basis function)},\\
&
\bar{Q}_k(x,y,z)=l(x,y)\bar{\Phi}_k(z), \quad\text{(test function)},
\end{split}
\right.
\end{equation}
where $l(x,y)$ denotes an arbitrary function in the $xy$ plane.
Define $\nabla = \left(\nabla_{2D},\frac{\partial}{\partial z}\right)=\left(\frac{\partial}{\partial x},\frac{\partial}{\partial y},\frac{\partial}{\partial z}\right)$.
Let $f(x,y,z)$ denote a generc scalar field
and $\mbs u(x,y,z)=(\mbs u_{2D}(x,y,z),u_z(x,y,z))=(u_x(x,y,z),u_y(x,y,z),u_z(x,y,z))$ denote the velocity (or a generic vector) field.
Then the following relations hold,
\begin{equation}\small
\left\{
\begin{aligned}
&\int_{\Omega}f(x,y,z)\bar{Q}_k(x,y,z)d\Omega=L_z\int_{\Omega_{2D}}\hat{f}_k(x,y)l(x,y)d\Omega_{2D}\label{FFT_int}\\
&\int_{\Omega}\nabla f(x,y,z)\cdot\nabla \bar{Q}_k(x,y,z)d\Omega=L_z\int_{\Omega_{2D}}\left[\nabla_{2D}\hat{f}_k(x,y)\cdot \nabla_{2D}l(x,y)+\beta_k^2\hat{f}_k(x,y)l(x,y)\right]d\Omega_{2D}\\
&\int_{\Omega}\textbf{u}\cdot \nabla\bar{Q}_kd\Omega=L_z\int_{\Omega_{2D}}\left[\nabla_{2D}l(x,y)\cdot\hat{\textbf{u}}_{2D,k}-i\beta_kl(x,y)\hat{u}_{z,k}
\right]d\Omega_{2D}\\
\end{aligned}
\right.
\end{equation} 
where $\hat{\textbf{u}}_{2D,k}$ and $\hat{u}_{z,k}$ are the Fourier modes of $\textbf{u}_{2D}$ and $u_z$, respectively, and $d\Omega=d\Omega_{2D}dz$.

By using the above integral relations, we can reduce the 3D weak forms in~\eqref{eq_41}--\eqref{eq_52} into 2D weak forms about the Fourier modes. Let us assume  in the following that $\omega(x,y)$ denote an arbitrary 2D test function for the electric potential, the electric field and
the phase field functions, and $v(x,y)$ denote an arbitrary 2D test function about the pressure and velocity fields.  For simplicity, we will assume that $\omega(x,y)$ and $v(x,y)$  vanish on the corresponding Dirichlet type boundaries. 
We use the 2D function $\hat{(\cdot)}_k$ or $\hat{(\cdot)}_{,k}$ of $(x,y)$ to denote
the $k$-th Fourier mode of the 3D functon $(\cdot)$ of $(x,y,z)$.

Let $\textbf{R}=(\epsilon(\phi^{*,n+1})-\varepsilon_0)\nabla V^{*,n+1}=(\mbs R_{2D},R_z)$. 
The weak form~\eqref{eq_41} is reduced to, 
\begin{align}\label{eq_66}\small
&\int_{\Omega_{2D}}\varepsilon_0\nabla_{2D}\hat{V}^{n+1}_k\cdot\nabla_{2D}\omega+\beta_k^2\int_{\Omega_{2D}}\varepsilon_0\hat{V}^{n+1}_k\omega=-\int_{\Omega_{2D}}\hat{\textbf{R}}_k\cdot\nabla\omega-\int_{\Omega_{2D}}\hat{f}_{V,k}^{n+1}\omega, \quad\forall \omega(x,y),
\end{align}
where $\nabla\omega(x,y)=(\nabla_{2D}\omega,-i\beta_k\omega)$, and we have used the
following equation,
\begin{equation}\small
\begin{split}
\int_{\Omega}\textbf{R}\cdot\nabla \bar{Q}_k&=\sum_{m=-N_z/2}^{N_z/2-1}\left(\int_{\Omega_{2D}}\hat{\textbf{R}}_m\cdot\nabla \omega\right)\left(\int_0^{L_z}\Phi_m(z)\bar{\Phi}_k(z)\right) \\
&=L_z\int_{\Omega_{2D}}\hat{\textbf{R}}_k\cdot\nabla \omega 
    =L_z\int_{\Omega_{2D}}\left(
    \hat{\mbs R}_{2D,k}\cdot\nabla_{2D}\omega -i\beta_k \hat{R}_{z,k}\omega
    \right).
    \label{3D_integral_equation}
\end{split}\end{equation}

The 3D weak form~\eqref{eq_43} now becomes
\begin{equation}\small
    \int_{\Omega_{2D}}(\hat{E}_{x,k}^{n+1},\hat{E}_{y,k}^{n+1},\hat{E}_{z,k}^{n+1}) \omega=\int_{\Omega_{2D}} (\partial_x \hat{V}_k^{n+1},\partial_y \hat{V}_k^{n+1},-i\beta_k\hat{V}_k^{n+1})\omega, 
    \quad \forall \omega(x,y),
    \label{eq 70}
\end{equation}
where $\mbs{\hat{E}}^{n+1}_k=(\hat{E}_{x,k}^{n+1},\hat{E}_{y,k}^{n+1},\hat{E}_{z,k}^{n+1})$.

The weak forms~\eqref{eq_46}--\eqref{eq_47} are reduced to,
\begin{align}\small
&\int_{\Omega_{2D}}\nabla_{2D}\hat{\psi}_k^{n+1}\cdot \nabla_{2D}\omega 
+\left(\alpha+\dfrac{S}{\eta^2}+\beta_k^2\right)\int_{\Omega_{2D}}\hat{\psi}_k^{n+1}\omega=\int_{\Omega_{2D}}(\beta_k^2\hat{Q}_{2,k}-\hat{Q}_{1,k})\omega+\int_{\Omega_{2D}}\nabla_{2D}\hat{Q}_{2,k}\cdot\nabla\omega\notag\\
&\quad\quad\quad\quad+\int_{\partial\Omega_{o}^{2D}}\left[\hat{g}_{1,k}^{n+1}+\left(\alpha+\dfrac{S}{\eta^2}\right)\hat{g}_{2,k}^{n+1}\right]\omega +\int_{\partial\Omega_{s}^{2D}} \left[\hat{g}_{1,k}^{n+1}
   +\left(\alpha+\dfrac{S}{\eta^2}\right)\hat{U}_k\right]\omega, \quad\forall \omega(x,y);
\end{align}
\begin{align}\label{eq_69}
&\int_{\Omega_{2D}}\nabla_{2D}\hat{\phi}_k^{n+1}\cdot \nabla_{2D}\omega +(-\alpha+\beta_k^2)\int_{\Omega_{2D}}\hat{\phi}_k^{n+1}\omega=\int_{\Omega_{2D}}\hat{\psi}_k^{n+1}\omega+\int_{\partial\Omega_{o}^{2D}}\hat{g}_{2,k}^{n+1}\omega
+\int_{\partial\Omega_{s}^{2D}} \hat{U}_k\omega, \notag\\
&\qquad\qquad\forall \omega(x,y),
\end{align}
where $U=-g_{3}^{n+1}-\frac{\Theta'(\phi^{*,n+1})}{\lambda}$, and $\hat{U}_k$ denotes the Fourier modes of $U$.

Let 
\begin{equation}\small
\left\{
\begin{split}
&
\textbf{G}=(\textbf{G}_{2D},G_z)=\textbf{T}+\nabla\left(\dfrac{\mu^{n+1}}{\rho^{n+1}}\right)\times\boldsymbol{\omega}^{*,n+1}, \qquad
\textbf{Y}=\textbf{G}-\dfrac{1}{\rho_0}\nabla P^{n+1}, \\
&
\textbf{J}=(\textbf{J}_{2D},J_z)=\dfrac{\mu^{n+1}}{\rho^{n+1}}\textbf{n}\times\boldsymbol{\omega}^{*,n+1}, \quad
\textbf{K}=\left(\dfrac{\mu^{n+1}}{\rho^{n+1}}-\nu_m\right)\boldsymbol{\omega}^{*,n+1}, \quad
\textbf{L}=\left(\dfrac{\mu^{n+1}}{\rho^{n+1}}-\nu_m\right)\textbf{n}\times\boldsymbol{\omega}^{*,n+1}.
\end{split}
\right.
\end{equation}
The weak form~\eqref{eq_48} for the pressure is reduced to,
\begin{align}\small
\label{3D_u}
&\int_{\Omega_{2D}}\nabla_{2D}\hat{P}_k^{n+1}\cdot \nabla_{2D}v +\beta_k^2\int_{\Omega_{2D}}\hat{P}_k^{n+1} v=\rho_0\int_{\Omega_{2D}}\hat{\textbf{G}}_{2D,k}\cdot\nabla_{2D}v-i\beta_k\rho_0\int_{\Omega_{2D}}\hat{G}_{z,k}v\notag\\
   &\quad\quad\quad\quad -\rho_0\int_{\partial\Omega_{o}^{2D}}\hat{\textbf{J}}_{2D,k}\cdot\nabla_{2D}v+i\beta_k\rho_0\int_{\partial\Omega_{o}^{2D}}\hat{J}_{z,k}v
   -\dfrac{\rho_0\gamma_0}{\Delta t}\int_{\partial\Omega_{s}^{2D}}\textbf{n}_{2D}\cdot \hat{\mbs w}_{2D,k}^{n+1}v, \quad \forall v(x,y).
\end{align}
The weak form~\eqref{eq_52} for the velocity is reduced to,
\begin{align}\small
&\int_{\Omega_{2D}}\nabla_{2D} v\cdot\nabla_{2D}\hat{\textbf{u}}_k^{n+1} +\left(\beta_k^2+\dfrac{\gamma_0}{\nu_m\Delta t}\right)\int_{\Omega_{2D}}\hat{\textbf{u}}_k^{n+1} v=\dfrac{1}{\nu_m}\int_{\Omega_{2D}}\hat{\textbf{Y}}_k v-\dfrac{1}{\nu_m}\int_{\Omega_{2D}}\hat{\textbf{K}}_k\times\nabla v\notag\\
    &\qquad \qquad \qquad \qquad \qquad \qquad \qquad +\int_{\partial \Omega_{o}^{2D}}\hat{\mbs f}_{1,k}^{n+1}v-\dfrac{1}{\nu_m}\int_{\partial \Omega_{o}^{2D}}\hat{\textbf{L}}_k v, \quad \forall v(x,y).
    \label{3D_v}
\end{align}
Note that  the terms $i\beta_k\rho_0\int_{\Omega_{2D}}\hat{G}_{z,k}v$ and $i\beta_k\rho_0\int_{\partial\Omega_{o}^{2D}}\hat{J}_{z,k}^{n+1}v$ in equation~\eqref{3D_u} and  the term $\dfrac{1}{\nu_m}\int_{\Omega_{2D}}\hat{\textbf{K}}_k\times\nabla v$ in equation~\eqref{3D_v} 
 mixes up the imaginary and real parts, which calls for special attention  in the implementation.


To formulate the fully discretized equations in 3D, we partition the domain $\Omega_{2D}$ in the $xy$ plane by a mesh of $C^0$ spectral elements. Let $\Omega_{2Dh}$ denote the discretized $\Omega_{2D}$,
$\Omega_{2Dh}=\cup_{e=1}^{N_e}\Omega_{2Dh}^e$, where $\Omega_{2Dh}^e$ denotes the element $e$ in the $xy$ plane.
Let $\partial\Omega_{2Dh}$, $\partial\Omega_{oh}^{2D}$, and $\partial\Omega_{sh}^{2D}$ denote the discretized versions of $\partial\Omega_{2D}$, $\partial\Omega_{o}^{2D}$, and $\partial\Omega_{s}^{2D}$, respectively.
Let $\partial\Omega_{seh}^{2D}$ and
$\partial\Omega_{sgh}^{2D}$ denote the discretized solid-electrode and solid-gap boundaries in $\Omega_{2D}$,
$\partial\Omega_{sh}^{2D}= \partial\Omega_{seh}^{2D}\cup \partial\Omega_{sgh}^{2D}$.
Let $\Pi_{K}(\Omega_{2Dh}^e)$ denote the polynomial space defined on $\Omega_{2Dh}^e$ with their degrees characterized by $K$. We define
\begin{equation}\small
\left\{
\begin{split}
&
\mbb Y_h = \{\ v\in H^1(\Omega_{2Dh})\ :\ v|_{\Omega_{2Dh}^e}\in\Pi_K(\Omega_{2Dh}^e),
\ 1\leqslant e\leqslant N_e \ \}, \\
&
\mbb Y_{h0}^E = \{\  v\in \mbb Y_h\ :\ 
v|_{\partial\Omega_{seh}^{2D}}=0
\ \}, \\
&
\mbb Y_{h0}^{P} = \{\ v\in \mbb Y_h\ :\ 
v|_{\partial\Omega_{oh}^{2D}}=0
\ \}, \\
&
\mbb Y_{h0}^u = \{\  v\in \mbb Y_h\ :\ 
v|_{\partial\Omega_{sh}^{2D}}=0
\ \}.
\end{split}
\right.
\end{equation}
In the following the subscript $h$ denotes the discretized version of a variable.

Then the fully discretized system in 3D consists of the following equations:\\[5pt]
\noindent\underline{For $V_h^{n+1}$:} 
find $\hat{V}_{kh}^{n+1}\in\mbb Y_{h}$ such that
(for $-N_z/2\leqslant k\leqslant N_z/2-1$)
\begin{subequations}\label{eq_74}
\begin{align}\small
&\int_{\Omega_{2Dh}}\varepsilon_0\nabla_{2D}\hat{V}_{kh}^{n+1}\cdot\nabla_{2D}\omega_h+\beta_k^2\int_{\Omega_{2Dh}}\varepsilon_0\hat{V}_{kh}^{n+1}\omega_h=-\int_{\Omega_{2Dh}}\hat{\textbf{R}}_{kh}\cdot\nabla\omega_h-\int_{\Omega_{2Dh}}\hat{f}_{V,kh}^{n+1}\omega_h, \notag\\
&\qquad\forall \omega_h\in\mbb Y_{h0}^E; \label{eq_74a}
\\
& \hat{V}_{kh}^{n+1} = \left\{
\begin{array}{ll}
\mathcal V_e, & k=0, \\
0, & k\neq 0.
\end{array}
\right. \quad\text{on}\ \partial\Omega_{seh}^{2D}.
\label{eq_74b}
\end{align}
\end{subequations}
\noindent\underline{For $\mbs E_h^{n+1}$:}
find $\mbs{\hat{E}}_{kh}^{n+1}=(\hat{E}_{x,kh}^{n+1},\hat{E}_{y,kh}^{n+1},\hat{E}_{z,kh}^{n+1})\in [\mbb Y_h]^3  $, such that
(for $-N_z/2\leqslant k\leqslant N_z/2-1$)
\begin{equation}\small
    \int_{\Omega_{2Dh}}(\hat{E}_{x,kh}^{n+1},\hat{E}_{y,kh}^{n+1},\hat{E}_{z,kh}^{n+1}) \omega_h=\int_{\Omega_{2Dh}} (\partial_x \hat{V}_{kh}^{n+1},\partial_y \hat{V}_{kh}^{n+1},-i\beta_k\hat{V}_{kh}^{n+1})\omega_h,\quad \forall \omega_h\in \mbb Y_h.
    \label{eq 78}
\end{equation}
\\
\noindent\underline{For $\psi_h^{n+1}$:}
find $\hat{\psi}_{kh}^{n+1}\in \mbb Y_h$ such that
(for $-N_z/2\leqslant k\leqslant N_z/2-1$)
\begin{align}\small
\label{eq_75}
&\int_{\Omega_{2Dh}}\nabla_{2D}\hat{\psi}_{kh}^{n+1}\cdot \nabla_{2D}\omega_h 
+\left(\alpha+\dfrac{S}{\eta^2}+\beta_k^2\right)\int_{\Omega_{2Dh}}\hat{\psi}_{kh}^{n+1}\omega_h \notag \\
&= \int_{\Omega_{2Dh}}(\beta_k^2\hat{Q}_{2,kh}-\hat{Q}_{1,kh})\omega_h 
+\int_{\Omega_{2Dh}}\nabla_{2D}\hat{Q}_{2,kh}\cdot\nabla\omega_h\notag\\
&\quad+\int_{\partial\Omega_{oh}^{2D}}\left[\hat{g}_{1,kh}^{n+1}+\left(\alpha+\dfrac{S}{\eta^2}\right)\hat{g}_{2,kh}^{n+1}\right]\omega_h +\int_{\partial\Omega_{sh}^{2D}} \left[\hat{g}_{1,kh}^{n+1}
   +\left(\alpha+\dfrac{S}{\eta^2}\right)\hat{U}_{kh}\right]\omega_h, \quad\forall \omega_h\in \mbb Y_h.
\end{align}
\noindent\underline{For $\phi_h^{n+1}$:}
find $\hat{\phi}_{kh}^{n+1}\in \mbb Y_h$ such that
(for $-N_z/2\leqslant k\leqslant N_z/2-1$)
\begin{align}\small
\label{eq_76}
&\int_{\Omega_{2Dh}}\nabla_{2D}\hat{\phi}_{kh}^{n+1}\cdot \nabla_{2D}\omega_h +(-\alpha+\beta_k^2)\int_{\Omega_{2Dh}}\hat{\phi}_{kh}^{n+1}\omega_h \notag \\
&=\int_{\Omega_{2Dh}}\hat{\psi}_{kh}^{n+1}\omega_h+\int_{\partial\Omega_{oh}^{2D}}\hat{g}_{2,kh}^{n+1}\omega_h
+\int_{\partial\Omega_{sh}^{2D}} \hat{U}_{kh}\omega_h,\quad\forall \omega_h\in\mbb Y_h.
\end{align}
\noindent\underline{For $P_h^{n+1}$:}
find $\hat{P}_{kh}^{n+1}\in \mbb Y_h$ such that
(for $-N_z/2\leqslant k\leqslant N_z/2-1$)
\begin{subequations}
\begin{align}\small
&\int_{\Omega_{2Dh}}\nabla_{2D}\hat{P}_{kh}^{n+1}\cdot \nabla_{2D}v_h +\beta_k^2\int_{\Omega_{2Dh}}\hat{P}_{kh}^{n+1} v_h=\rho_0\int_{\Omega_{2Dh}}\hat{\textbf{G}}_{2D,kh}\cdot\nabla_{2D}v_h-i\beta_k\rho_0\int_{\Omega_{2Dh}}\hat{G}_{z,kh}v_h\notag\\
   &\ \ -\rho_0\int_{\partial\Omega_{oh}^{2D}}\hat{\textbf{J}}_{2D,kh}\cdot\nabla_{2D}v_h+i\beta_k\rho_0\int_{\partial\Omega_{oh}^{2D}}\hat{J}_{z,kh}v_h
   -\dfrac{\rho_0\gamma_0}{\Delta t}\int_{\partial\Omega_{sh}^{2D}}\textbf{n}_{2Dh}\cdot \hat{\mbs w}_{2D,kh}^{n+1}v_h, \quad \forall v_h\in\mbb Y_{h0}^P;
   \label{eq_77a} \\
   &
   \hat{P}_{kh}^{n+1} = \hat{f}_{2,kh}^{n+1}, \quad
   \text{on}\ \partial\Omega_{oh}^{2D}.
   \label{eq_77b}
\end{align}
\end{subequations}
\noindent\underline{For $\mbs u_{h}^{n+1}$:}
find $\hat{\mbs u}_{kh}^{n+1}\in [\mbb Y_h]^3$ such that
(for $-N_z/2\leqslant k\leqslant N_z/2-1$)
\begin{subequations}\label{eq_78}
\begin{align}\small
&\int_{\Omega_{2Dh}}\nabla_{2D} v_h\cdot\nabla_{2D}\hat{\textbf{u}}_{kh}^{n+1}+\left(\beta_k^2+\dfrac{\gamma_0}{\nu_m\Delta t}\right)\int_{\Omega_{2Dh}}\hat{\textbf{u}}_{kh}^{n+1} v_h=\dfrac{1}{\nu_m}\int_{\Omega_{2Dh}}\hat{\textbf{Y}}_{kh} v_h \notag \\
&\quad\quad
-\dfrac{1}{\nu_m}\int_{\Omega_{2Dh}}\hat{\textbf{K}}_{kh}\times\nabla v_h
    +\int_{\partial \Omega_{oh}^{2D}}\hat{\mbs f}_{1,kh}^{n+1}v_h-\dfrac{1}{\nu_m}\int_{\partial \Omega_{oh}^{2D}}\hat{\textbf{L}}_{kh} v_h, \quad \forall v_h\in\mbb Y_{h0}^u.
    \label{eq_78a} \\
    &
    \hat{\mbs u}_{kh}^{n+1} = \hat{\mbs w}_{kh}^{n+1},
    \quad \text{on}\ \partial\Omega_{sh}^{2D}.
    \label{eq_78b}
\end{align}
\end{subequations}

Given $(V^{n},\phi^n,P^n,\mbs u^n)$ in 3D, the field variables $V^{n+1}$, $\mbs E^{n+1}$, $\psi^{n+1}$, $\phi^{n+1}$, $P^{n+1}$ and $\mbs u^{n+1}$ are computed by solving the equations~\eqref{eq_74}--\eqref{eq_78} individually and successively in an un-coupled fashion.
Algorithm~\ref{alg_2} summarizes the solution procedure for 3D problems.

\begin{algorithm}[tb]\small
  \DontPrintSemicolon
  \SetKwInOut{Input}{input}\SetKwInOut{Output}{output}

  \Input{$V^n$, $\phi^n$, $P^n$, $\mbs u^n$.}
  \Output{ $V^{n+1}$, $\mbs E^{n+1}$, $\psi^{n+1}$, $\phi^{n+1}$, $P^{n+1}$, $\mbs u^{n+1}$. }
  \BlankLine  \BlankLine
  solve equations~\eqref{eq_74a}--\eqref{eq_74b}
  for $\hat{V}_k^{n+1}$ ($-N_z/2\leqslant k\leqslant N_z/2$), with Fourier transform to attain $V^{n+1}$  \;
  solve equation~\eqref{eq 78} for $\hat{\mbs E}_k^{n+1}$ ($-N_z/2\leqslant k\leqslant N_z/2$),
  with Fourier transform to attain $\mbs E^{n+1}$\;
  solve equation~\eqref{eq_75} for $\hat{\psi}_k^{n+1}$ ($-N_z/2\leqslant k\leqslant N_z/2$),
  with Fourier transform to attain $\psi^{n+1}$\;
  solve equation~\eqref{eq_76} for $\hat{\phi}_k^{n+1}$ ($-N_z/2\leqslant k\leqslant N_z/2$),
  with Fourier transform to attain $\phi^{n+1}$\;
  solve equations~\eqref{eq_77a}--\eqref{eq_77b}
  for $\hat{P}_k^{n+1}$ ($-N_z/2\leqslant k\leqslant N_z/2$), with Fourier transform to attain $P^{n+1}$\;
  solve equations~\eqref{eq_78a}--\eqref{eq_78b}
  for $\hat{\mbs u}_k^{n+1}$ ($-N_z/2\leqslant k\leqslant N_z/2$), with Fourier transform to attain $\mbs u^{n+1}$\;
  
  \caption{Solution Procedure within a Time Step for 3D Dielectric Flows}
  \label{alg_2}
\end{algorithm}

\begin{remark}
The Algorithms~\ref{alg_1} and~\ref{alg_2}, respectively for 2D and 3D two-phase dielectric flows, share a common characteristic.
The resultant linear algebraic systems for the dynamic variables ($V^{n+1}$, $\mbs E^{n+1}$, $\psi^{n+1}$, $\phi^{n+1}$, $P^{n+1}$, $\mbs u^{n+1}$)
all involve a constant and time-independent coefficient matrix upon discretization, which only needs to be computed once and thus can be pre-computed and saved for later use, despite the variable permittivity/density/viscosity field involved in the system on the continuum level. Because of this property, the current method is  computationally very efficient for simulating two-phase dielectric flow problems.

\end{remark}

\section{Representative Numerical Simulations}
\label{sec:tests}

\subsection{Convergence Test}

\begin{figure}[tb]
\centering 
\subfigure[domain and mesh]{\includegraphics[width=0.4\textwidth]{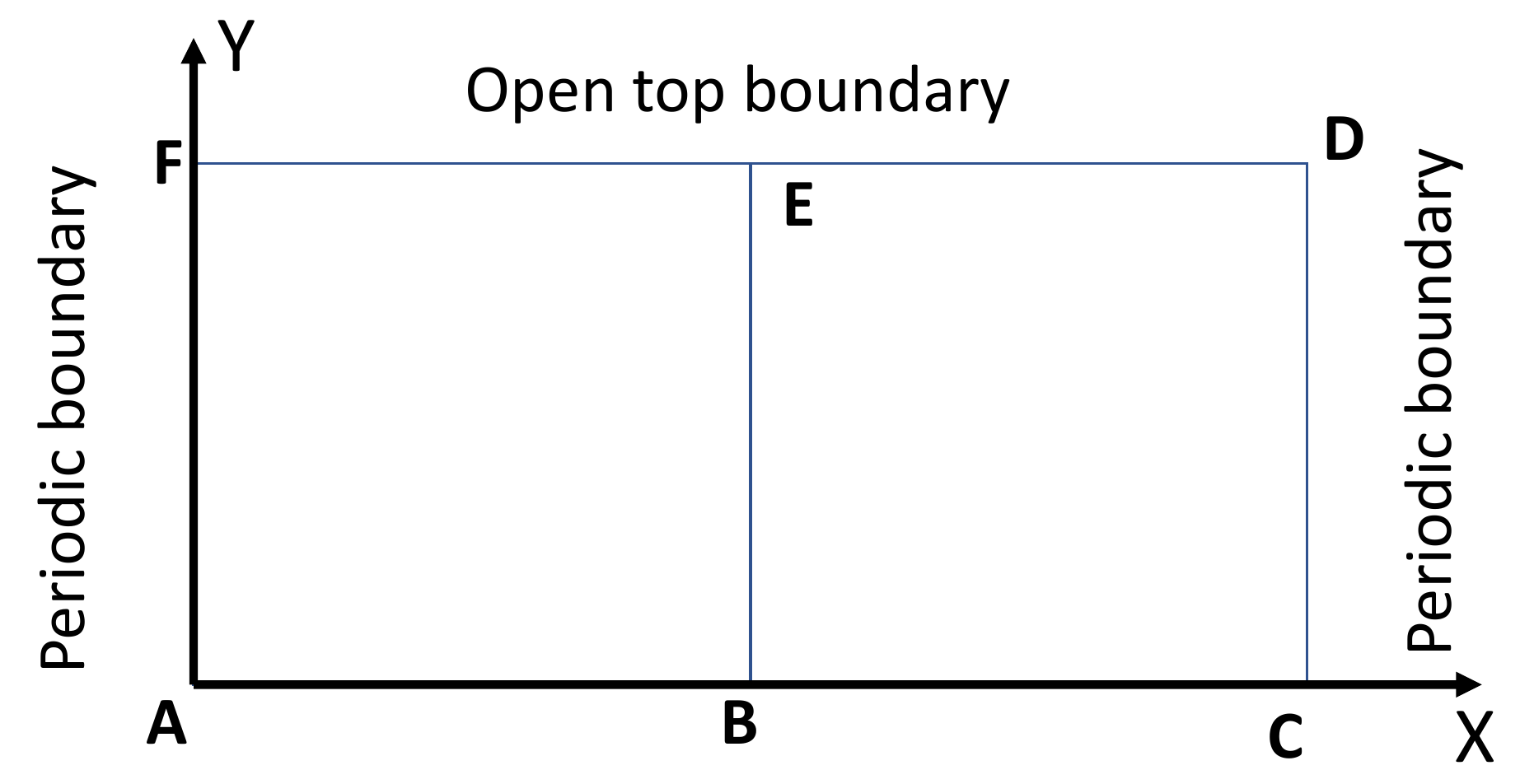}} \\
\subfigure[spatial convergence]{
\includegraphics[width=0.4\textwidth]{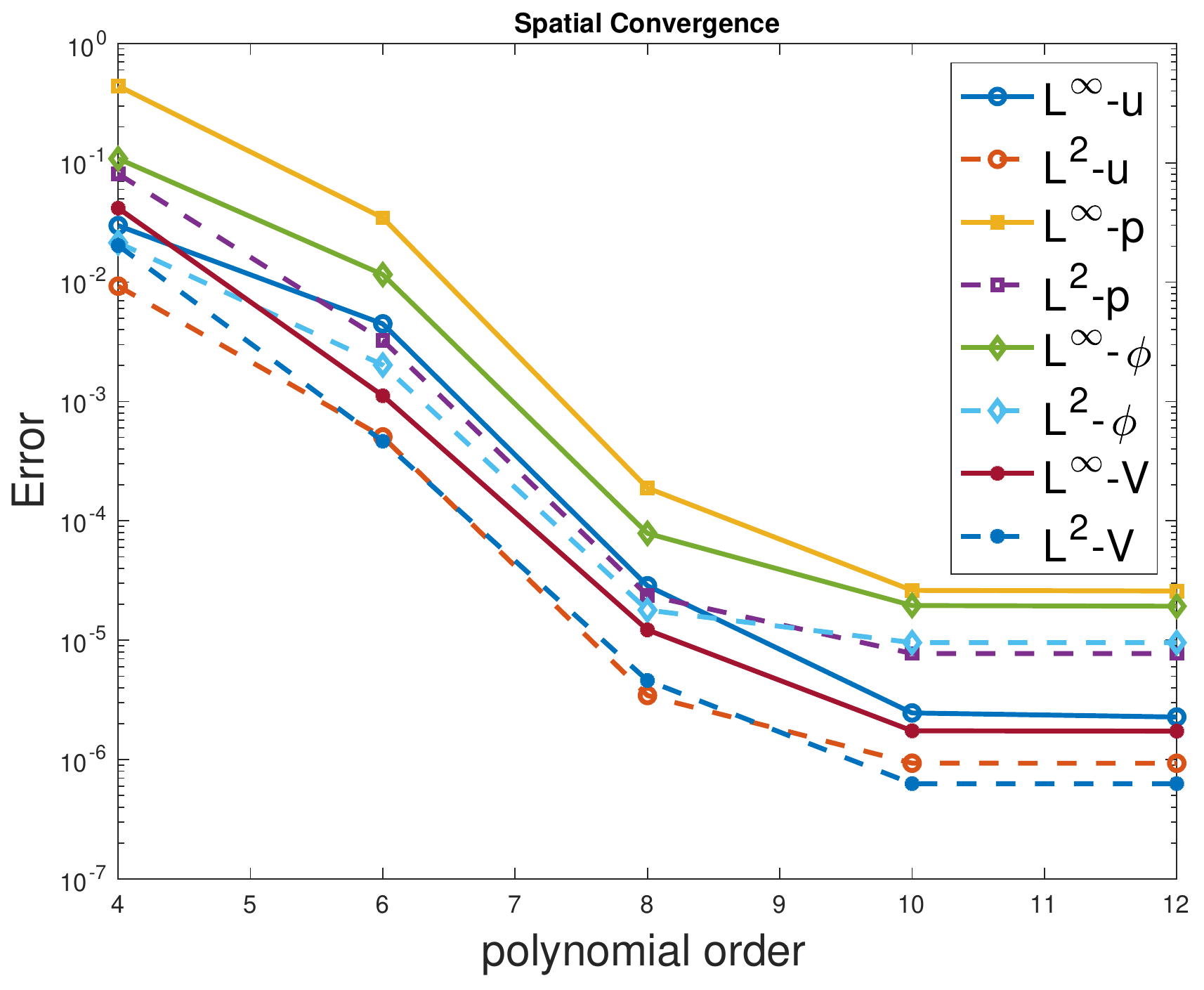}}
\subfigure[temporal convergence]{
\includegraphics[width=0.4\textwidth]{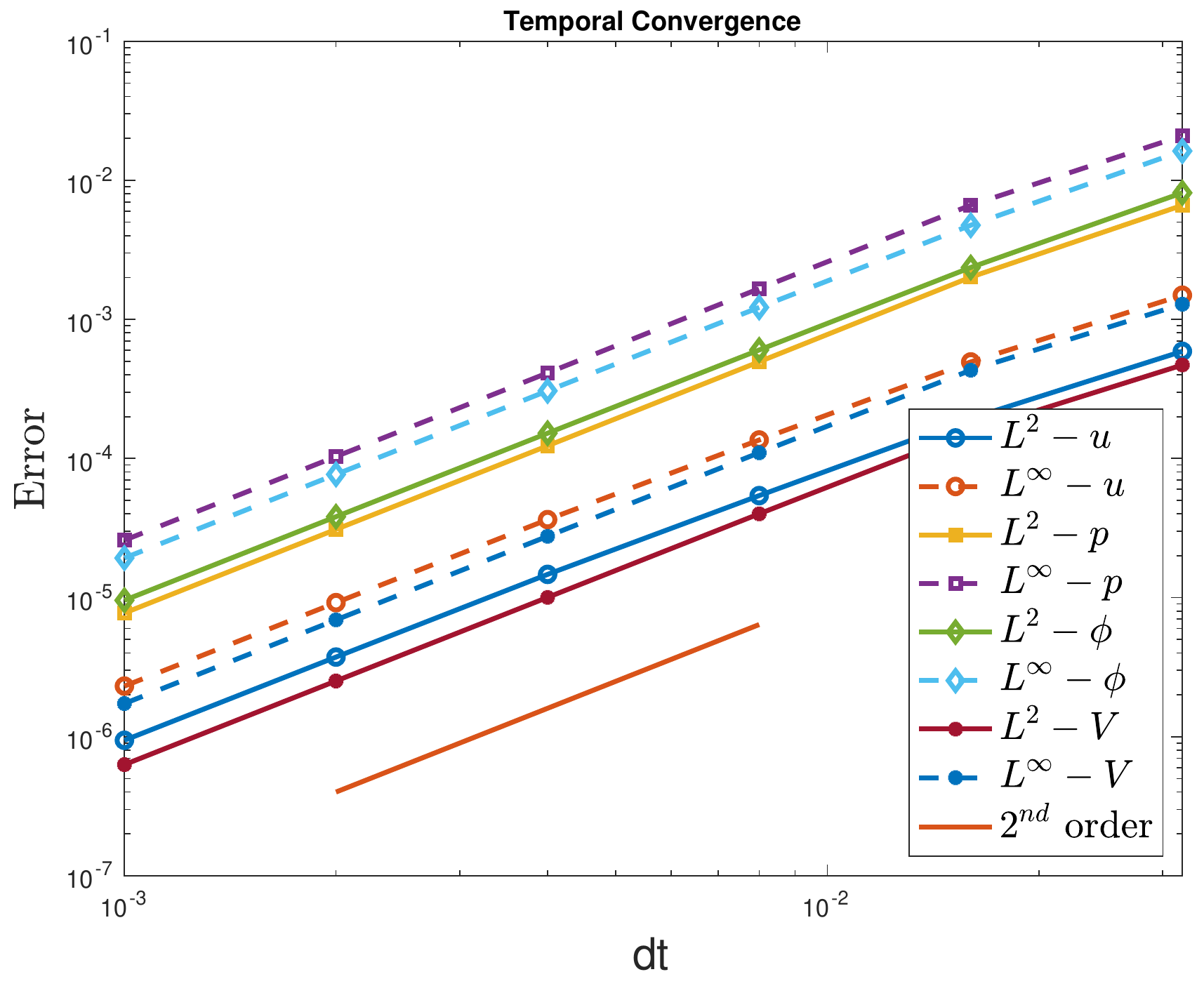}}
\caption{\small 2D convergence test: (a) Sketch of the computational domain and configuration.
(b) $L^{\infty}$ and $L^2$ errors of the dynamic variables versus the element order (fixed $\Delta t=0.001$), showing spatial exponential convergence.
(c) $L^{\infty}$ and $L^2$ errors versus $\Delta t$ (fixed element order $=14$), showing temporal second-order convergence rate.
} 
\label{fg_2} 
\end{figure}

We next employ a manufactured analytic solution to the governing equations to demonstrate the spatial and temporal convergence rates of the numerical method presented in Section~\ref{section3}.

We first look into the convergence for 2D problems. Consider a  domain $\Omega=\{(x,y):0\leq x \leq 2,0\leq y \leq 1\}$ (see Figure~\ref{fg_2}(a)), and the two-phase dielectric governing equations and boundary/initial conditions on $\Omega$ as given by  equations~\eqref{eq_a27}--\eqref{eq_28}, \eqref{eq_8}, \eqref{eq_a8}, \eqref{eq_29}--\eqref{eq_33}, \eqref{eq_a25}, \eqref{eq_21}.
We employ the following manufactured  solution to this problem:
\begin{equation}\small
\label{eq_80}
\left\{
\begin{aligned}
u&=\cos(\pi y)\sin(\pi x)\sin(t),\quad
v=-\sin(\pi y)\cos(\pi x)\sin(t),\quad
P=\sin(\pi y)\cos(\pi x)\cos(t),\\
\phi &=\cos(\pi x)\cos(\pi y)\sin(t),\quad
V=\sin(\pi x)\cos(\pi y),\\
\end{aligned}
\right.
\end{equation} 
where $\mbs u=(u,v)$. All the source terms involved in the governing equations and boundary/initial conditions are chosen such that the field distributions given in~\eqref{eq_80} satisfy the governing equations and boundary/initial conditions.


To simulate this problem, we discretize the domain using two spectral elements of the same size, as shown in Figure~\ref{fg_2}(a).
On the left/right boundaries ($x=0,2$) we impose the periodic condition for all the dynamic variables. The bottom boundary ($y=0$) is assumed to be a wall, and we impose the Dirichlet condition for the velocity and the electric potential (see equations~\eqref{eq_33} and~\eqref{eq_21}),
and the boundary condition~\eqref{eq_32} for
the phase field function.
The top boundary ($y=1$) is assumed to be open, and we impose the boundary conditions~\eqref{eq_a25},~\eqref{eq_29} and~\eqref{eq_30} for the electric potential, the phase field function and the velocity/pressure, respectively.

Figure~\ref{fg_2}(b) shows the $L^{\infty}$ and $L^2$ errors of the velocity, pressure, phase field function, and the electric potential versus the element order in the simulations. Here the time step size is fixed at $\Delta t=0.001$, and the governing equations are integrated from $t=0$ to $t=t_f=0.2$. Shown in this figure are the errors of dynamic variables at $t=t_f$. 
The errors decrease exponentially with increasing element order (when below  $10$), and they stagnate when the element order increases beyond $10$ due to the dominance of the temporal truncation error.  

Figure~\ref{fg_2}(c) illustrates the temporal convergence of the  method. The $L^{\infty}$ and $L^2$ errors of the dynamic variables at $t=t_f=0.4$ as a function of $\Delta t$
are shown. In this group of tests the element order is fixed at $14$.
We observe a second-order convergence rate for the velocity and the electric potential, as well as with the $L^2$ errors for the pressure ($P$) and the phase field function ($\phi$). The $L^{\infty}$ errors of the pressure and the phase field function exhibit an approximate second-order rate, with some irregularities on the error curves.

\begin{figure}[tb]
\centerline{ 
\subfigure[domain and mesh]{\includegraphics[width=0.36\textwidth]{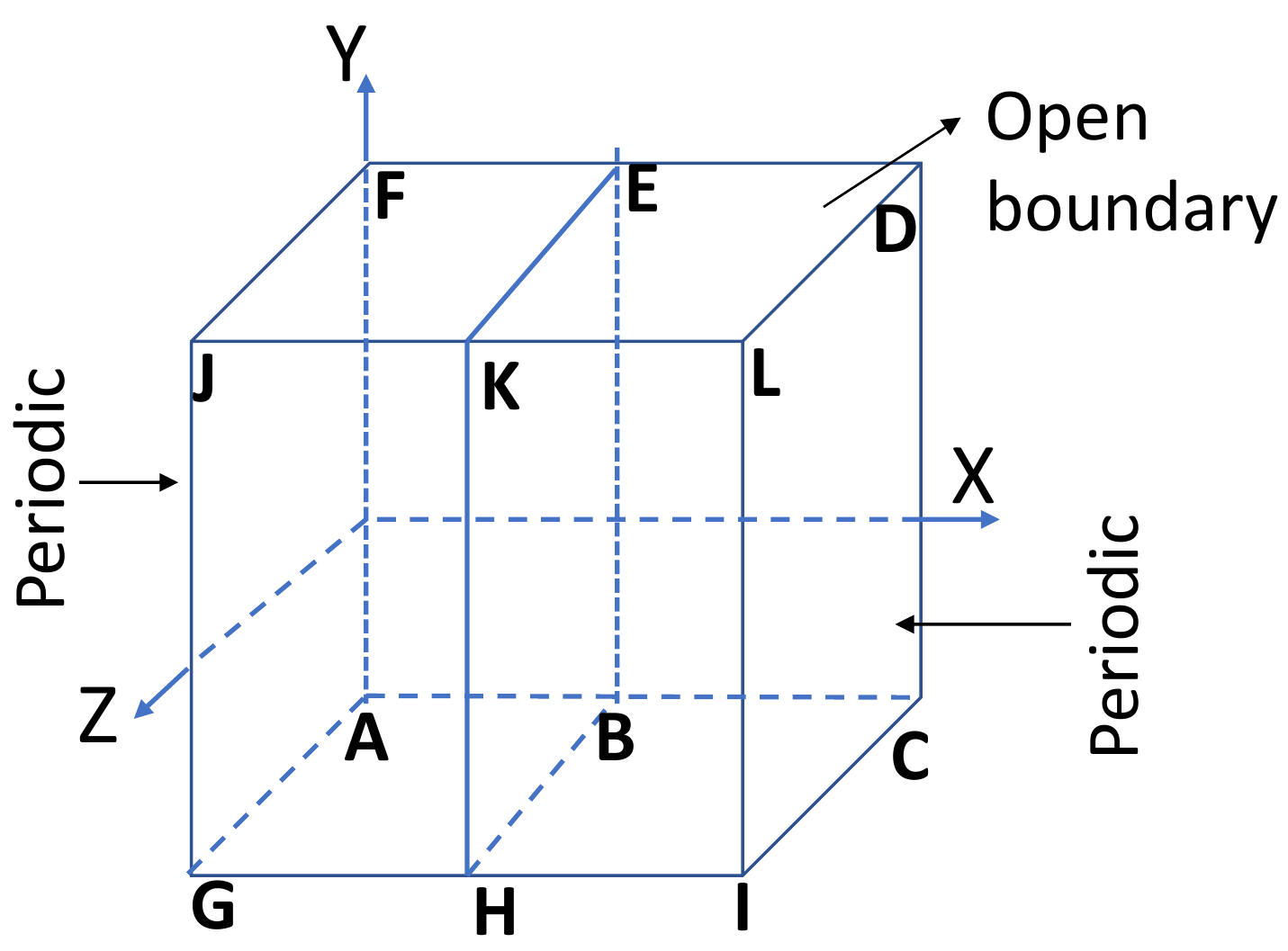}}
\subfigure[spatial convergence]{
\includegraphics[width=0.3\textwidth]{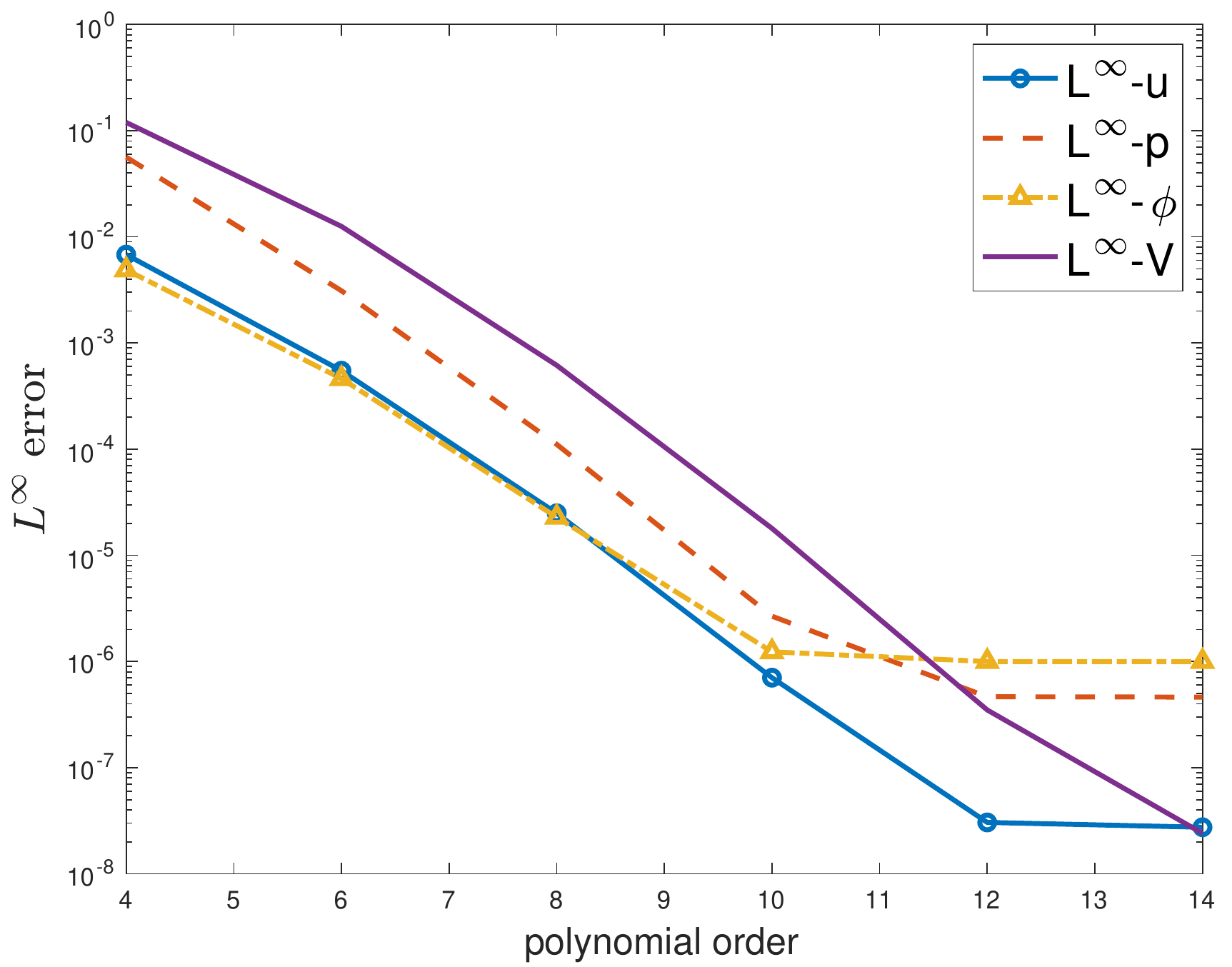}}
\subfigure[temporal convergence]{
\includegraphics[width=0.3\textwidth]{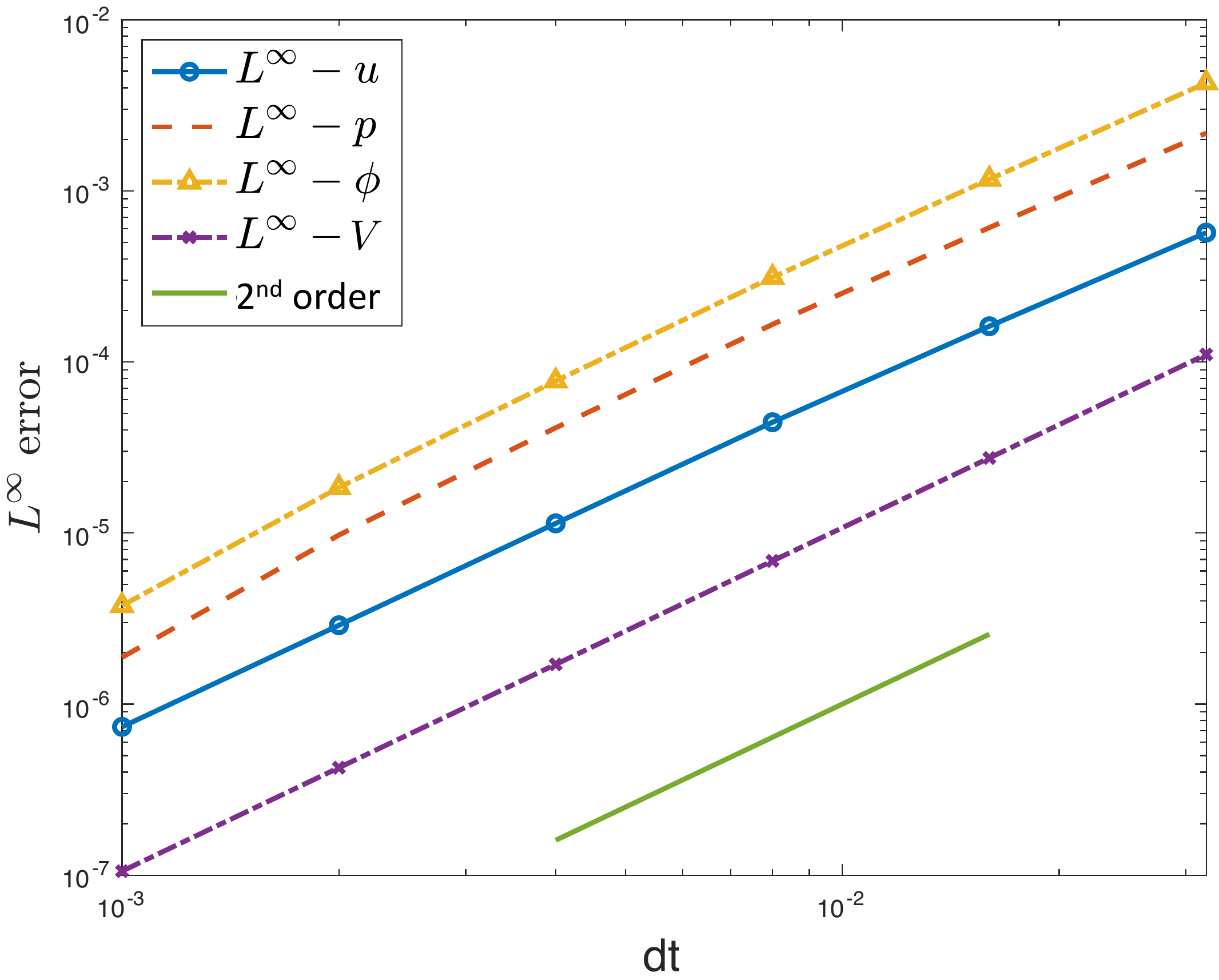}}
}
\caption{\small 3D convergence test:
(a) Domain and configuration. (b) $L^{\infty}$ errors of the dynamic variables versus the element order (fixed $\Delta t=0.001$ and $N_z=8$), showing spatial exponential convergence.
(b) $L^{\infty}$ errors versus $\Delta t$ (fixed element order $=14$ and $N_z=8$), showing temporal second-order convergence rate.
} 
\label{fg_3} 
\end{figure}

To test the spatial/temporal convergence of the 3D algorithm, we consider the domain $\Omega=\{(x,y,z):0\leq x \leq 2,-1\leq y \leq 1, 0\leq z \leq 2\}$, as sketched in Figure~\ref{fg_3}(a). The plane $\overline{HBEK}$ ($x=1$) partitions the domain into two equal sub-domains. The domain $\Omega$ and all the flow variables are assumed to be homogeneous along  $z$. The top boundary ($y=1$) is open.  
The boundaries along the $x$ direction ($x=0$ and $2$) are periodic.
On the bottom face  $\overline{ACIG}$, we impose the Dirichlet boundary condition for the velocity $\mbs u$, and the wall boundary conditions~\eqref{eq_32} for the phase field function $\phi$. 
For the electric potential $V$, we impose the Dirichlet condition (second equation in~\eqref{eq_21}) on the region~$\overline{ABHG}$ and the Neumann condition (first equation in~\eqref{eq_21}) on the region~$\overline{BCIH}$.

We employ the following manufactured analytic solution on $\Omega$ for the 3D convergence tests,
\begin{equation}\small\label{eq_81}
\left\{
\begin{aligned}
u&=\cos(\pi x)\cos(\pi y)cos(\pi z)\sin(t),\quad
v=0,\quad
w=\sin(\pi x)\cos(\pi y)\sin(\pi z)\sin(t),\\
P&=\sin(\pi x)\sin(\pi y)\sin(\pi z)\cos(t),\quad
\phi =\cos(\pi x)\cos(\pi y)\cos(\pi z)\sin(t),\quad
V=\sin(\pi x)\cos(\pi y)\cos(\pi z),\\
\end{aligned}
\right.
\end{equation} 
where $\mbs u=(u,v,w)$. The source terms in the governing equations and the non-homogeneous boundary conditions are set according to these analytic expressions.
In the simulations we employ $N_z=8$ Fourier grid points along the $z$ direction, and two spectral elements in the $xy$ planes, as shown in Figure~\ref{fg_3}(a).

The spatial convergence of the 3D algorithm is illustrated by Figure~\ref{fg_3}(b), in which the $L^{\infty}$ errors of the dynamic variables are shown as a function of the element order.
Here the problem is simulated from
$t=0$ to $t=t_f=0.1$, and the time step size is fixed at $\Delta t=0.001$.
The exponential convergence in space
is evident from the results.

The temporal convergence of the 3D algorithm is illustrated by Figure~\ref{fg_3}(c). Here the $L^{\infty}$ errors of the dynamic variables are shown as a function of
$\Delta t$.  The element order has been fixed at $14$, and the final integration time is $t=t_f=0.1$. 
One can observe the second-order convergence rate with respect to $\Delta t$.


\subsection{Equilibrium Dielectric  Drop on a Wall}
\label{sec:drop}

\begin{figure}[tb]
\centerline{ 
\includegraphics[width=0.38\textwidth]{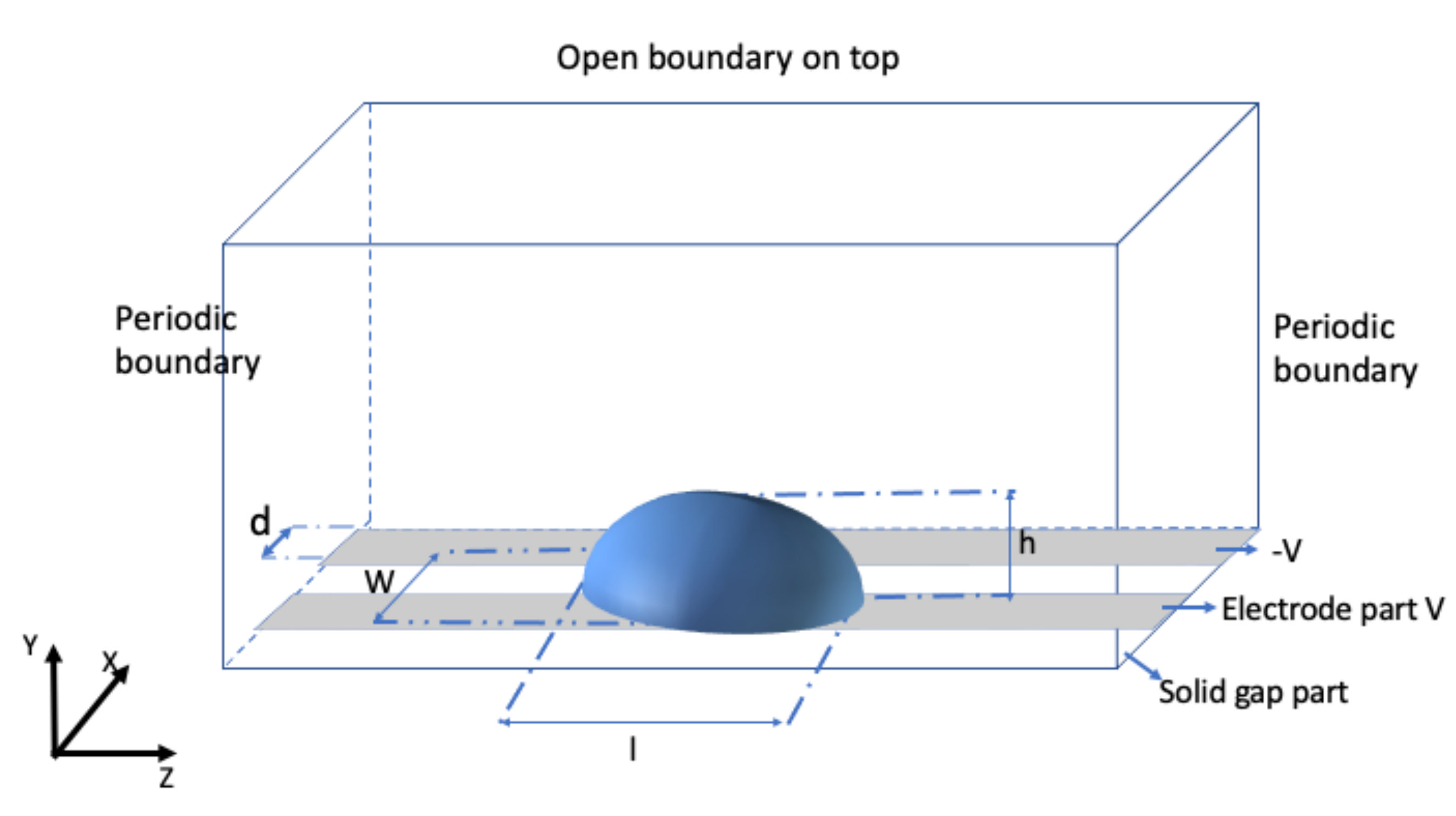}(a)
\includegraphics[width=0.45\textwidth]{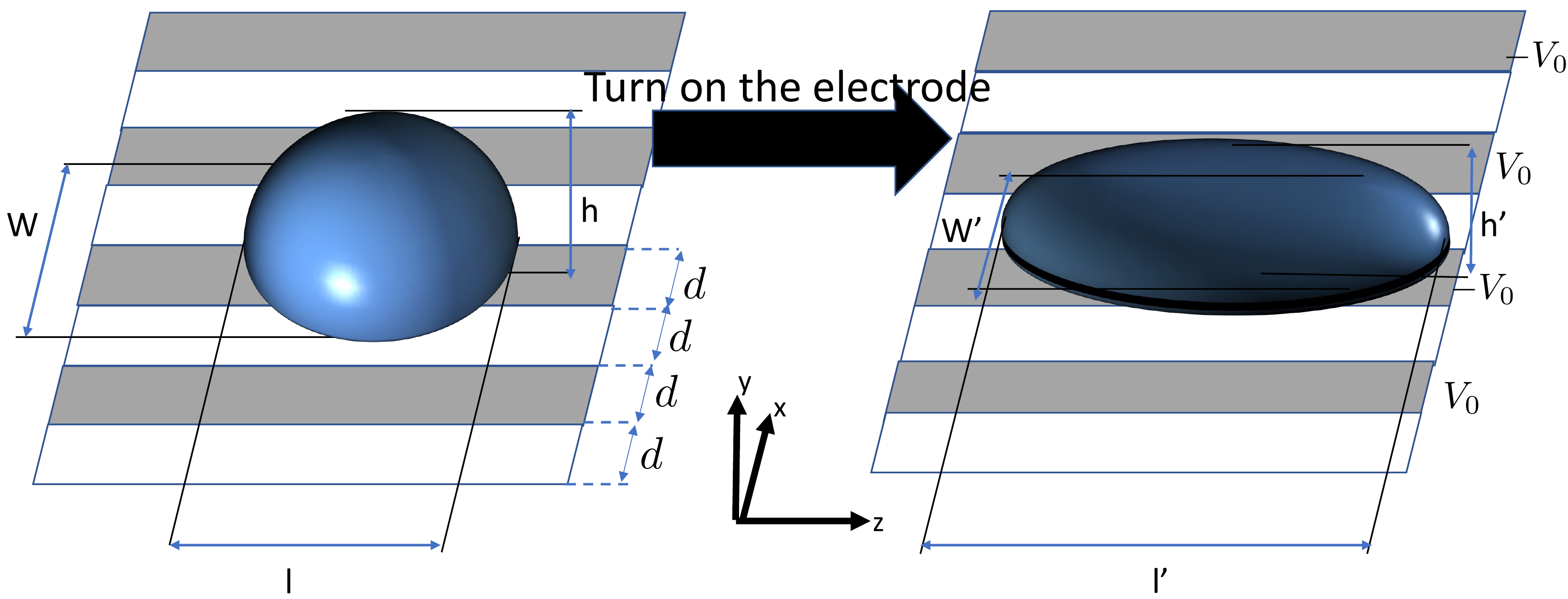}(b)
}
\centerline{
\includegraphics[width=0.45\textwidth]{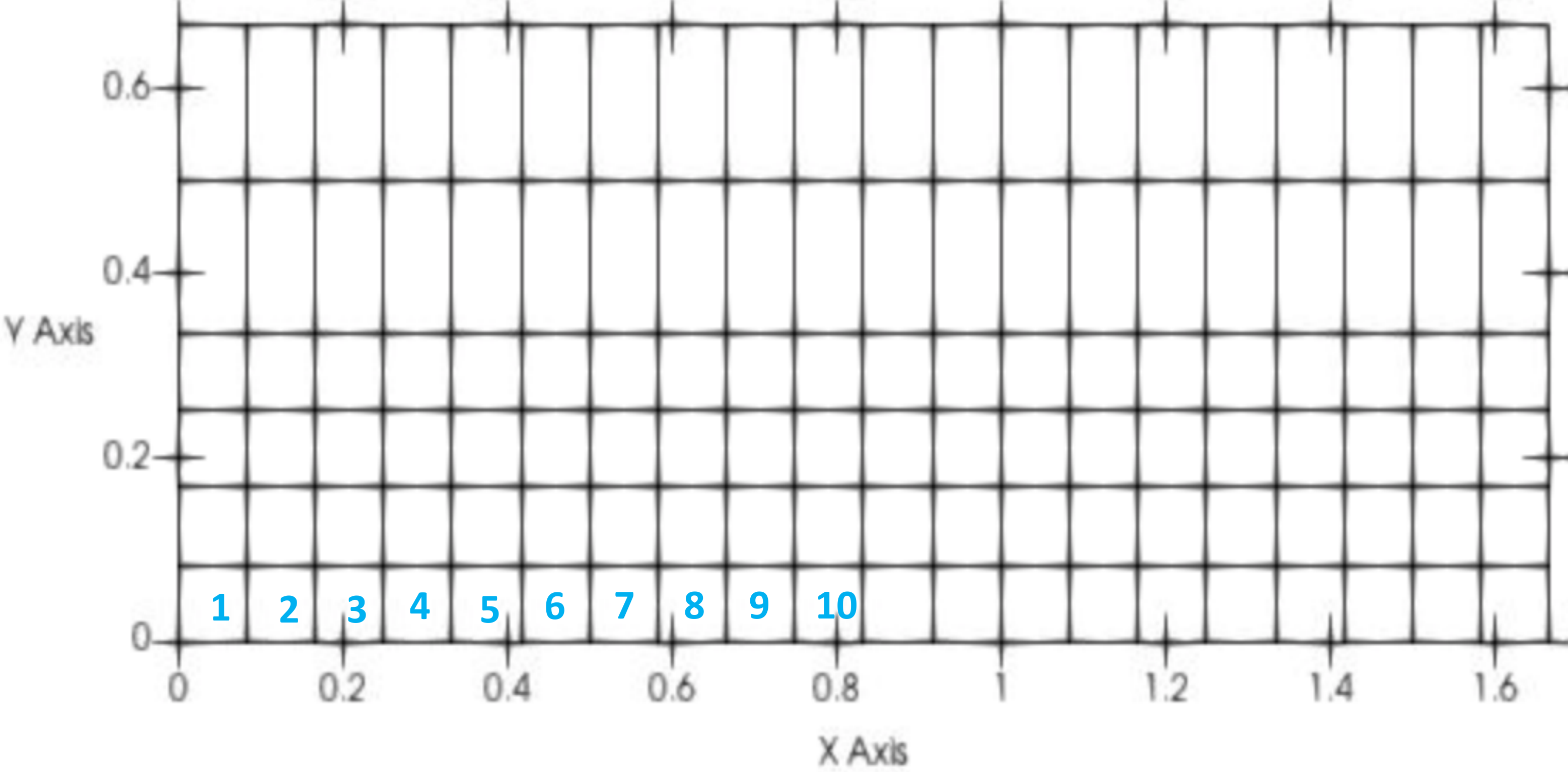}(c)
}
\caption{\small Dielectric liquid drop on a wall:
(a) Sketch of domain and flow configuration.
(b) Cartoon of the drop deformation when the electrodes are switched on.
(c) The spectral element mesh in the $xy$ plane used in the simulations.
On bottom wall, the shaded regions are the electrodes and the white regions denote the gaps between the electrodes. 
In (c), at the bottom wall, the voltage is $V_0$ in element 1,
element 2 is a gap, the voltage is $-V_0$ in element 3, element 4 is a gap, etc.
} 
\label{dew} 
\end{figure}

\begin{figure}[tb]
\centerline{  
\subfigure[$V_0$=100 volt]{
\includegraphics[width=0.25\textwidth]{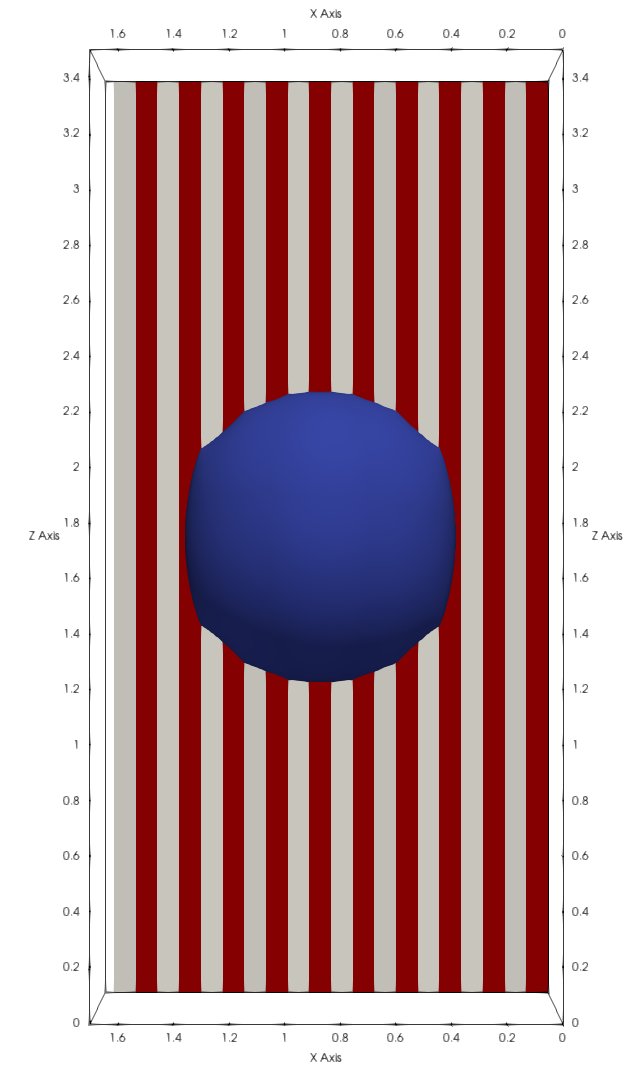}}\qquad\quad
\subfigure[$V_0$=150 volt]{
\includegraphics[width=0.25\textwidth]{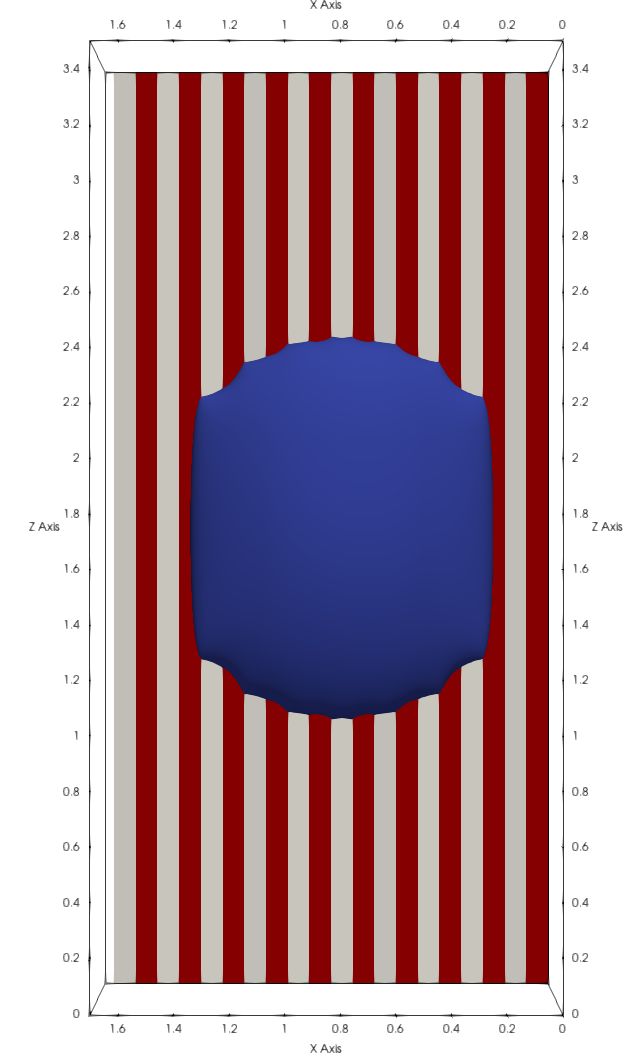}}\qquad\quad
\subfigure[$V_0$=200 volt]{
\includegraphics[width=0.24\textwidth]{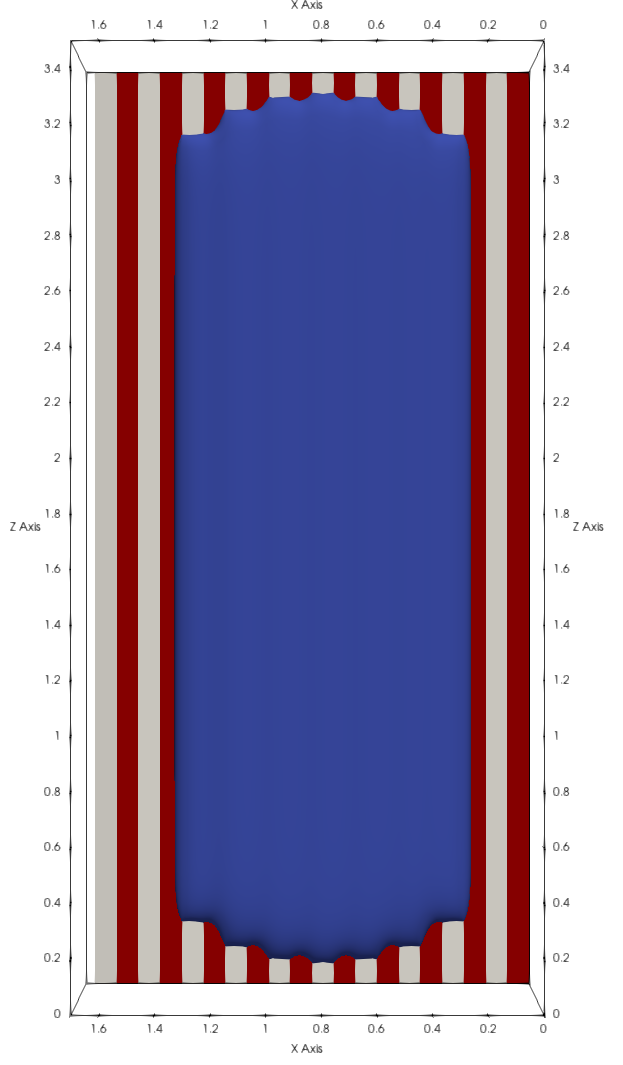}}
}
\centerline{  
\subfigure[$V_0$=100 volt]{
\includegraphics[width=0.3\textwidth]{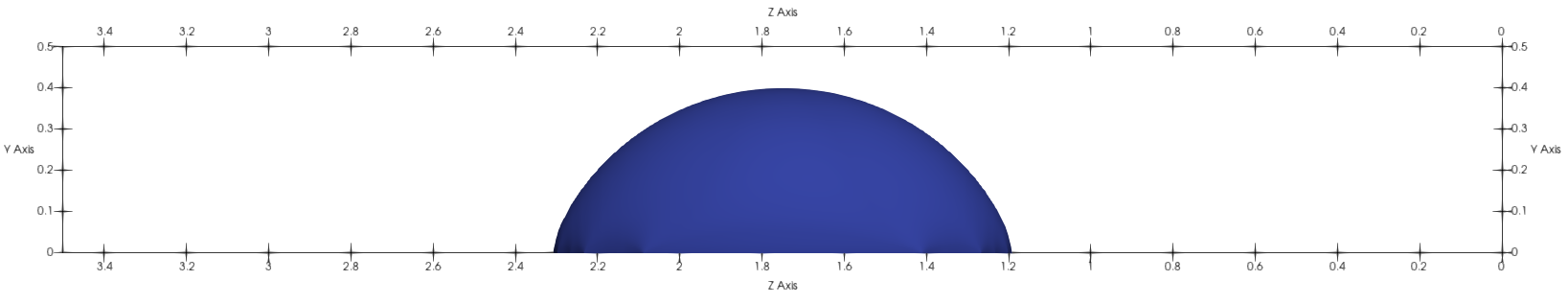}}
\subfigure[$V_0$=150 volt]{
\includegraphics[width=0.3\textwidth]{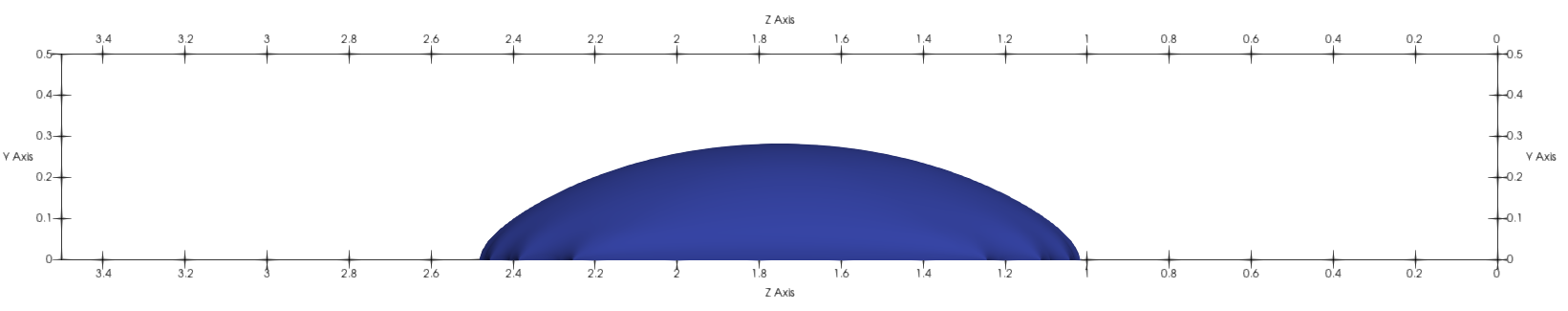}}
\subfigure[$V_0$=200 volt]{
\includegraphics[width=0.3\textwidth]{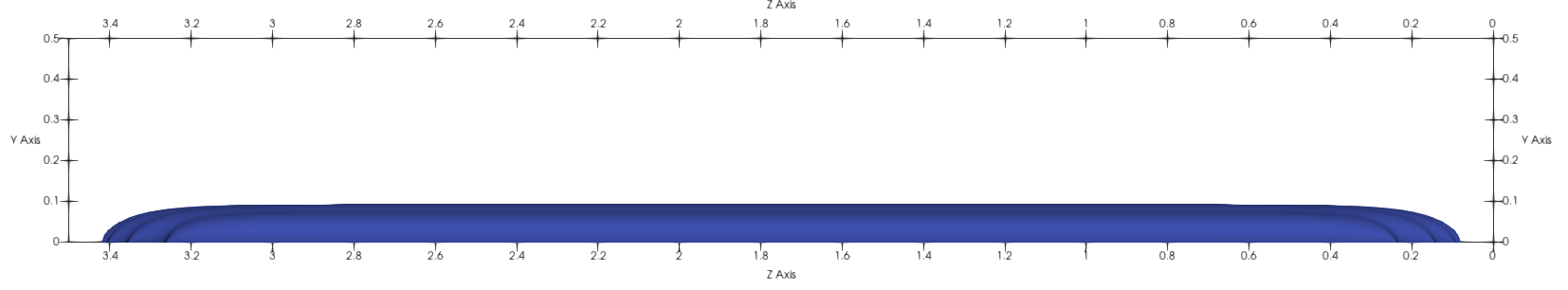}}
}
\centerline{  
\subfigure[$V_0$=100 volt]{
\includegraphics[width=0.3\textwidth]{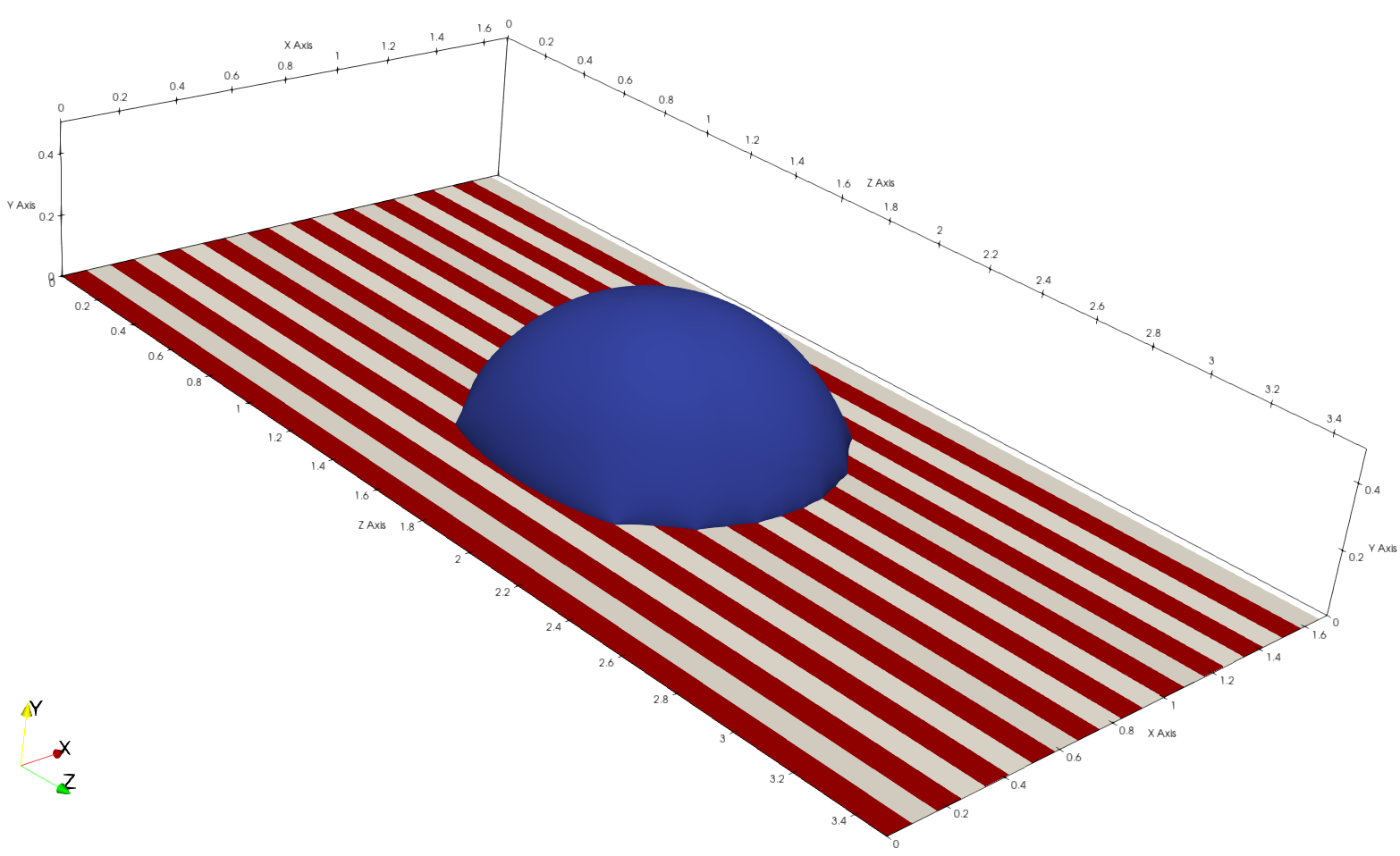}}
\subfigure[$V_0$=150 volt]{
\includegraphics[width=0.3\textwidth]{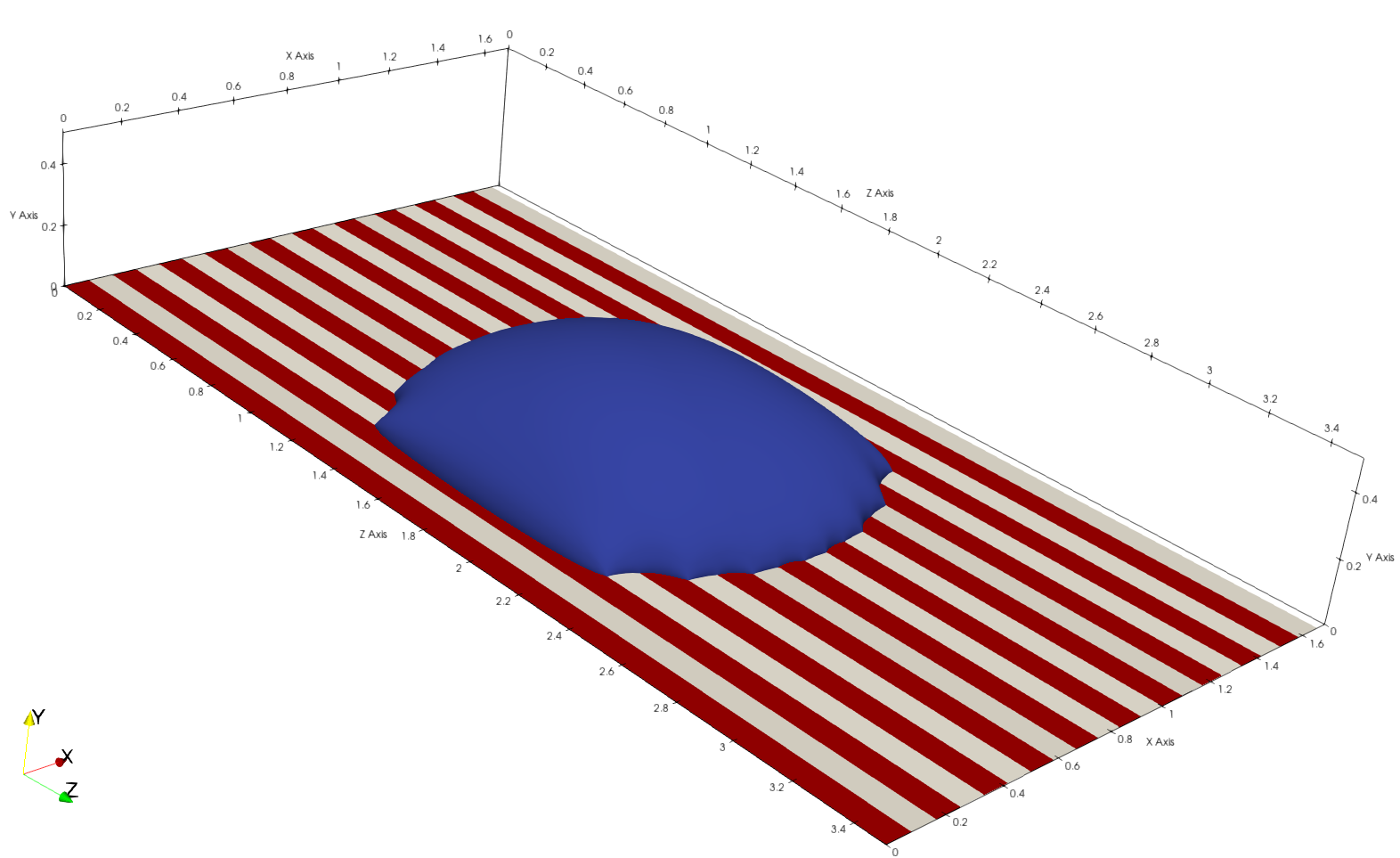}}
\subfigure[$V_0$=200 volt]{
\includegraphics[width=0.3\textwidth]{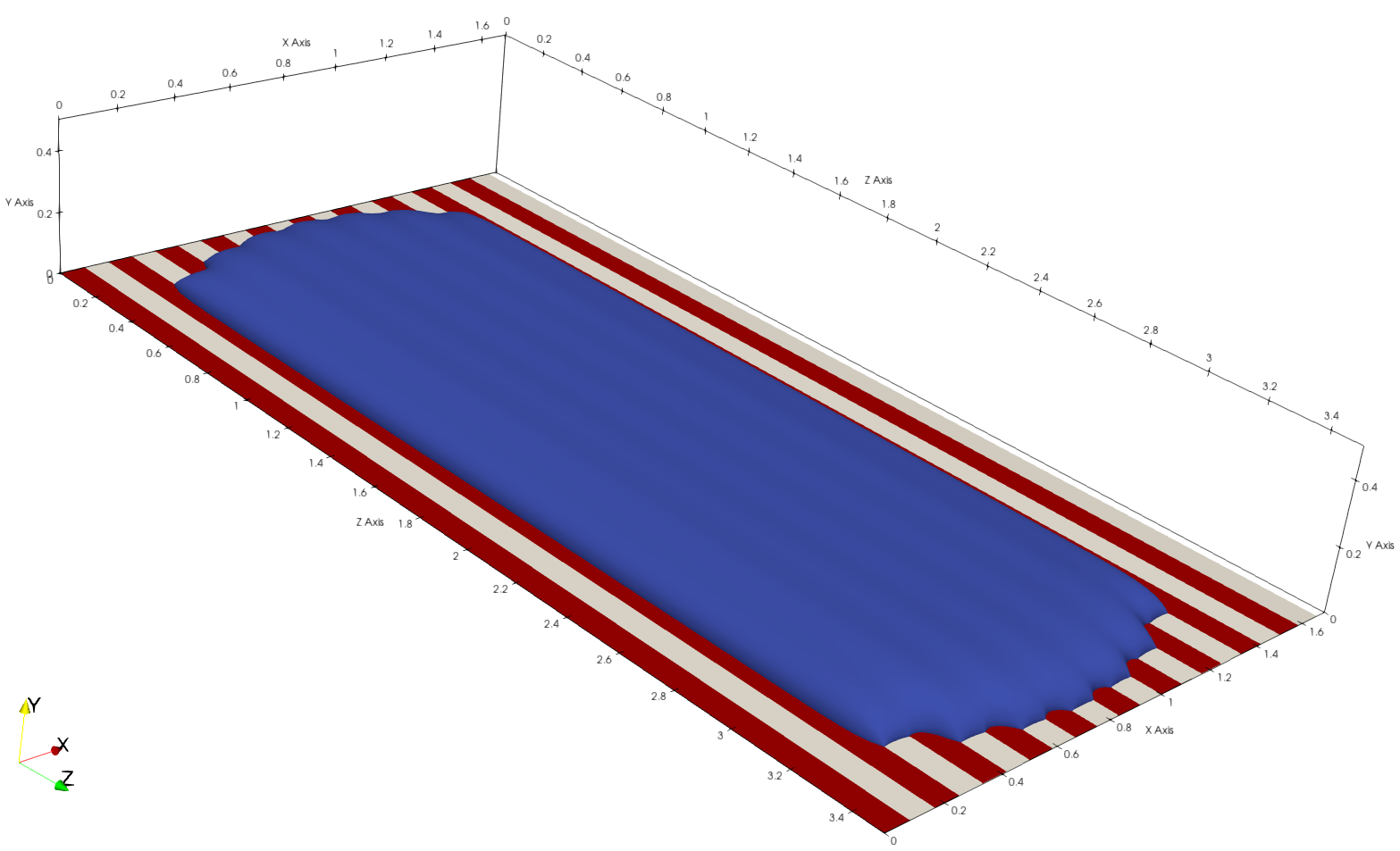}}
}
\caption{\small Dielectric  drop on the wall: 
equilibrium drop shapes under imposed electrode voltage $V_0=100$ volt (left column), $V_0=150$ volt (middle column), and $V_0=200$ volt (right column).
Top row: plan view (toward $-y$ direction);
Middle row: side view (toward $-x$ direction);
Bottom tow: perspective view.
}
\label{fg_5}
\end{figure}

\begin{figure}[tb]
\centering 
\includegraphics[width=0.4\textwidth]{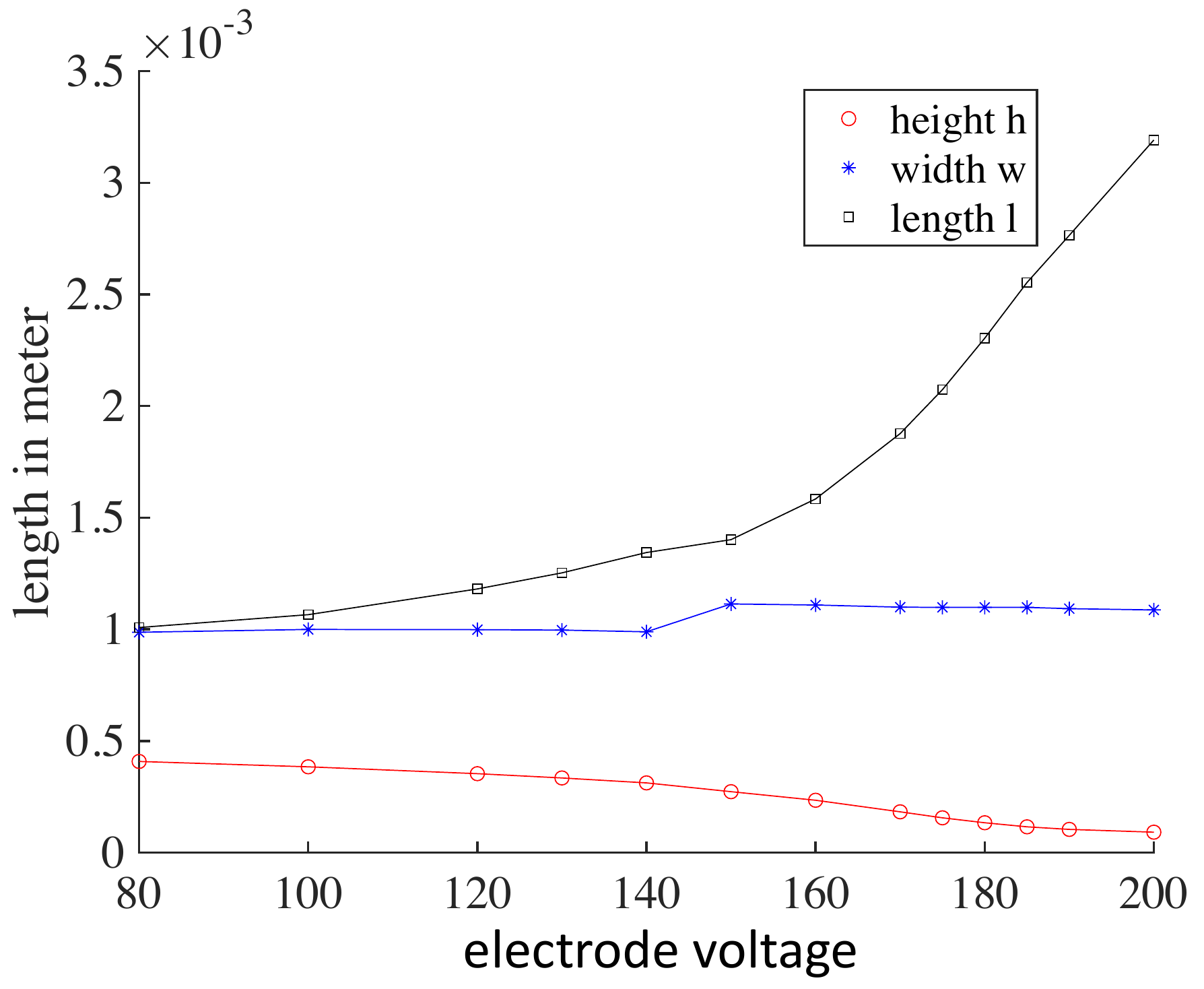}(a)\quad
\includegraphics[width=0.4\textwidth]{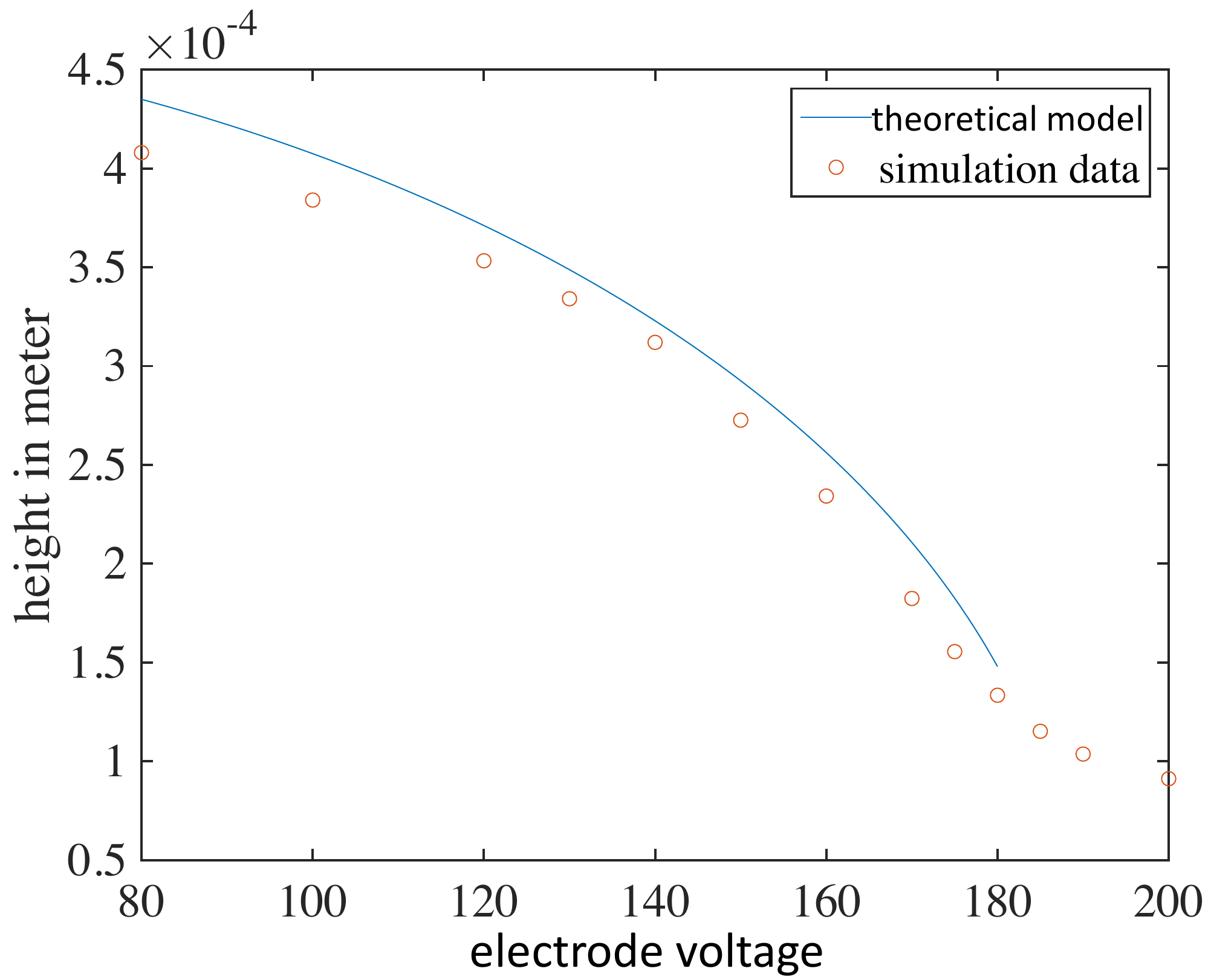}(b)
\caption{\small Dielectric  drop on a wall:
(a) the drop height/width/length as a function of the electrode voltage.
(b) Comparison between the theoretical model and the current simulation on the drop height as a function of the electrode voltage.
} 
\label{dew_hlw} 
\end{figure}

We study the 3D equilibrium shape of a dielectric liquid drop on a horizontal wall under an imposed electric field in this test. The problem setting is in accordance with the  experiment from~\cite{mchale2011dielectrowetting}; see  Figure~\ref{dew}(a).
When the electrodes on the wall are turned on,
the dielectric drop (an initial hemisphere)
deforms due to the imposed electric field, and eventually reaches an equilibrium state, as sketched in Figure~\ref{dew}(b).
We are interested in simulating the equilibrium shape of the dielectric  drop.


As discussed in Section~\ref{sec:steady}, 
the current phase field model allows us to 
compute the equilibrium state of the system 
by solving an alternative simpler system consisting of equations~\eqref{eq_16}, \eqref{eq_15c} and~\eqref{eq_15d}, with the corresponding boundary and initial conditions. After that, if needed, the pressure field can be computed by solving~\eqref{eq_15b}, and the velocity is
given by $\mbs u=0$. We will simulate the equilibrium shape of the dielectric liquid drop by this method.

\begin{table}[tb]
\centering
\begin{tabular}{cc| cc}
     \hline
     variable & normalization constant & variable & normalization constant \\ \hline
    $x,y,z,d,\eta$  & $L_0$ & $t$ & $1$\\
    $V$, $V_0$  & $V_d$  & $\mbs{E}$ & $V_d/L_0 $\\
    $\epsilon,\epsilon_1,\epsilon_2$ &  $L_0\gamma/V_d^2 $ &
    $\gamma_1$  & $L_0^3/\gamma $  \\
    $\lambda$   & $L_0\gamma$ & $\phi,\psi$ & 1\\
     \hline
\end{tabular}
\caption{\small Normalization constants for variables and parameters with the simpler system from Section~\ref{sec:steady} for computing the equilibrium solution. 
Choose $L_0$, $V_d$, and $\gamma$ (surface tension).
}
\label{tab_2}
\end{table}

We consider a computational domain  $(x,y,z)\in\Omega = [0,\frac53 L_0]\times [0,\frac23 L_0]\times [0,\frac72 L_0]$, where $L_0=1.2mm$, as shown in Figure~\ref{dew}(a). The electrodes embedded on the bottom wall each has a width $d=0.1mm$. Adjacent electrodes are $0.1mm$ apart on the wall, and the constant voltage imposed on adjacent electrodes have the same magnitude but with opposite signs ($V_0$ and $-V_0$), as sketched in Figure~\ref{dew}(b). The dielectric liquid drop (in ambient air) is initially shaped like a hemisphere, with a radius $R_0=\frac12 L_0$ and its center located at 
$(X_0,Y_0,Z_0)=(\frac56 L_0,0,\frac74 L_0)$.

We employ the following physical parameter values:
\begin{equation}\small
\left\{
\begin{split}
&
\text{surface tension:}\ \gamma=3.857\times 10^{-2} kg/s^2; \\
&
\text{permittivity: (air)}\ \epsilon_1=\epsilon_0, \quad
\text{(dielectric liquid)}\ \epsilon_2=32\epsilon_0;
\end{split}
\right.
\end{equation}
where $\epsilon_0=8.854\times 10^{-12}F/m$ is the vacuum permittivity. Note that 
the fluid density and viscosity play no role when we simulate the equilibrium state using the system consisting of~\eqref{eq_16}, \eqref{eq_15c} and~\eqref{eq_15d}.

All the dynamic variables and the simulation
parameters have been normalized consistently. The normalization constants used  for non-dimensionalizing the alternative system of equations from Section~\ref{sec:steady} for the equilibrium solution are provided in Table~\ref{tab_2}.
Note that they are a little different from those shown in Table~\ref{tab_1} for normalizing the full system of governing equations.
In particular, all the length variables are normalized by $L_0$. 
For brevity and convenience of presentation, in what follows we employ the same symbols to denote the dimensional and the normalized variables or parameters.
We employ
a Cahn number $\eta=0.02$, and the mobility 
is set by $\lambda\gamma_1=0.1$, where $\lambda=\frac{3}{2\sqrt{2}}\eta$.
The pseudo-time-step size is $\Delta t=2\times 10^{-6}$ in the simulations.

We solve the system consisting of equations~\eqref{eq_16},~\eqref{eq_15c} and~\eqref{eq_15d} by the hybrid spectral element/Fourier spectral method in 3D.
We employ $N_z=120$ Fourier grid points along the $z$ direction and a mesh of $120$ quadrilateral spectral 
elements (with element order $12$) in the $xy$ plane, with $20$ uniform elements along  $x$ and $6$ non-uniform elements along  $y$ (see Figure~\ref{dew}(c)).
We impose the periodic boundary condition in
 $x$  (at $x=0$ and $x=\frac53 L_0$), and the boundary conditions~\eqref{eq_a19} and \eqref{eq_a25} at
the top boundary $y=\frac23 L_0$.
On the bottom wall ($y=0$) we impose
the boundary conditions~\eqref{eq_19} and~\eqref{eq_21}, where the imposed voltage on adjacent electrodes  alternates between $V_0$ and $-V_0$ (see Figure~\ref{dew}(c)).
All the dynamic variables are homogeneous along the $z$ direction.
The initial distribution of the phase field
function is given by 
$\phi(x,y,z)=\tanh\left(\frac{\sqrt{(x-X_0)^2+(y-Y_0)^2+(z-Z_0)^2}-R_0}{\sqrt{2}\eta}\right)$.

Figure~\ref{fg_5} shows the deformed shape of the dielectric drop under three imposed electrode voltages ($V_0=100$volt, $150$volt, and $200$volt) obtained from the 3D simulations. The plots in the three rows show the plan view, the side view, and 
the perspective view of the drop, respectively. The drop deformation becomes increasingly pronounced with increasing electrode voltage.
At $V_0=200$volt, the dielectric drop becomes highly elongated along the $z$
direction (see Figures~\ref{fg_5}(c,f,i)).

Figure~\ref{fg_5} illustrates the asymmetric deformation of dielectric  drops, an important feature observed in experiments (see~\cite{edwards2018dielectrowetting}).
The dielectric droplet tends to stretch along the direction parallel to the electrodes, while in the direction perpendicular to the electrodes the droplet remains approximately the same in dimension. In other words, the width of the drop ($w$ in Figure~\ref{dew}(a))
remains approximately unchanged, while the length and  height of the drop
($l$ and $h$ in Figure~\ref{dew}(a)) can vary significantly with the electrode voltage.

The asymmetric deformation is further demonstrated by Figure~\ref{dew_hlw}(a),
in which we plot the length, width, and height of the deformed dielectric drop as a function of the electrode voltage from our simulations. It is evident that,
while the length and height exhibit a significant change,
the width of the deformed drop remains nearly  constant  as the electrode voltage increases.
This is because the electrodes serve as some potential walls, and so crossing those walls will increase the energy of the system. We refer to~\cite{edwards2018dielectrowetting} 
for more details on the experimental observation and the explanation of the asymmetric deformation.

In \cite{brown2015dielectrophoresis} a theoretical model was proposed on the dielectric drop deformation, and it leads to the following formula relating the drop height to the electrode voltage,
\begin{equation}\small\label{eq_89}
    h^2=h_0^2-\frac{\epsilon_0\Delta \epsilon V_0^2}{4\delta \gamma}\Omega.
\end{equation}
In this equation, $h$ is the deformed drop height, $V_0$ is the electrode voltage, and
$\Omega=h_0l_0$, with $h_0$ and $l_0$ denoting the initial height (in $y$ direction) and initial length (in $x$ direction) of the drop.
$\gamma$ is the surface tension. 
$\epsilon_0$ is the vacuum permittivity, and
$\Delta \epsilon$ is the difference in
the relative permittivity of the two fluids.  $\delta=\frac{4d}{\pi}$ is a geometry parameter.
In Figure \ref{dew_hlw}(b) we show a comparison of the deformed drop height as a function of the electrode voltage between our simulation results and the theoretical model~\eqref{eq_89}.
While there exist some discrepancies in the quantitative values, the simulation results and the model are generally in reasonable agreement in the range $100 \leq V_0 \leq 180$ for the electrode voltage.
It should be noted that the theoretical model~\eqref{eq_89} is only valid for a range of electrode voltage values (when $h\gg \delta$, see~\cite{brown2015dielectrophoresis}).
For the electrode voltage beyond about $180$, the simulation result and the model prediction are qualitatively different. This discrepancy is due to
the breakdown of the model equation~\eqref{eq_89}.
The trend exhibited by the simulation result in this region is similar to what has been observed in the experimental measurement~\cite{brown2015dielectrophoresis}.

\subsection{Equilibrium Dielectric  Film on a Surface}
\label{sec:film}

\begin{figure}[tb]
\centering 
\includegraphics[width=0.45\textwidth]{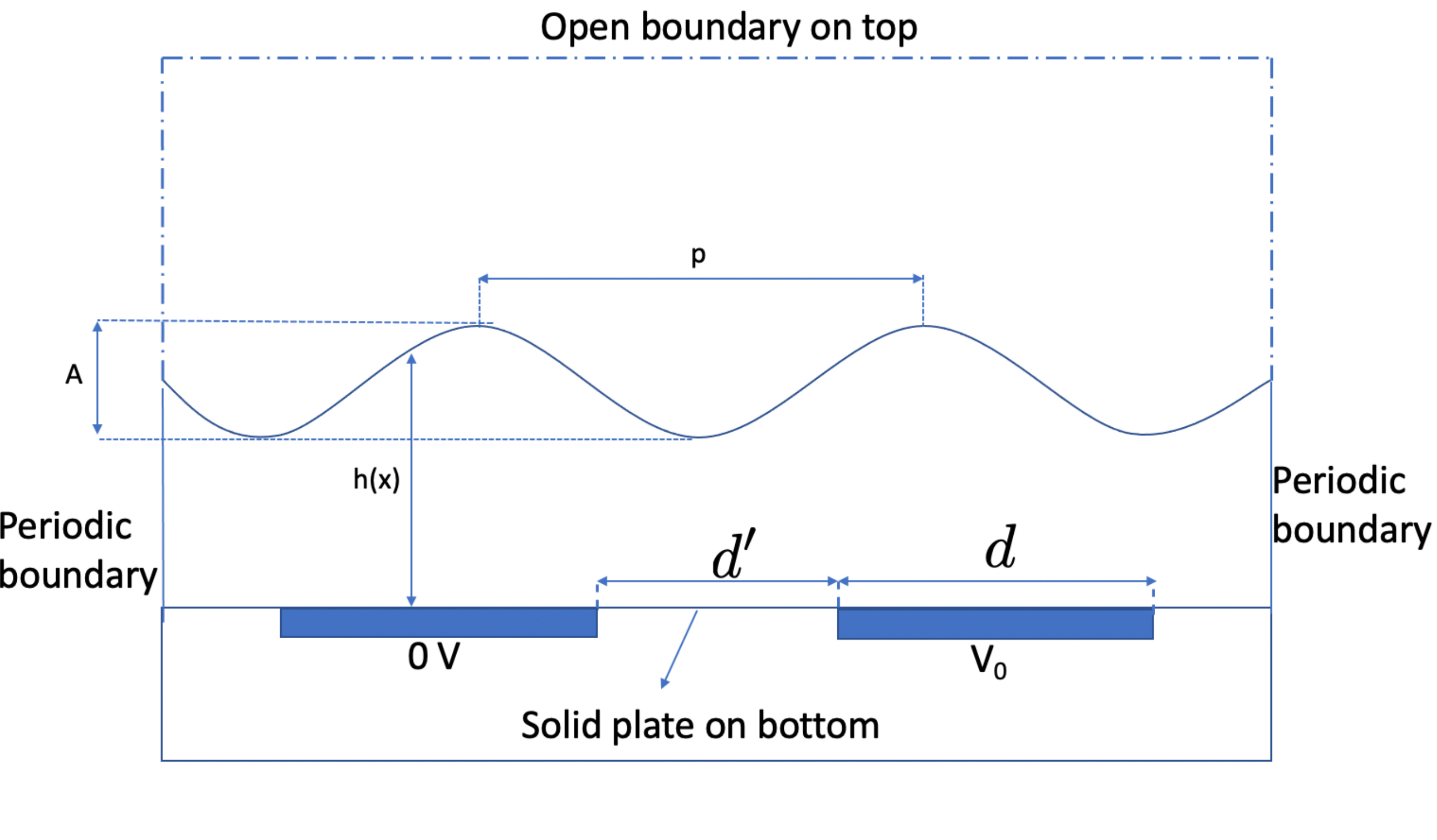}(a)\quad
\includegraphics[width=0.41\textwidth]{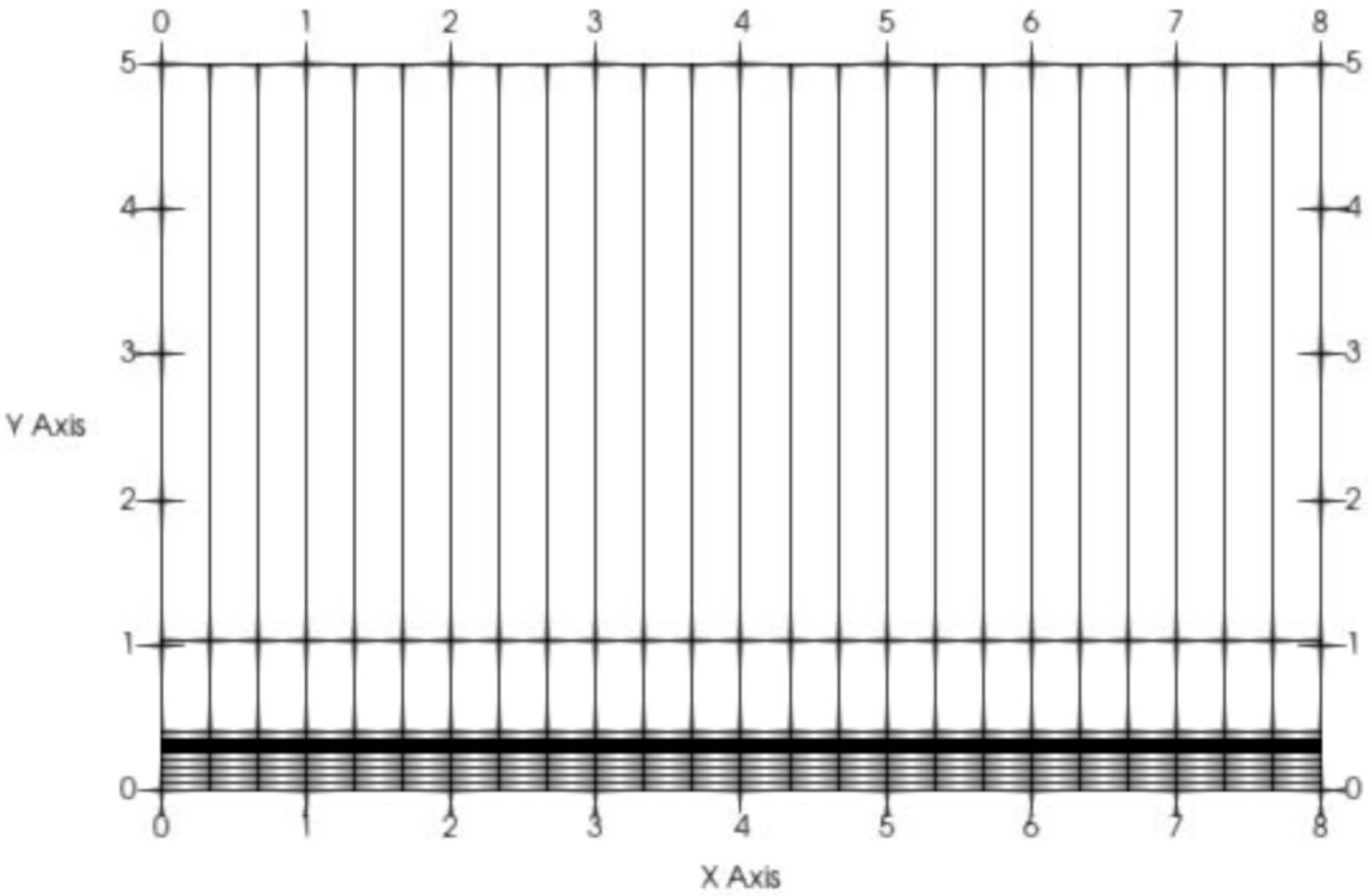}(b)
\caption{\small Dielectric thin film:
(a) flow configuration and settings,
(b) spectral-element mesh. The imposed voltage is $0$ on the left electrode ($1\leqslant x\leqslant 3$) and $V_0$ on the right electrode ($5\leqslant x\leqslant 7$).
} 
\label{optic} 
\end{figure}

In this subsection we study the equilibrium state of a thin dielectric liquid film on a solid wall in two dimensions using the methods developed herein. The dielectric film exhibits a wave-like profile under an imposed  electric field, as observed in the experiment~\cite{brown2009voltage}, in which this is referred to as an optical interface. 

The problem configuration and settings are illustrated in Figure~\ref{optic}(a).
We consider a rectangular domain, $(x,y)\in\Omega=[0,4d]\times [0,\frac52 d]$, where $d$ is the width of the electrode (see below). The domain and all the variables are assumed to be periodic  in the horizontal ($x$) direction.
The top of the domain is open, and the bottom of the domain is a solid wall. Two electrodes, each with a width $d$, are embedded on the bottom wall. The gap between the electrodes is $d'=d$.
The two electrodes specifically occupy the regions $x\in[d/2,3/2]$ and $x\in[5d/2,7d/2]$ on the wall.
The voltage imposed on the right electrode is $V_0$, and on the left electrode is $0$.
A thin layer of dielectric fluid, with a thickness $h_0$, is at rest on the bottom wall in an ambient fluid. When the electrodes are turned on, the fluid interface deforms under the imposed electric field and exhibits a wave-like profile at equilibrium. Our goal is to simulate the equilibrium  dielectric fluid interface.

In what follows we provide two sets of simulations. The first set is obtained using the method from Section~\ref{sec:steady}, based on the simpler system consisting of equations~\eqref{eq_16}, \eqref{eq_15c} and~\eqref{eq_15d}. 
The second set, for comparison, is based on the full model consisting of equations~\eqref{eq_6}--\eqref{eq_a8},
together with appropriate boundary/initial conditions.

\begin{figure}[tb]
\centering 
\includegraphics[width=0.4\textwidth]{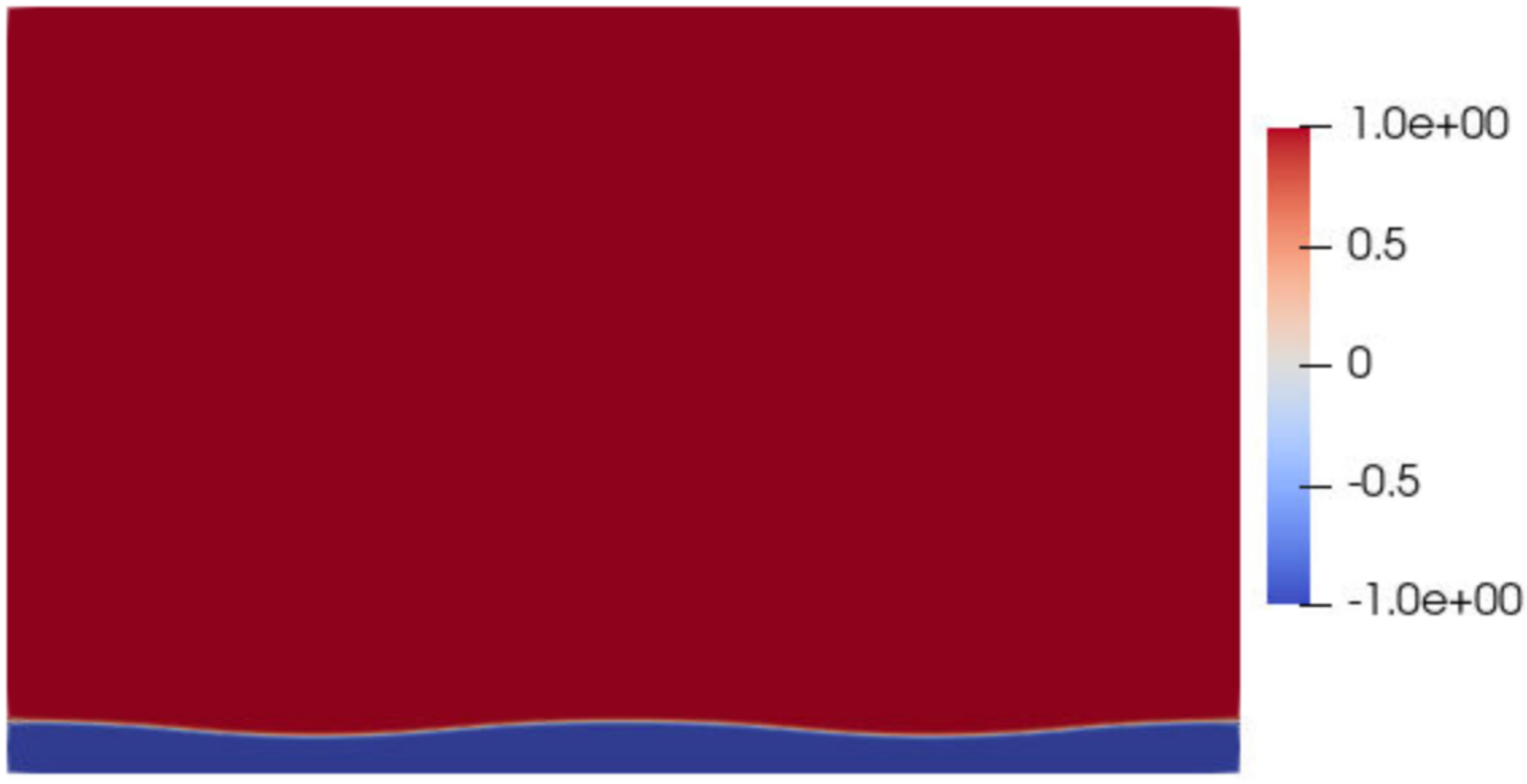}(a)\quad
\includegraphics[width=0.39\textwidth]{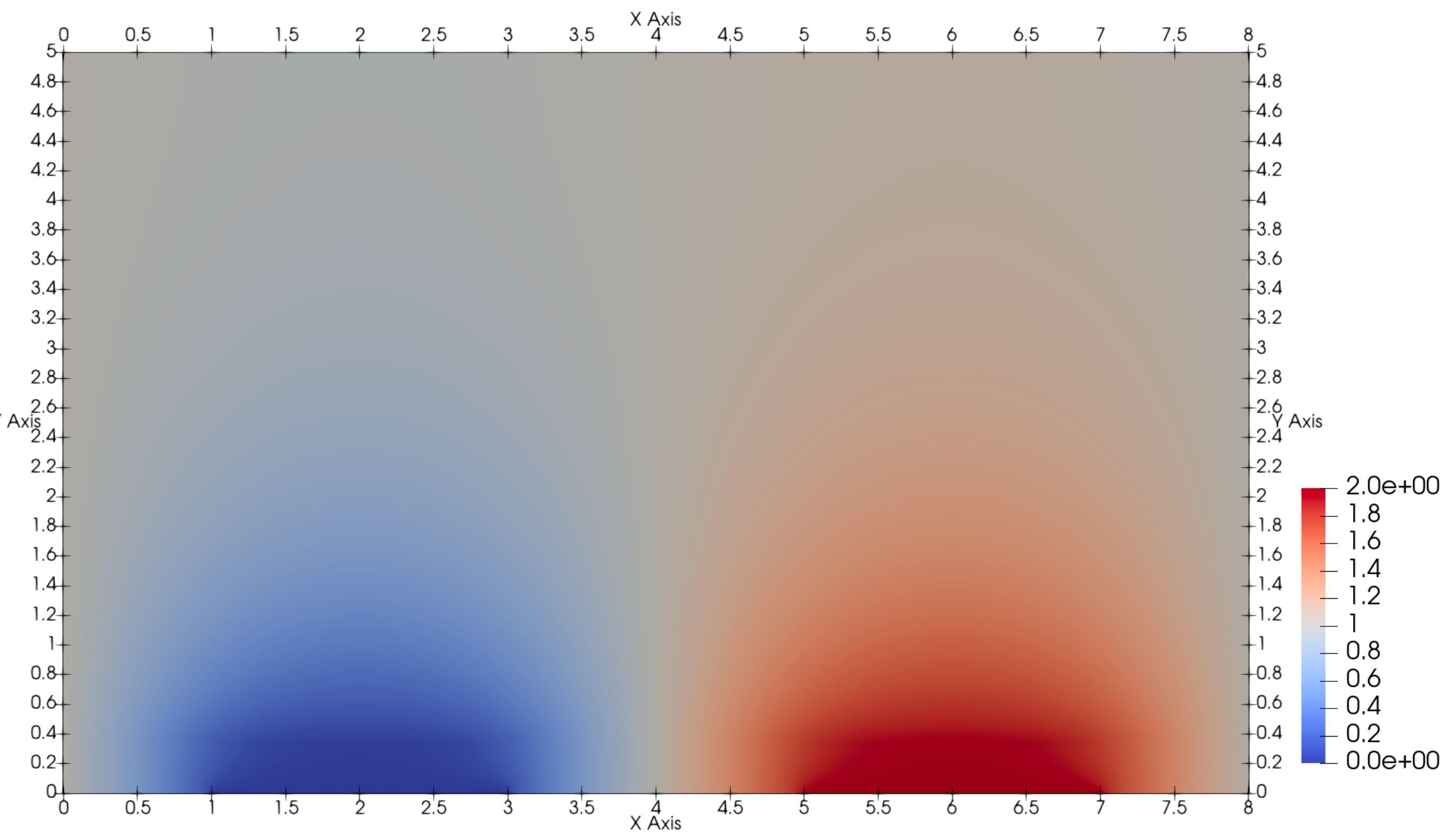}(b)
\caption{\small Dielectric thin film:
Distributions of (a) the phase field function showing the fluid interface, and 
(b) the electric potential in the domain. 
} 
\label{optic_field} 
\end{figure}

\subsubsection{Equilibrium Simulation Using
the Simpler System}
\label{sec:che_film}

We first simulate the equilibrium profile of the dielectric fluid interface using the method from Section~\ref{sec:steady}, by solving the simpler system
of~\eqref{eq_16}, \eqref{eq_15c} and~\eqref{eq_15d}, with the boundary conditions as outlined in the above paragraphs.

We employ a surface tension $\gamma=2.84\times 10^{-2}kg/s^2$, and a permittivity for the ambient fluid the same as the vacuum permittivity, $\epsilon_1=\epsilon_0$.
The permittivity for the dielectric film ($\epsilon_2$) is varied and will be specified below. All the variables and parameters are normalized  based on the normalization constants in Table~\ref{tab_2}.
Here we choose the length scale as $L_0=\frac{d}{2}$, and the voltage scale as $V_d=100$volt.
We use $h(x)$ to denote the thickness of the equilibrium film at $x$.

Figure~\ref{optic}(b) shows a spectral element mesh employed in the current simulations.
The elements are uniform in the $x$ direction, and are generally non-uniform in $y$. 
Along the $y$ direction we divide the domain into three regions: (i) near-wall region ($0\leq y \leq h_0-A/2$), (ii) wave region ($h_0-A/2 \leq y \leq h_0+A/2$), and (iii) upper region ($y\geq h_0+A/2$),
where $A$ is the peak-to-valley amplitude of the wave profile (see Figure~\ref{optic}(a)). 
For setting up the simulations,  the  amplitude $A$ in the above is estimated based on the following theoretical model formula from~\cite{brown2009voltage},
\begin{equation}\small
    A=\frac{16\epsilon_0}{3\pi^4 \gamma}(\epsilon_1-\epsilon_2 )\exp\left(-\frac{2\pi h_0}{p}\right) V_0^2, \label{eq_91}
\end{equation}
where $p=d+d'=2d$.
We employ $N_{y_1}$, $N_{y_2}$ and
$N_{y_3}$ spectral elements in these three regions respectively along the $y$ direction.
The mesh is uniform in the near-wall and wave regions, and is non-uniform in the upper region (Figure~\ref{optic}(b)).
The specific values for $N_{y_1}$, $N_{y_2}$ and $N_{y_3}$ will be provided below when discussing different simulation cases.


In all the simulations we employ a pseudo-time step size $\Delta t=2.0\times 10^{-6}$,
Cahn number $\eta=0.01$, and a mobility $\gamma_1$ by $\lambda\gamma_1=0.1$, where $\lambda=\frac{3}{2\sqrt{2}}\eta$.
The initial phase field distribution is
\begin{equation}\small\label{eq_a92}
\phi(x,y)=\tanh\left(\frac{y-h_0}{\sqrt{2}\eta}\right).
\end{equation}
It should be noted that, while the physical length scale may be different for different simulation cases,
the normalized computational domain is fixed
due to the choice $L_0=d/2$ and is always
$(x,y)\in\Omega=[0,8]\times [0,5]$.


\begin{figure}
\centerline{
\includegraphics[width=0.45\textwidth]{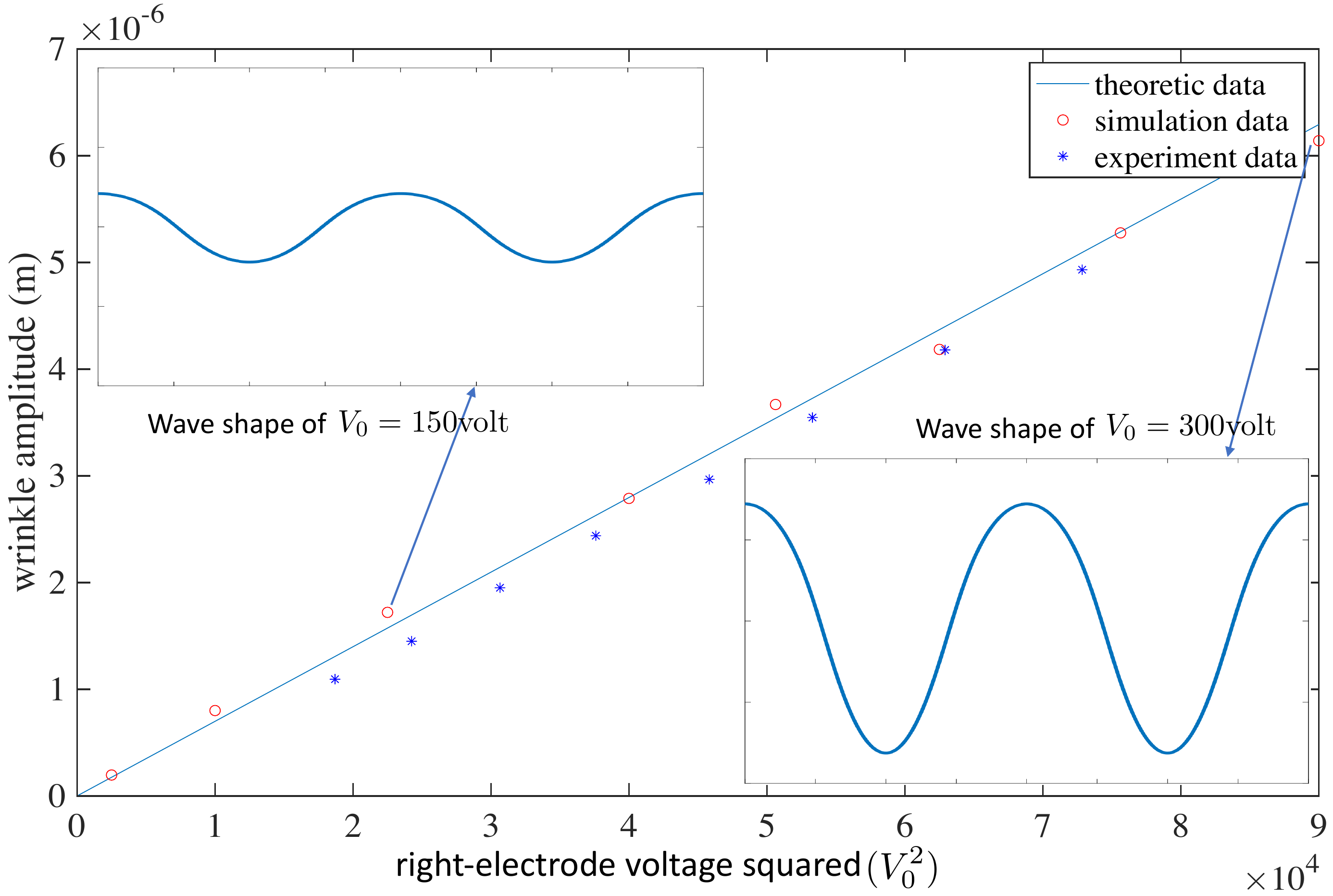}(a)
\includegraphics[width=0.45\textwidth]{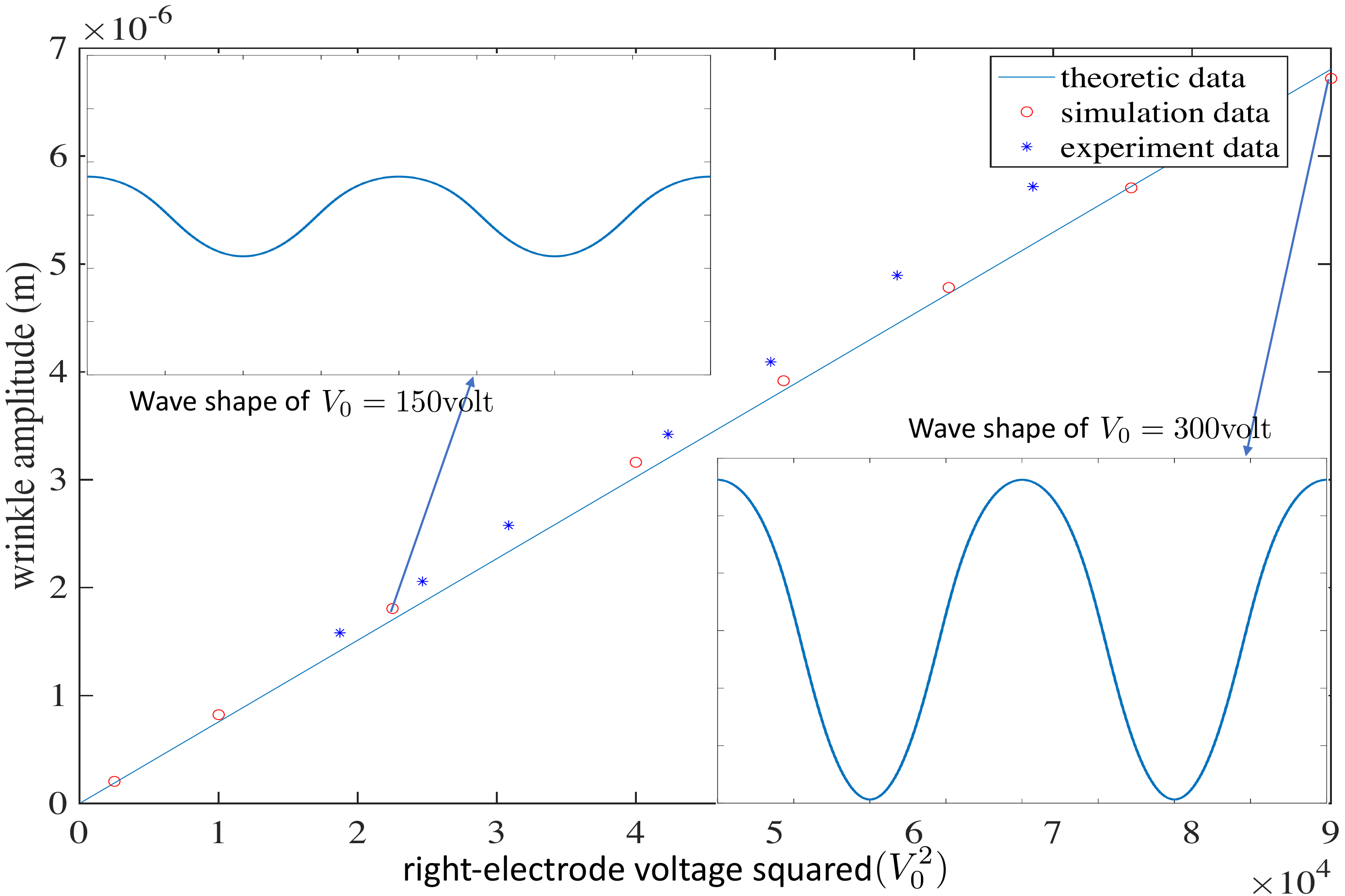}(b)
}
\caption{\small Dielectric thin film: Comparison of the amplitude ($A$) as a function of the electrode voltage squared ($V_0^2$) from the current simulations, the theoretical model (equation~\eqref{eq_91}), and the experimental measurement~\cite{brown2009voltage}, for two cases with (a) $h_0=14\mu m$ and $p=160\mu m$, and
(b) $h_0=18\mu m$ and $p=240\mu m$. The insets of these plots show two typical interface profiles. 
}
\label{fg_9} 
\end{figure}

\begin{figure}
\centering 
\includegraphics[width=0.45\textwidth]{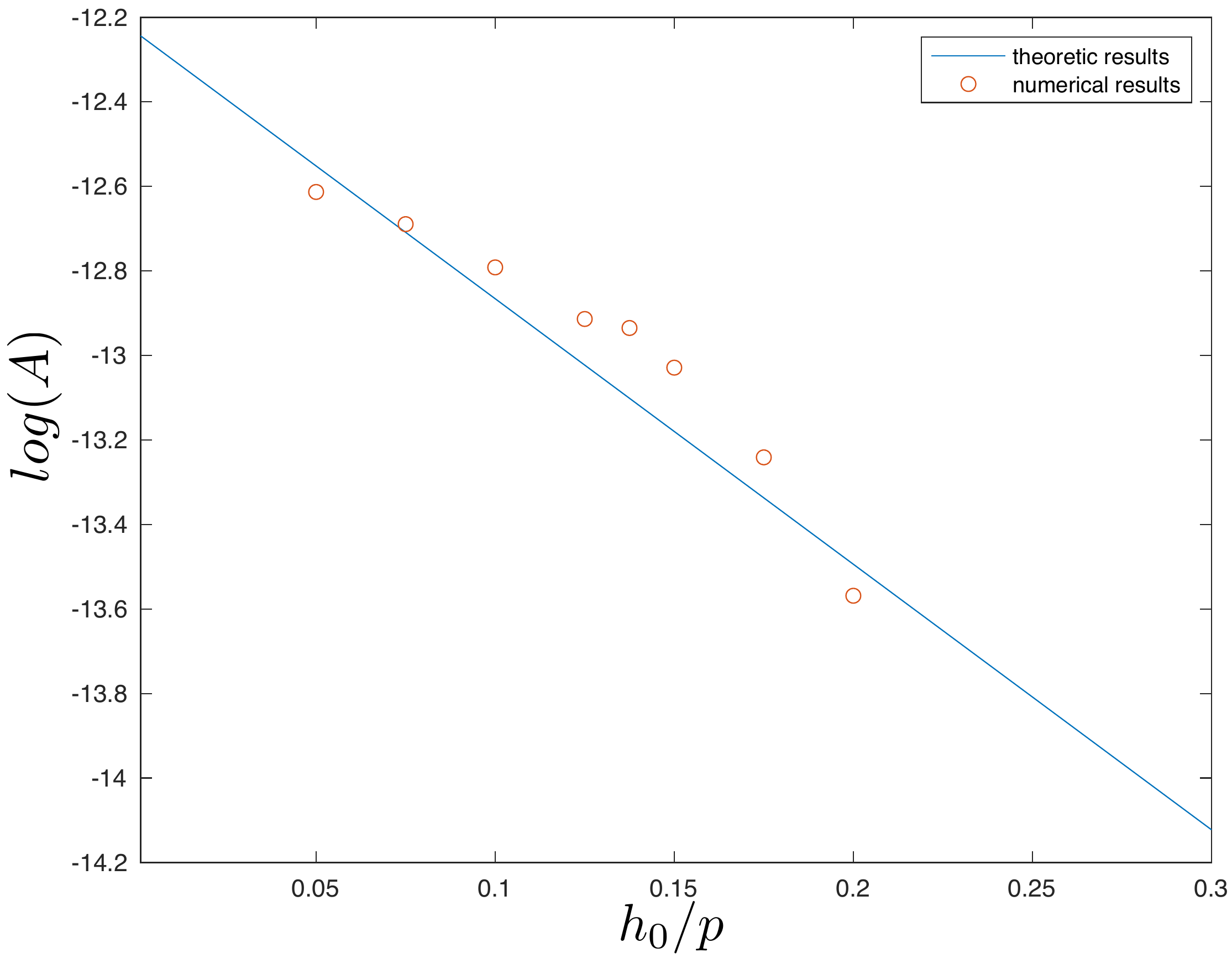}
\caption{\small Dielectric thin film: Comparison of $\log(A)$ (interfacial wave amplitude) versus $h_0/p$ (initial film thickness) from the current simulations and the theoretical model equation~\eqref{eq_91}. 
}
\label{fg_10} 
\end{figure}

\begin{figure}
\centerline{ 
\includegraphics[width=0.45\textwidth]{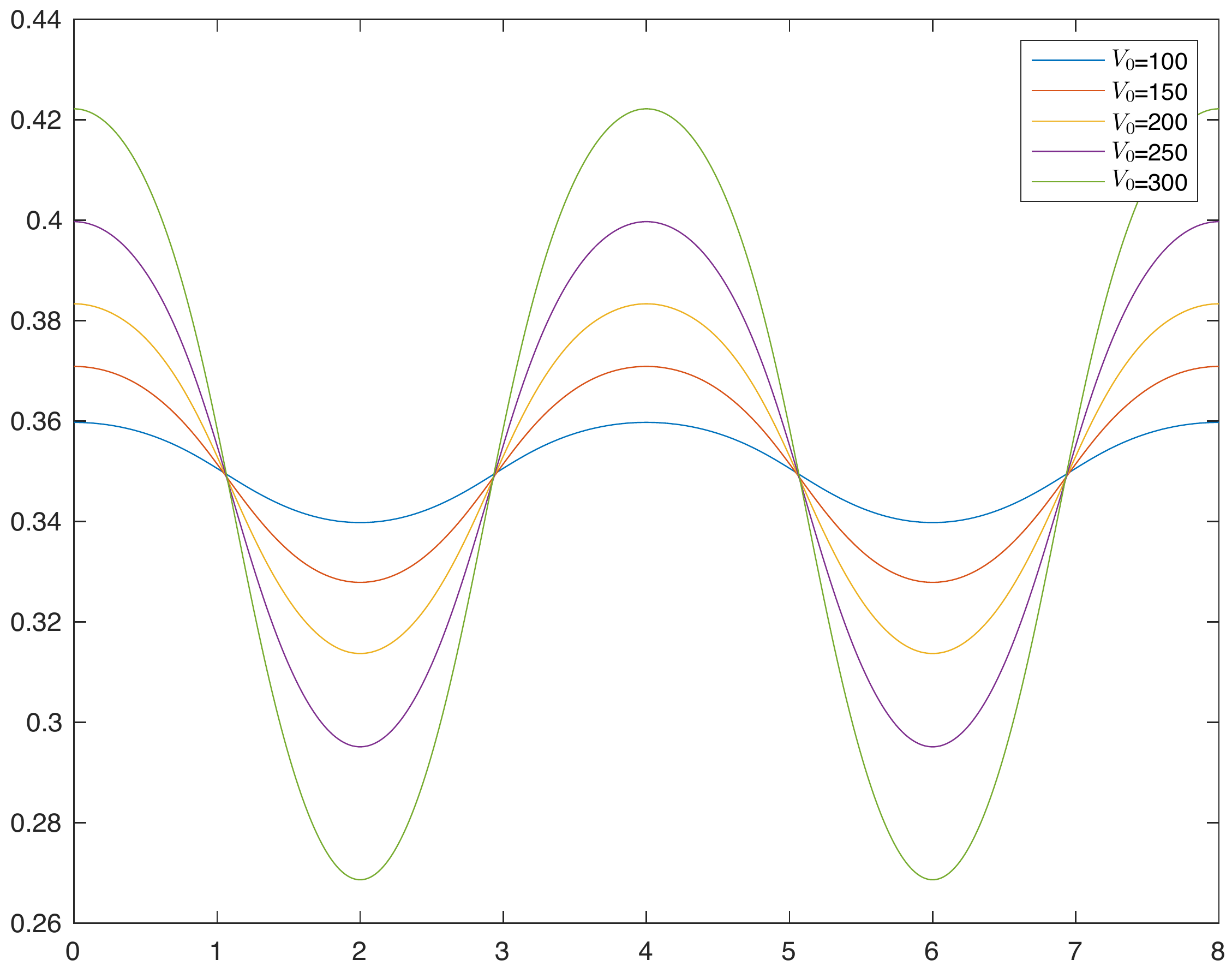}(a)
\includegraphics[width=0.45\textwidth]{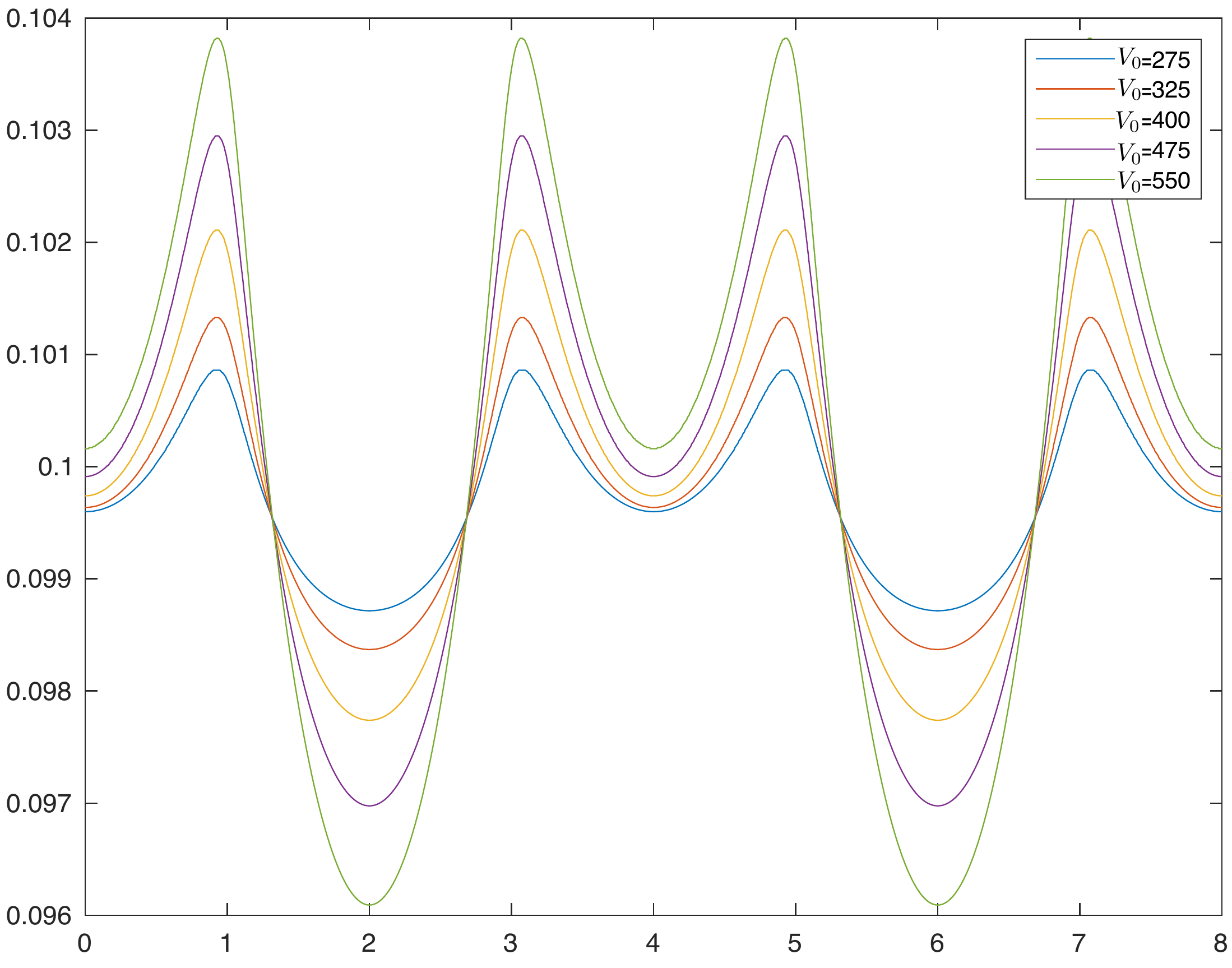}(b)
}
\caption{\small Dielectric thin film: Interfacial  profiles at different electrode voltages: (a) sinusoidal-like profiles, (b) non-sinusoidal like profiles.
In (a), $h_0=14\mu m$, $p=160\mu m$, and
$\frac{\epsilon_2}{\epsilon_1}=8$.
In (b), $h_0=6\mu m$, $p=240\mu m$, and
$\frac{\epsilon_2}{\epsilon_1}=2$.
}
\label{fg_11} 
\end{figure}


Figure~\ref{optic_field} provides an overview of the equilibrium distributions of
the phase field function $\phi$ (plot (a)) and the electric potential field $V$ (plot (b)). 
The wavy fluid interface
is unmistakable from Figure~\ref{optic_field}(a).
This figure also illustrates that the domain dimension in $y$ is much larger than the dielectric film thickness in our simulations.
This is necessary  because in the upper open boundary we have used the boundary condition~\eqref{eq_a25},
which is  accurate only when  the height of the computational domain is sufficiently large compared with the size of the electrode.

Figure~\ref{fg_9} is a comparison of the 
equilibrium interfacial amplitudes obtained 
from the current simulations, the theoretical model formula~\eqref{eq_91},
and the experimental measurement of~\cite{brown2009voltage}.
Here the initial film thickness ($h_0$) and
the electrode width ($d$) are fixed, while the voltage on the right electrode ($V_0$) is varied systematically. The permittivity ratio is $\frac{\epsilon_2}{\epsilon_1}=8$.
We employ $(N_{y_1},N_{y_2},N_{y_3})=(5,4,4)$ elements along the $y$ direction in this set of simulations, with an element order $12$ for all the elements.
The two plots in this figure show 
the equilibrium interfacial amplitude as a function of  $V_0^2$ for two cases,
corresponding to $h_0=14\mu m$ and $p=160\mu m$ (Figure~\ref{fg_9}(a)) and 
$h_0=18\mu m$ and $p=240\mu m$ (Figure~\ref{fg_9}(b)), respectively.
The insets of these plots depict two typical  interfacial profiles at equilibrium
corresponding to $V_0=150$ and $300$ volts.
It can be observed that the simulation results agree with the theoretical model
and with the experimental data reasonably well.

Figure~\ref{fg_10} shows another comparison between the current simulation and the theoretical model~\eqref{eq_91}.
In this set of simulations we have a fixed $V_0=200$volt, $p=160\mu m$ and $\epsilon_2/\epsilon_1=8$, while the initial thickness of the film is varied systematically. We again employ $(N_{y_1},N_{y_2},N_{y_3})=(5,4,4)$ elements along the $y$ direction.
This figure plots the $\log(A)$ as a function of $h_0/p$ from these tests.
While there are some discrepancies, the simulation results overall are close to the predictions of the theoretical model equation~\eqref{eq_91}.

As observed in the experiments of~\cite{brown2009voltage} and in the boundary integral model of~\cite{chappell2020numerical},
the interfacial profiles that are  sinusoidal-like  or non-sinusoidal-like can occur under the imposed electric field. We have observed both types of profiles in our simulations.
Figure~\ref{fg_11} shows examples of these two types of interfacial profiles attained from our simulations, corresponding to several electrode voltage values.
Figure~\ref{fg_11}(a) corresponds to the case in Figure~\ref{fg_9}(a) (with $h_0=14\mu m$, $p=160\mu m$ and $\epsilon_2/\epsilon_1=8$), exhibiting a sinusoidal wave-like profile.
Figure~\ref{fg_11}(b) corresponds to the parameter values $h_0=6\mu m$, $p=240\mu m$
and $\epsilon_2/\epsilon_1=2$, exhibiting an apparently non-sinusoidal wave-like profile.

\subsubsection{Comparison with  Full-Model Simulation}

\begin{figure}
\centering  
\subfigure[$V_0=150$volt]{
\includegraphics[width=0.4\textwidth]{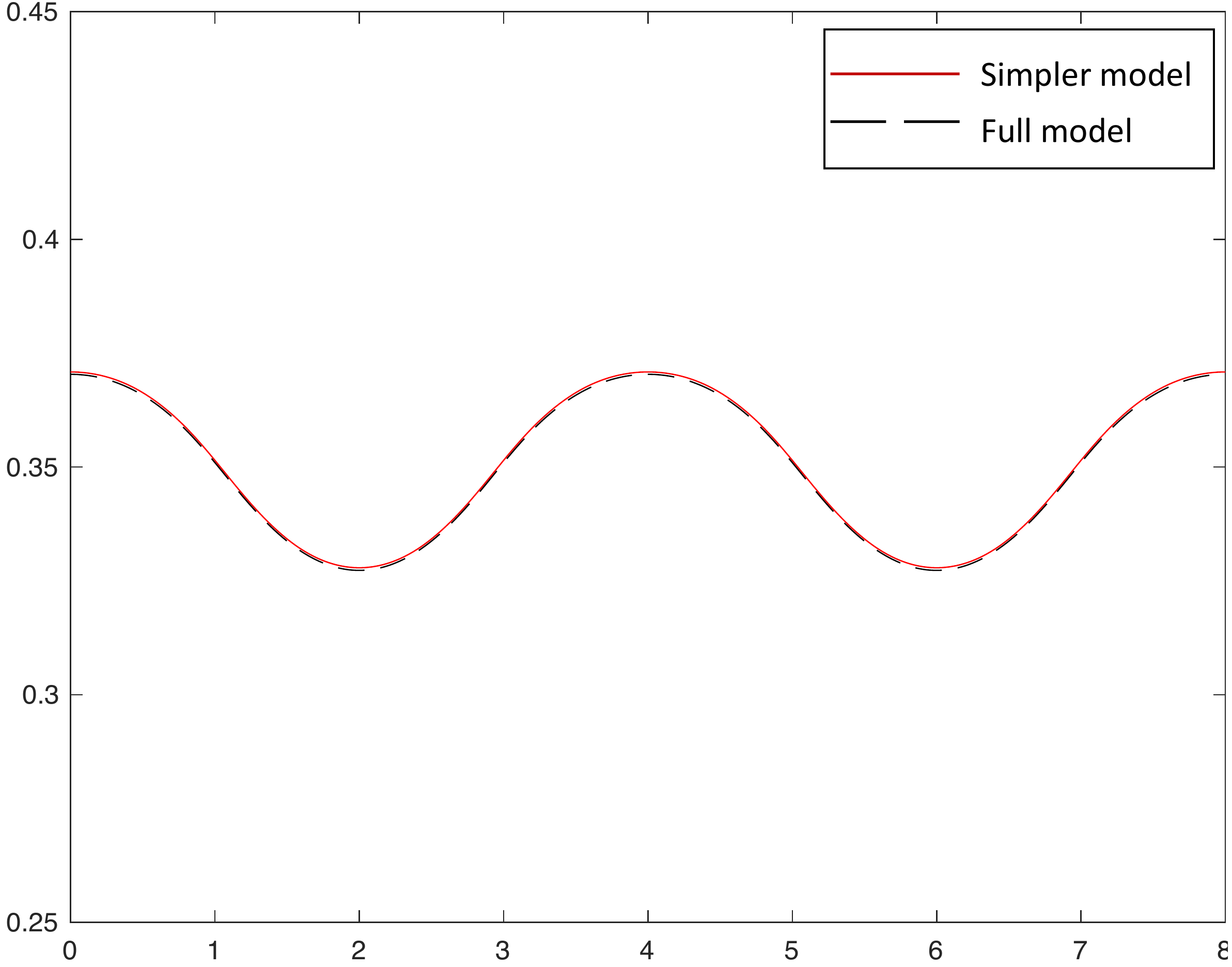}}\qquad
\subfigure[$V_0=200$volt]{
\includegraphics[width=0.4\textwidth]{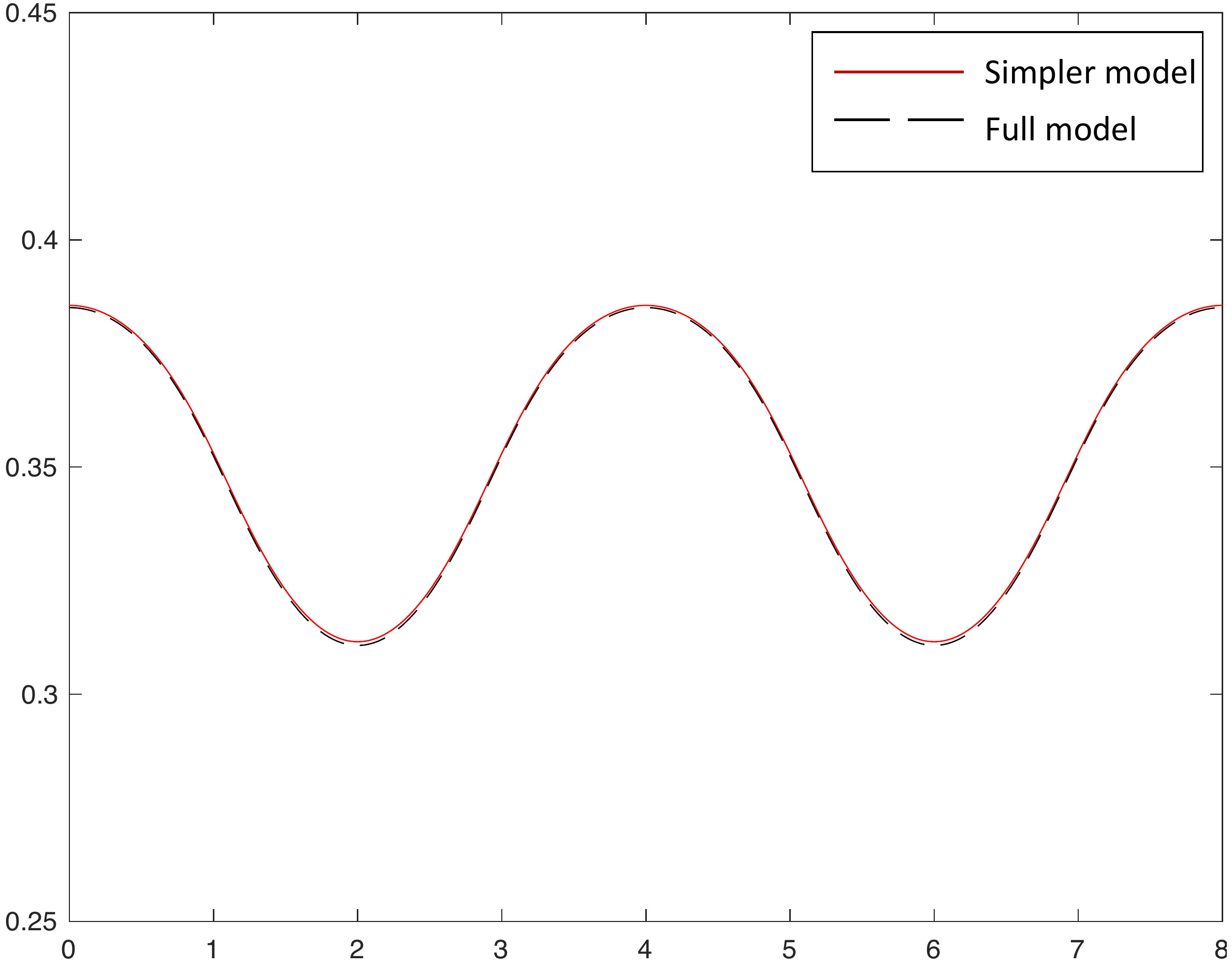}}
\subfigure[$V_0=225$volt]{
\includegraphics[width=0.4\textwidth]{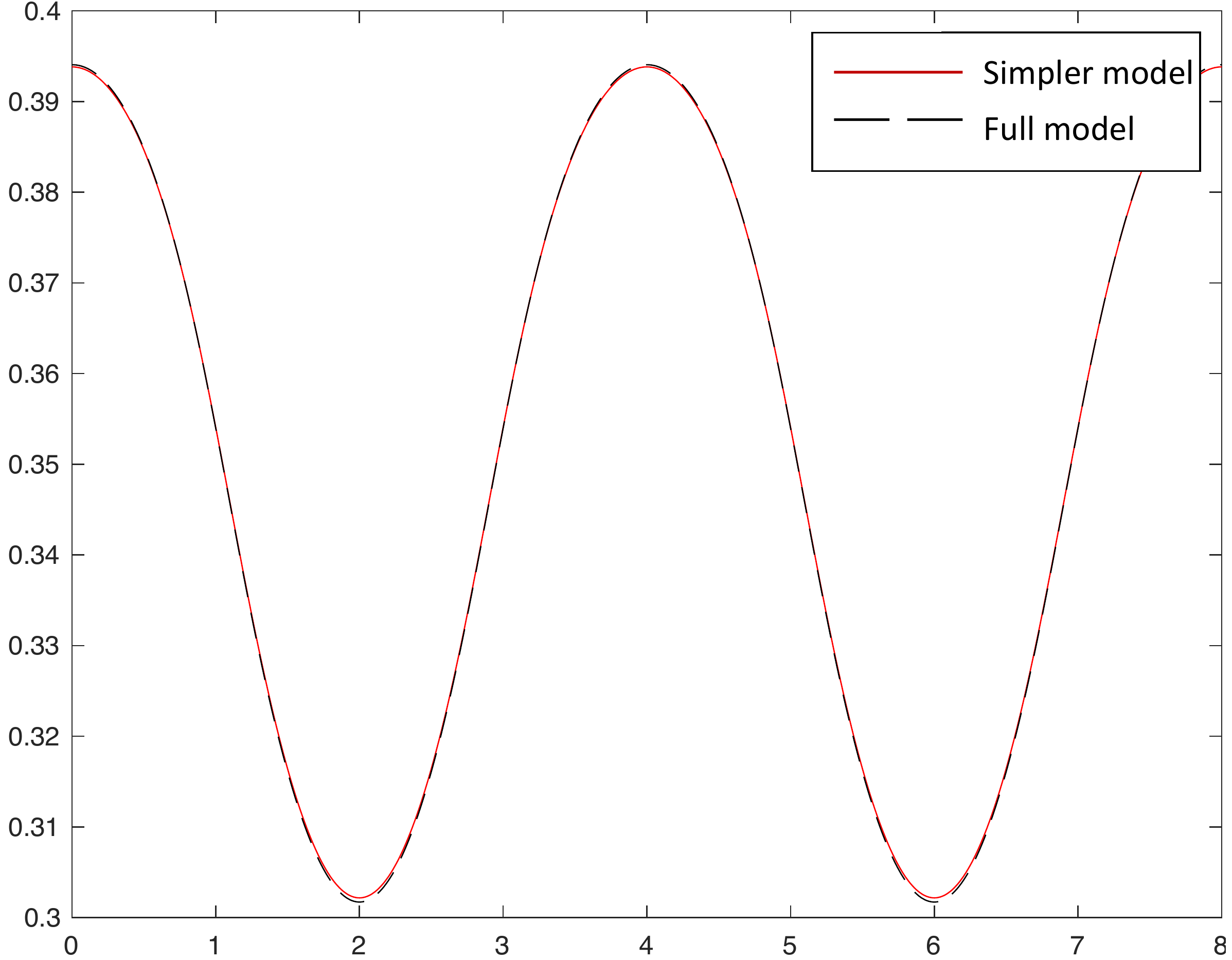}}\qquad
\subfigure[$V_0=275$volt]{
\includegraphics[width=0.4\textwidth]{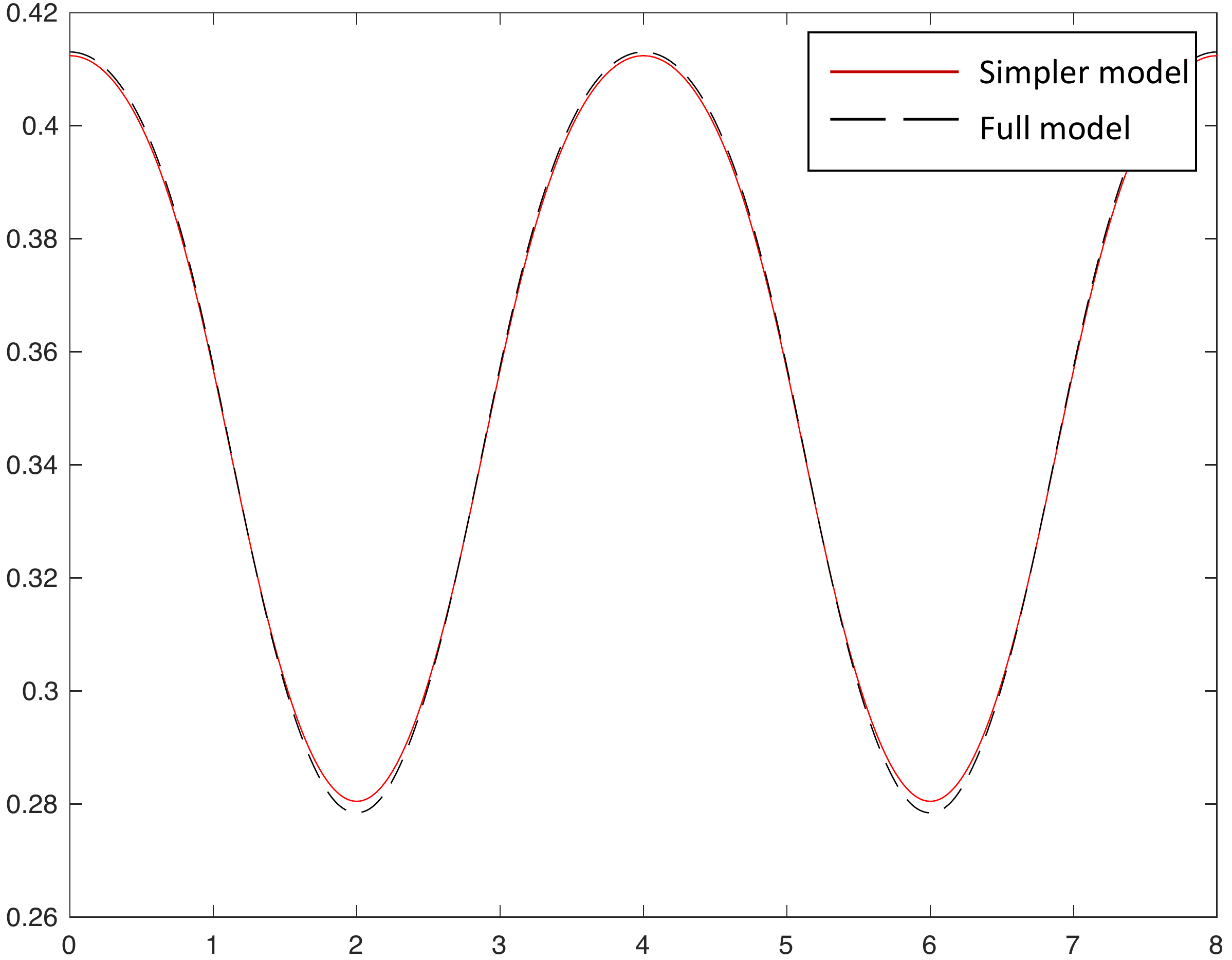}}
\caption{\small Dielectric thin film: Comparison of equilibrium interfacial profiles at several electrode voltages obtained from equilibrium simulations based on the simpler model of Section~\ref{sec:steady} and based on the full model.
}
\label{fg_12}
\end{figure}

\begin{figure}
\centering 
\includegraphics[width=0.45\textwidth]{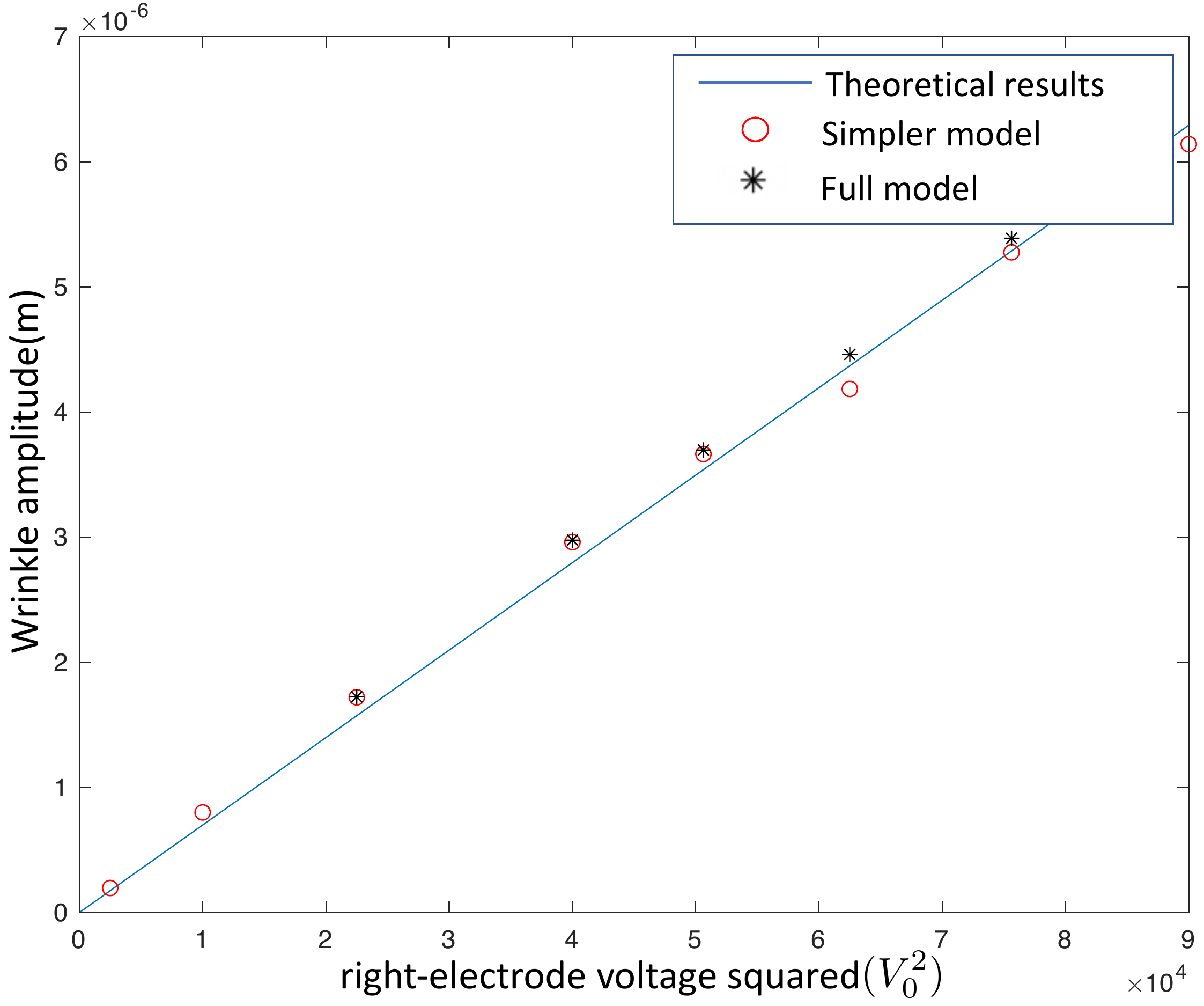}
\caption{\small Dielectric thin film: Comparison of the interfacial amplitude versus the electrode voltage squared obtained from the equilibrium simulations based on the simpler model and the full model, and from the theoretical model~\eqref{eq_91}. 
}
\label{fg_13} 
\end{figure}

We now simulate the equilibrium profile of the dielectric fluid
interface using the full model as given by the Equations~\eqref{eq_6}--\eqref{eq_a8}, together with the boundary and initial conditions.
The flow configuration and the problem setting follow those of Section~\ref{sec:che_film}, as given in Figure~\ref{optic}. 

We consider the same group of tests as in Figure~\ref{fg_9}(a).
The values for the physical and geometric parameters, such as the surface tension $\gamma$, the permittivities ($\epsilon_1$ and $\epsilon_2$), $h_0$ and $p$, are taken to be the same as in Section~\ref{sec:che_film} (specifically Figure~\ref{fg_9}(a)).
The only difference lies in the fluid densities and the
dynamic viscosities, which are needed in the full model  but do not appear in the simpler model of Section~\ref{sec:che_film}. 
Here in the full model we employ $\rho_1=\rho_2=830 kg/m^3$ for the two densities, and
$\mu_1=1.2048\times 10^{-5}kg/(m\cdot s)$
and $\mu_2=2\mu_1$ for the two dynamic viscosities.
Employing the same density for the two fluids  apparently does not correspond to realistic situations. Since we are seeking the equilibrium solution, employing the same density in principle will not alter the solution at equilibrium, but will make the computation considerably easier.
All the physical variables and parameters have been normalized consistently.

In the full-model simulations, we employ the following simulation parameter values (non-dimensional): Cahn number $\eta=0.01$, mobility $\gamma_1=0.05$, $\Delta t=2\times 10^{-6}$,
the number of elements in the three regions along $y$ 
$(N_{y_1}, N_{y_2}, N_{y_3})=(5,4,4)$, and an element order $12$. The initial phase field profile is given by~\eqref{eq_a92}. The electrode voltage $V_0$
is varied in the tests.
The simulations have been performed for a sufficiently long time until the velocity becomes very  small.


Figure~\ref{fg_12} shows a comparison of the equilibrium interfacial profiles obtained by the simpler model of Section~\ref{sec:che_film} and by the full model here. 
These profiles correspond to several elctrode voltages ranging from $V_0=150$volt to $V_0=275$volt.
The results from the simpler model and the full model in general agree very well, with their profiles essentially overlapping with each other.
At larger electrode voltages (e.g.~$V_0=275$volt), some discrepancy in the valley (or peak) of the interfacial profile can be noticed between these two models.

Figure~\ref{fg_13} is another comparison between the simpler model and the full model. It shows the interfacial amplitude $A$ (see Figure~\ref{optic}(a))
as a function of the electrode voltage squared ($V_0^2$) obtained from the simpler model, the full model, and the theoretical model~\eqref{eq_91}.
It can be observed that the results from the simpler model and the full model agree well with each other, and that both are in good agreement with the theoretical model~\eqref{eq_91}.

\subsection{Dynamic Simulations}
\label{sec:transport}

\begin{figure}[tb]
\centerline{
\includegraphics[width=0.48\textwidth]{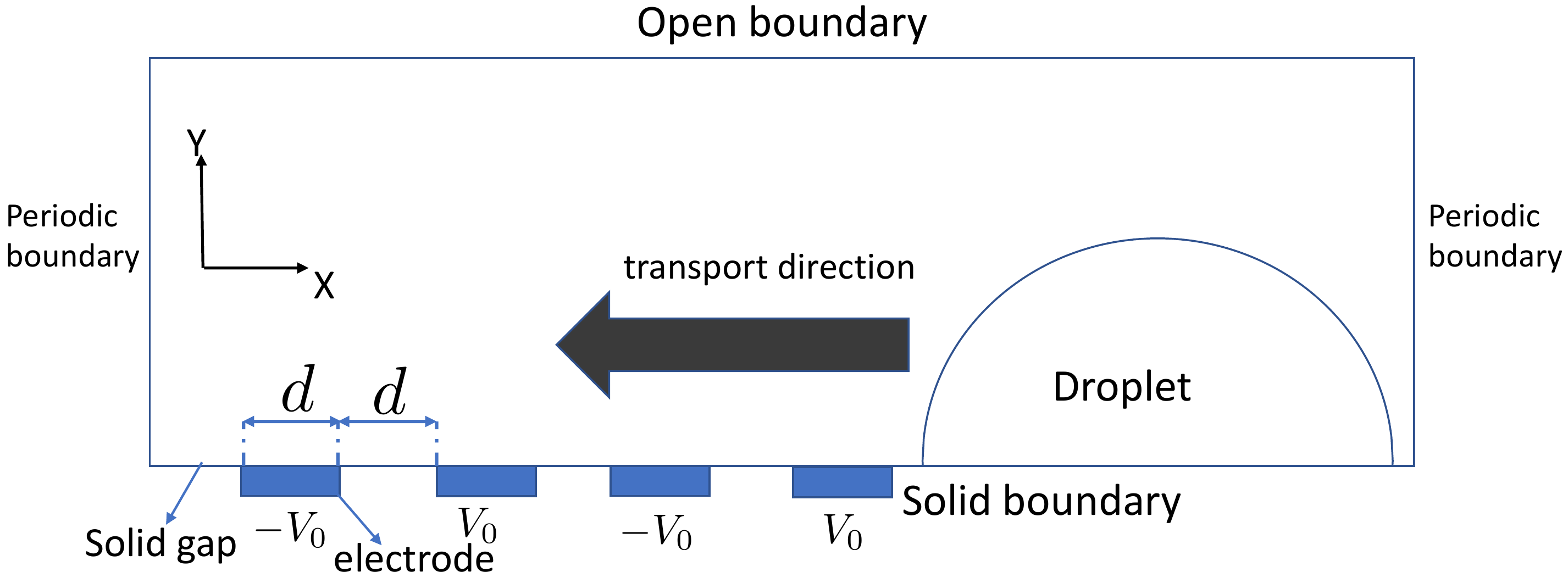}(a)
\includegraphics[width=0.48\textwidth]{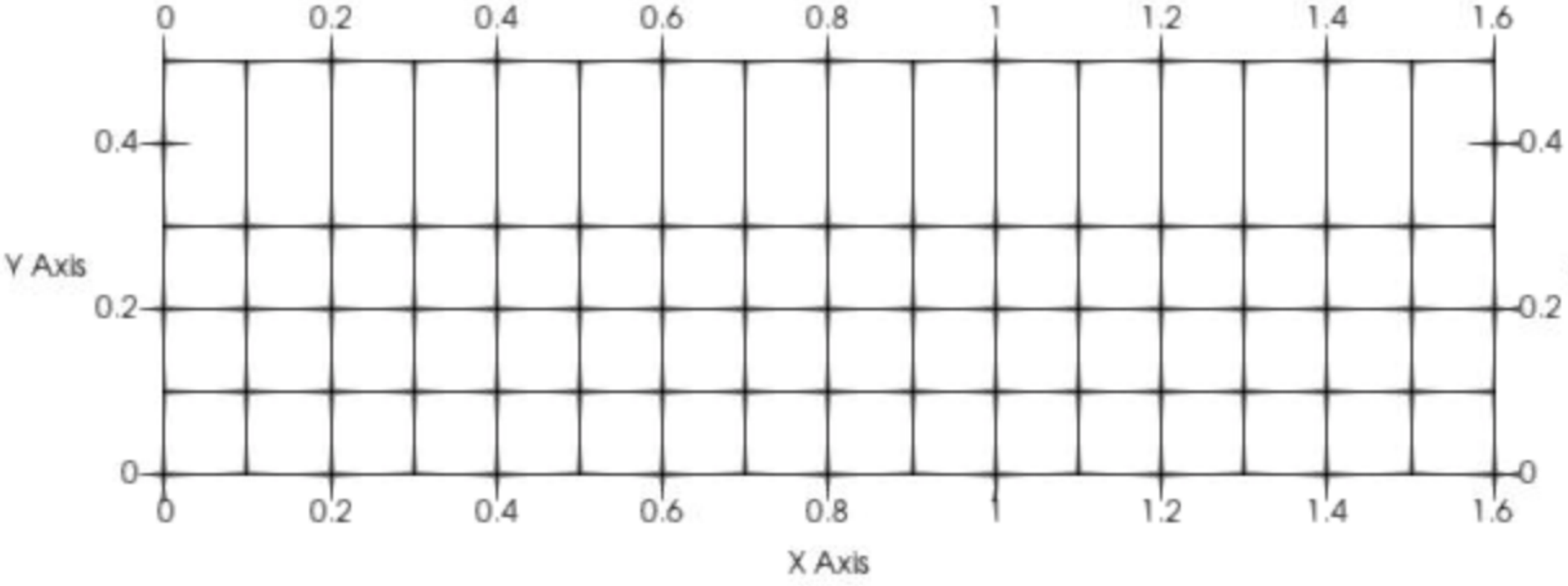}(b)
}
\caption{\small Drop transport:
(a) Flow configuration and settings.
(b) Mesh of $64$ spectral elements. 
} 
\label{fg_14} 
\end{figure}

The simulations in Sections~\ref{sec:drop} and~\ref{sec:film} are for steady-state problems. In this section we  further test the proposed method using  dynamic problems with  two-phase dielectric flows.

\subsubsection{Transport of a Dielectric  Drop on a Wall}
\label{sec:transdrop}

We look into the transport of a dielectric fluid drop on a horizontal wall in two dimensions in this subsection. The problem setting is illustrated by Figure~\ref{fg_14}(a). We consider a rectangular domain, which is periodic in the horizontal direction, open on the top, and has a solid wall at the bottom. An array of  electrodes is embedded  on the left half of the bottom wall, while the right half of the wall is free of electrodes. 
A dielectric liquid drop is initially at rest in the electrode-free region of the  wall. When the electrodes are switched on, the drop will be pulled leftwards to the electrode-embedded region of the wall, due to its interaction with the nonuniform electric field. 
The goal of this problem is to simulate the motion of the liquid drop.

\begin{figure}[tb]
\centerline{
\includegraphics[width=0.6\textwidth]{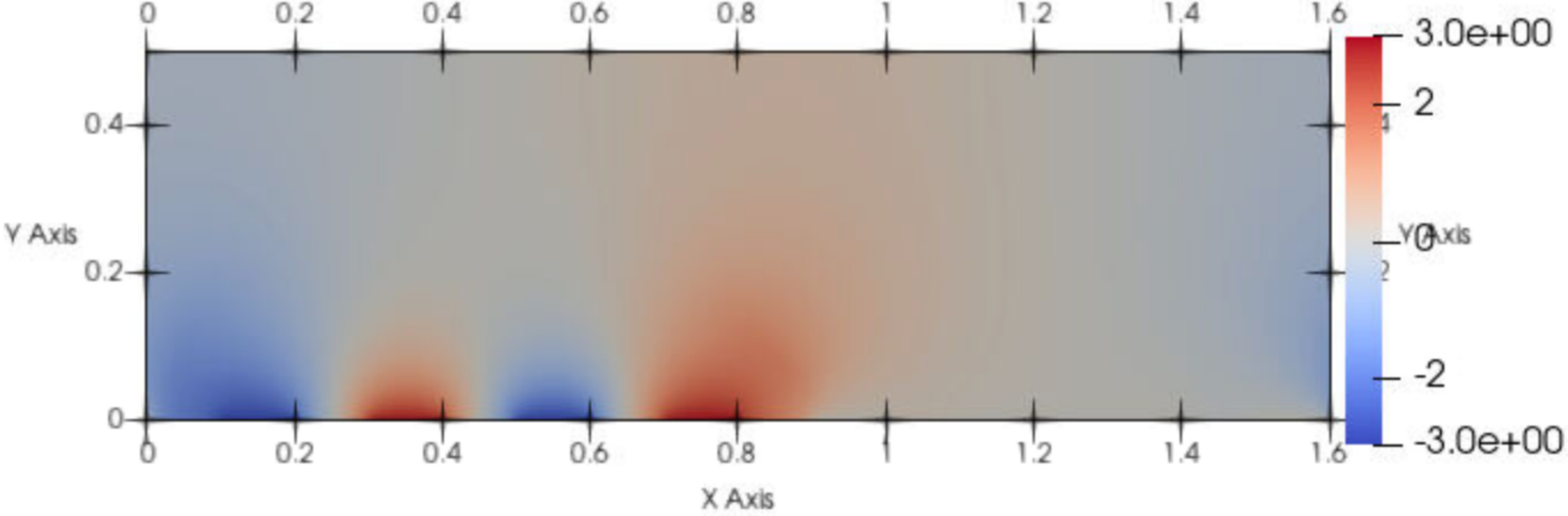}
}
\caption{\small Drop transport: Distribution of the
electric potential ($t=0.002$).
} 
\label{fg_15} 
\end{figure}


We employ the model given by equations~\eqref{eq_6}--\eqref{eq_a8} to simulate this problem, with the boundary and initial conditions as described below.
We consider a computational domain, $(x,y)\in\Omega=[0,\frac85 L_0]\times[0,\frac12 L_0]$, where $L_0=1mm$.
Figure~\ref{fg_14}(b) shows the mesh of $64$ quadrilateral spectral elements employed in the simulations, with $16$ and $4$ elements along the $x$ and $y$ directions, respectively.
Four electrodes are embedded on the bottom wall, with a voltage $V_0$ or $-V_0$, where $V_0=300$volt.
The electrode-embedded regions on the wall are:
$x/L_0\in[0.1,0.2]$ (voltage $-V_0$),
$x/L_0\in[0.3,0.4]$ (voltage $V_0$),
$x/L_0\in[0.5,0.6]$ (voltage $-V_0$),
and $x/L_0\in[0.7,0.8]$ (voltage $V_0$).
%
We impose the boundary conditions~\eqref{eq_19}--\eqref{eq_21} on the bottom wall ($y/L_0=0$), with  a static contact angle $\theta_s=90^0$. The boundary conditions~\eqref{eq_a19}--\eqref{eq_a25} are imposed on the top domain boundary ($y/L_0=0.5$).
We impose periodic boundary conditions for all the dynamic variables on the horizontal boundaries ($x/L_0=0, 1.6$).
%
The drop is assumed to be semi-circular initially, with a radius $R_0=0.35L_0$ and its center located at $(x_0,y_0)=(1.2L_0,0)$.
We employ an initial phase field profile,
\begin{equation}\small
\phi(x,y,t=0)=\tanh\frac{\sqrt{(x-x_0)^2+(y-y_0)^2}-R_0}{\sqrt{2}\eta},
\end{equation}
and zero initial velocity in the simulations.


The following physical parameters are employed for this problem:
\begin{equation}\small
\left\{
\begin{split}
&
\text{surface tension:}\ \gamma=2.84\times 10^{-2}kg/s^2; \\
&
\text{densities:}\ \rho_1=\rho_2=429.7kg/m^3;
\ (\text{ambient fluid}\ \rho_1,\ \text{drop}\ \rho_2) \\
&
\text{dynamic viscosities:}\
\text{(ambient fluid)}\ \mu_1=12.048\times 10^{-4}kg/(m\cdot s), \quad
\text{(drop)}\ \mu_2=2\mu_1; \\
&
\text{permittivities:}\ 
\text{(ambient fluid)}\ \epsilon_1=\epsilon_0=8.854\times 10^{-12}F/m, \quad
\text{(drop)}\ \epsilon_2=8.1\epsilon_0; \\
&
\text{model parameters:}\
\eta = 0.01L_0, \quad
\lambda = \frac{3}{2\sqrt{2}}\gamma\eta, \quad
\gamma_1 = 5\times 10^{-6}\frac{L_0^2}{\mu_1}, \quad
\Delta t=1\times 10^{-6}\frac{L_0^2\mu_1}{\epsilon_0 V_0^2}.
\end{split}
\right.
\end{equation}
All the physical variables and  parameters are normalized consistently based on the normalization constants listed in Table~\ref{tab_1},
with $L_0$ and $\epsilon_0$ as given above and the $\mu_0$ and $V_d$ therein given by
$\mu_0=\mu_1$ and $V_d=V_0$ for this problem.
We employ an element order $12$ in the simulations.

\begin{figure}[tb]
\centerline{  
\subfigure[$t=0.002$]{
\includegraphics[width=0.24\textwidth]{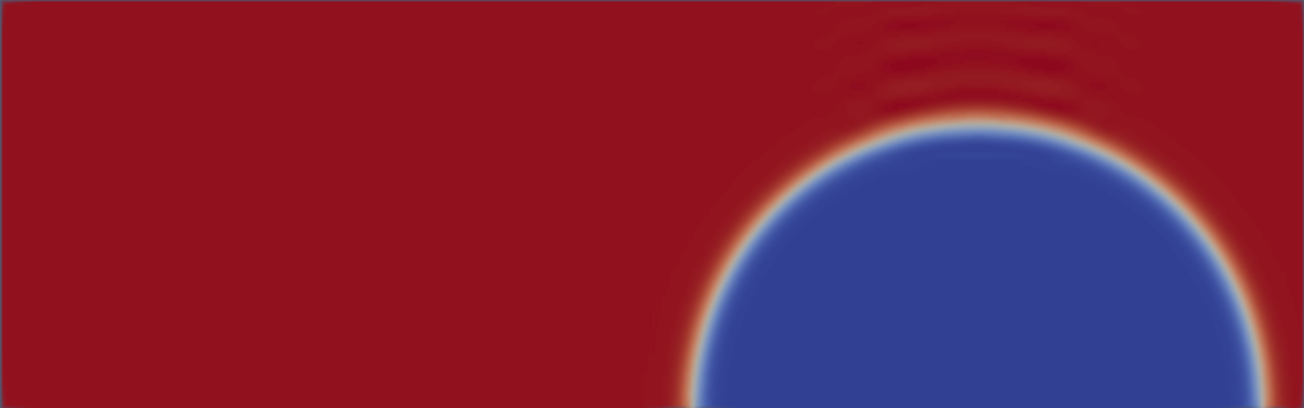}}
\subfigure[$t=0.09$]{
\includegraphics[width=0.24\textwidth]{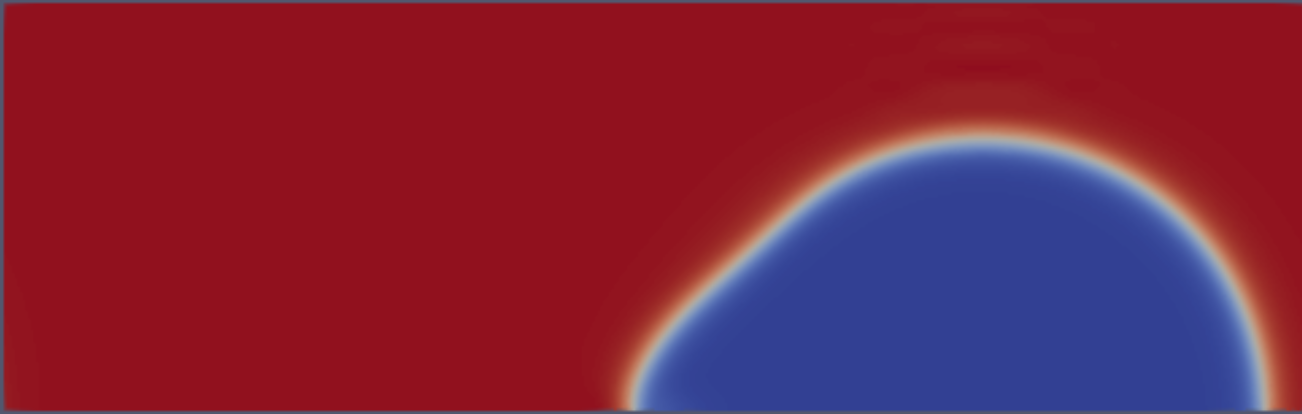}}
\subfigure[$t=0.21$]{
\includegraphics[width=0.24\textwidth]{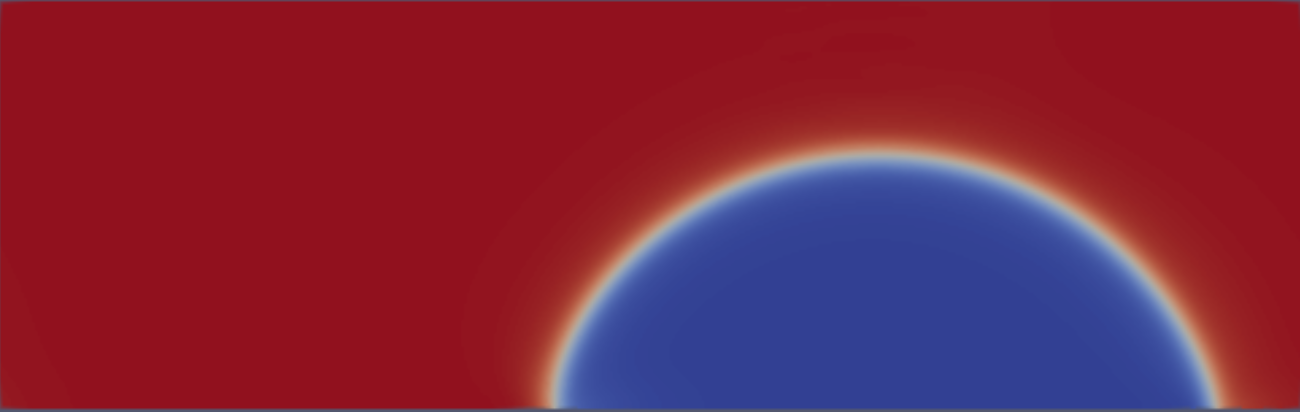}}
\subfigure[$t=0.38$]{
\includegraphics[width=0.24\textwidth]{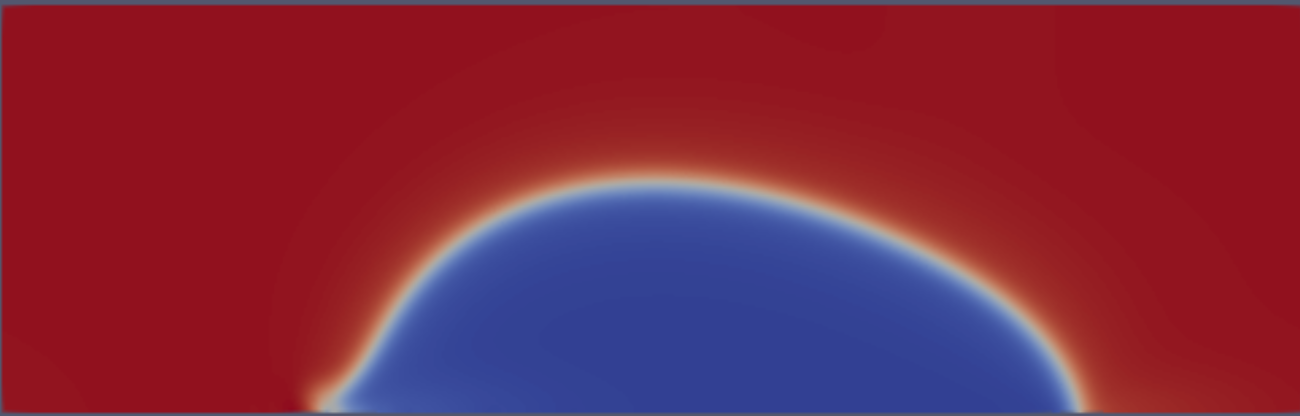}}
}
\centerline{
\subfigure[$t=0.43$]{
\includegraphics[width=0.24\textwidth]{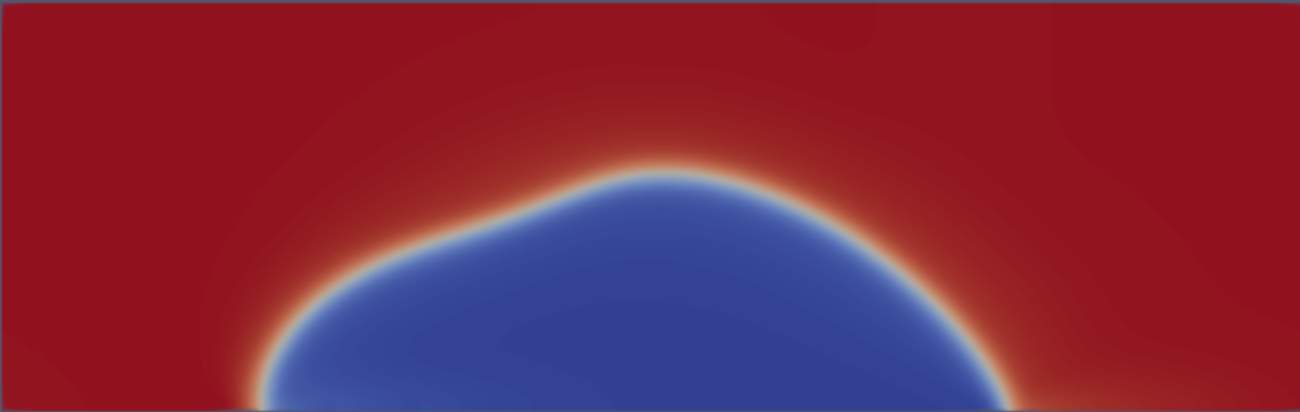}}
\subfigure[$t=0.55$]{
\includegraphics[width=0.24\textwidth]{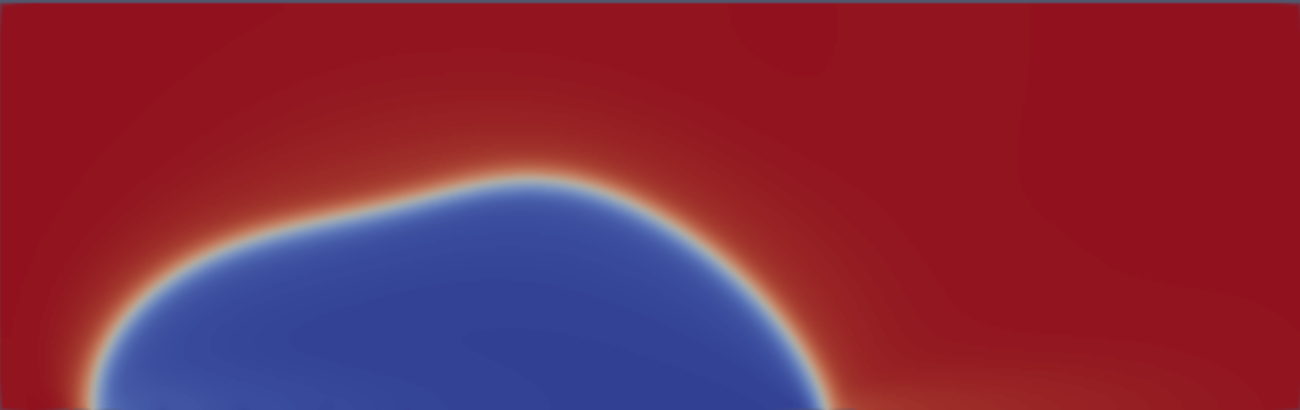}}
\subfigure[$t=0.64$]{
\includegraphics[width=0.24\textwidth]{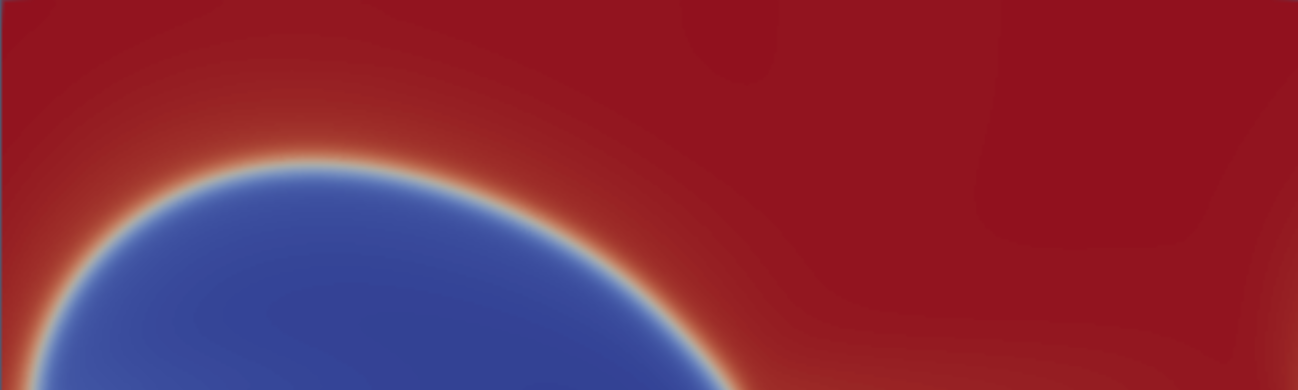}}
\subfigure[$t=1.2$]{
\includegraphics[width=0.24\textwidth]{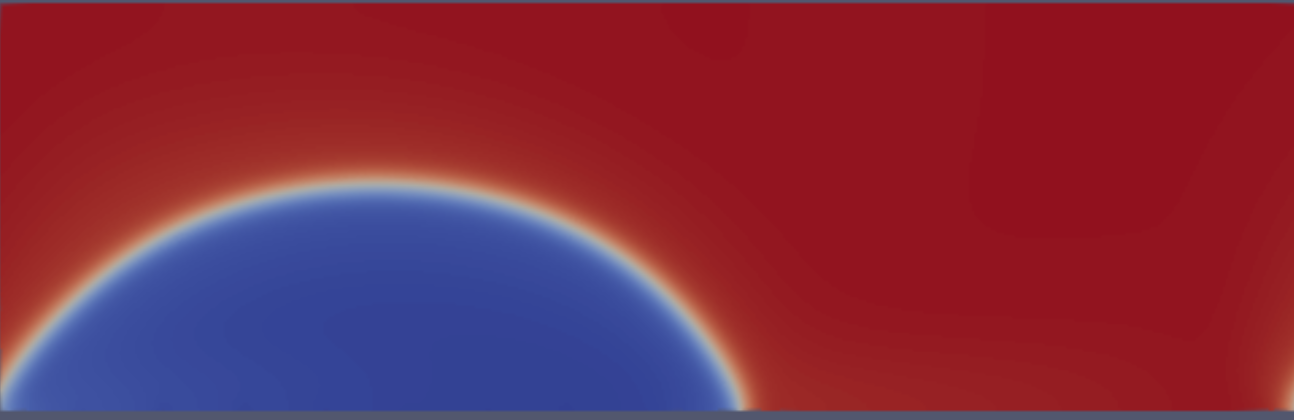}}
}
\caption{\small Drop transport: A temporal sequence of snapshots of the dielectric drop showing its motion on the wall. Shown are the distributions of the phase field function $\phi$ at different time instants.
}
\label{fg_16}
\end{figure}


Figures~\ref{fg_15} and~\ref{fg_16} provide an overview of the electric potential distribution in the domain and the motion of the dielectric drop on the bottom wall.
Shown in Figure~\ref{fg_16} are a temporal sequence of snapshots of the phase field function $\phi(x,y,t)$ in the domain.
One can observe that the dielectric drop moves leftward along the wall due to the interaction with the imposed electric field, and approaches an equilibrium state resting on top of the electrodes.

\subsubsection{Coalescence of Two Dielectric Liquid Drops}
\label{sec:merger}

\begin{figure}[tb]
\centerline{
\includegraphics[width=0.45\textwidth]{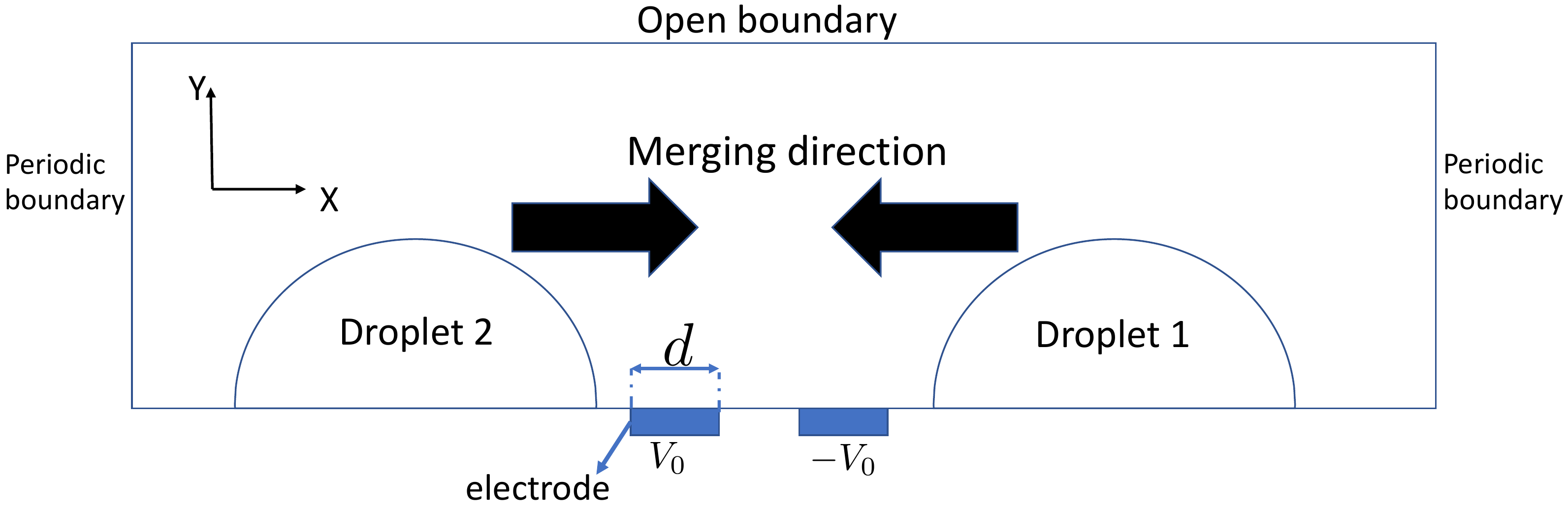}(a)
\includegraphics[width=0.48\textwidth]{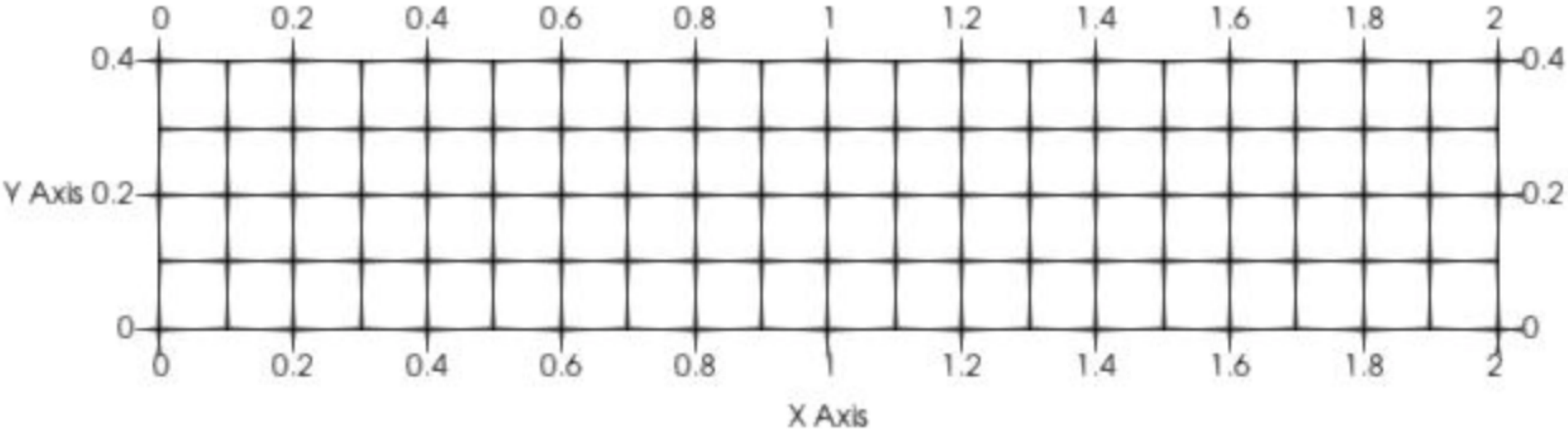}(b)
}
\caption{\small Drop coalescence:
(a) Flow configuration and settings.
(b) Mesh of $80$ quadrilateral spectral elements. 
} 
\label{fg_17} 
\end{figure}

\begin{figure}[tb]
\centering
\includegraphics[width=0.6\textwidth]{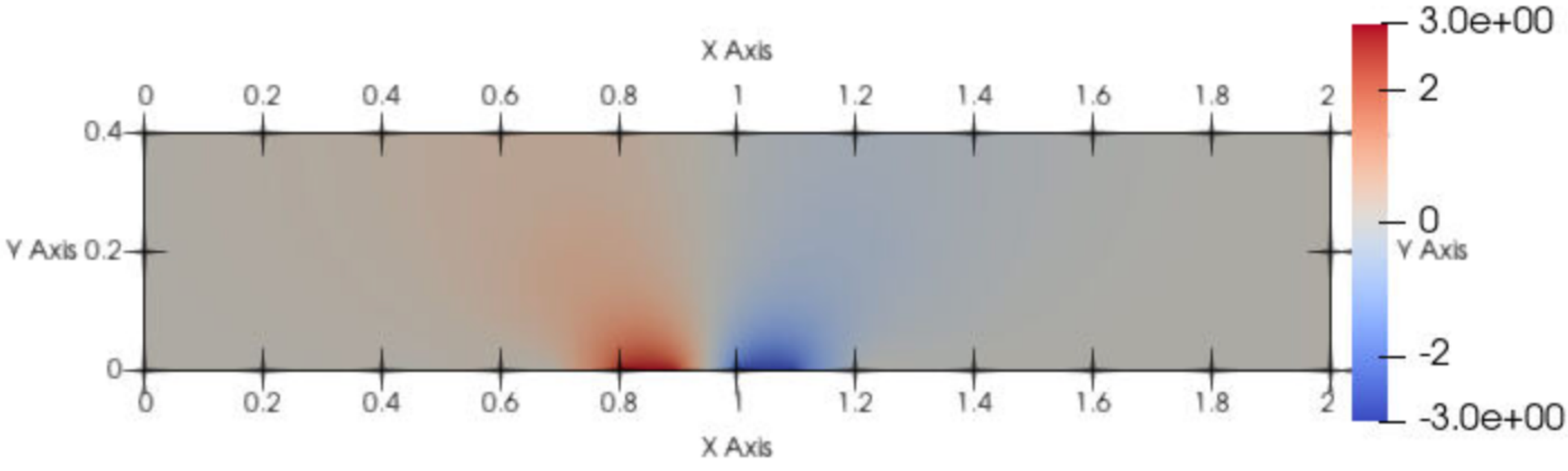}
\caption{\small Drop coalescence: distribution of the electric potential field ($t=0.01$).
} 
\label{fg_18} 
\end{figure}

\begin{figure}[tb]
\centerline{ 
\subfigure[$t=0.01$]{
\includegraphics[width=0.24\textwidth]{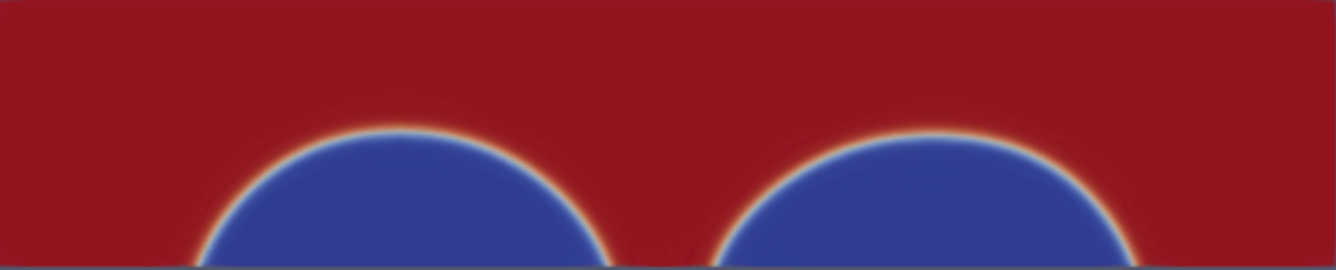}}
\subfigure[$t=0.063$]{
\includegraphics[width=0.24\textwidth]{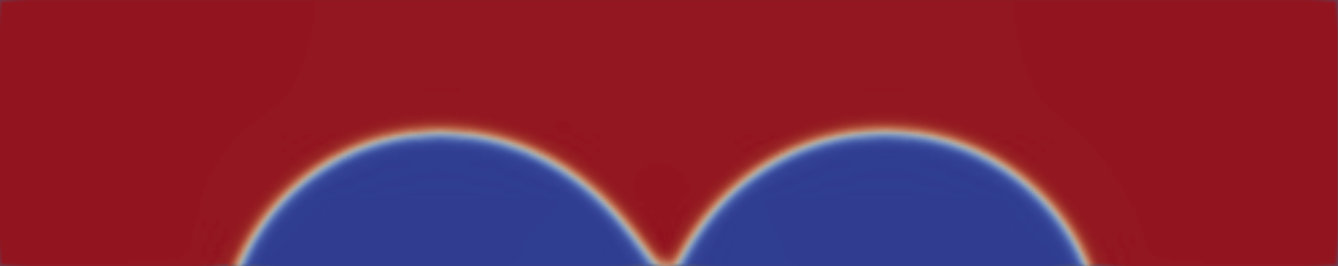}}
\subfigure[$t=0.065$]{
\includegraphics[width=0.24\textwidth]{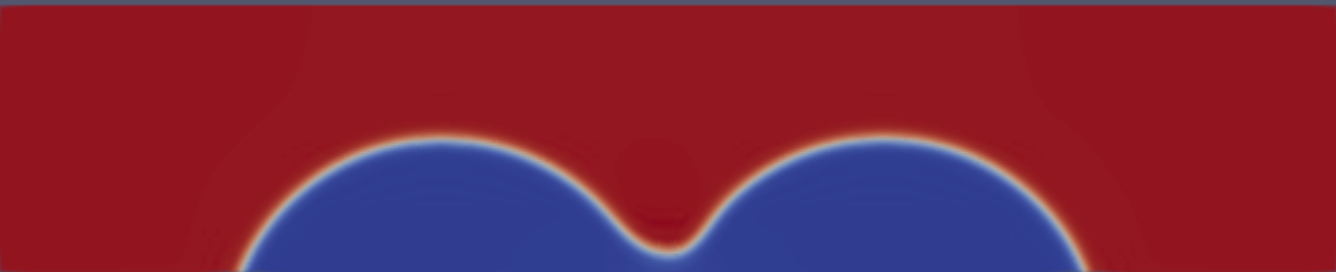}}
\subfigure[$t=0.0675$]{
\includegraphics[width=0.24\textwidth]{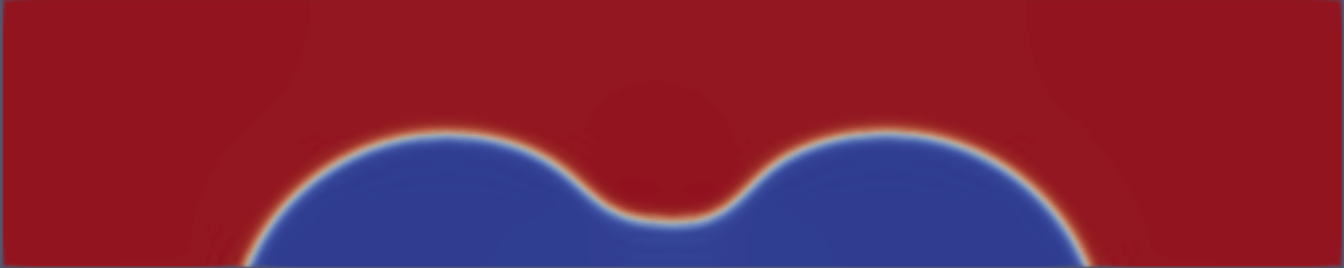}}
}
\centerline{
\subfigure[$t=0.07$]{
\includegraphics[width=0.24\textwidth]{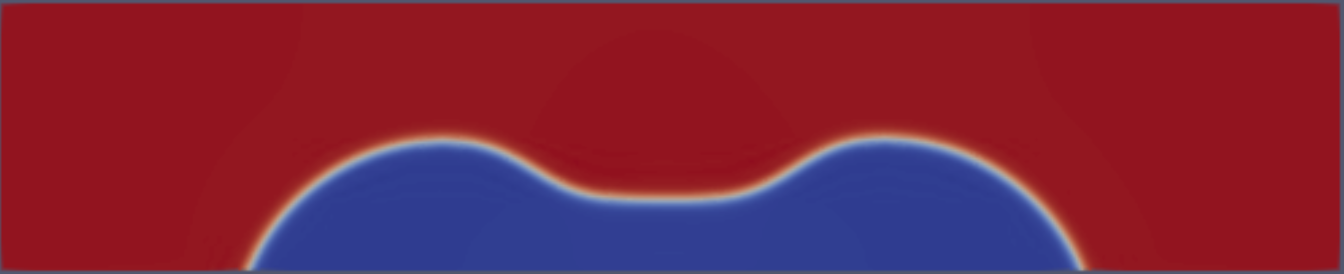}}
\subfigure[$t=0.075$]{
\includegraphics[width=0.24\textwidth]{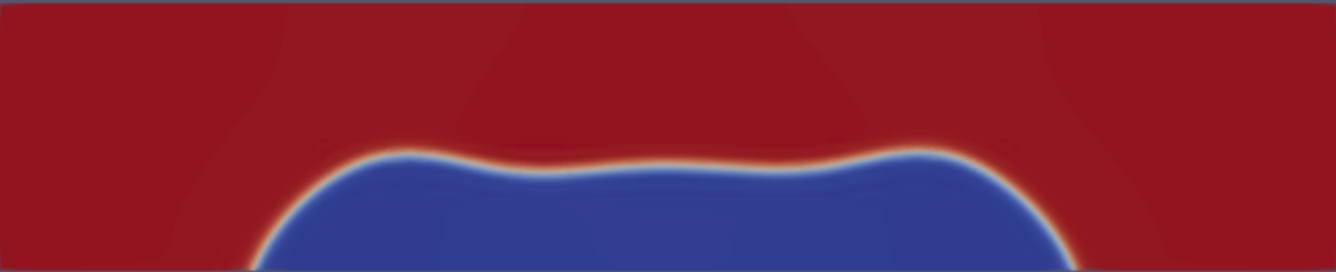}}
\subfigure[$t=0.12$]{
\includegraphics[width=0.24\textwidth]{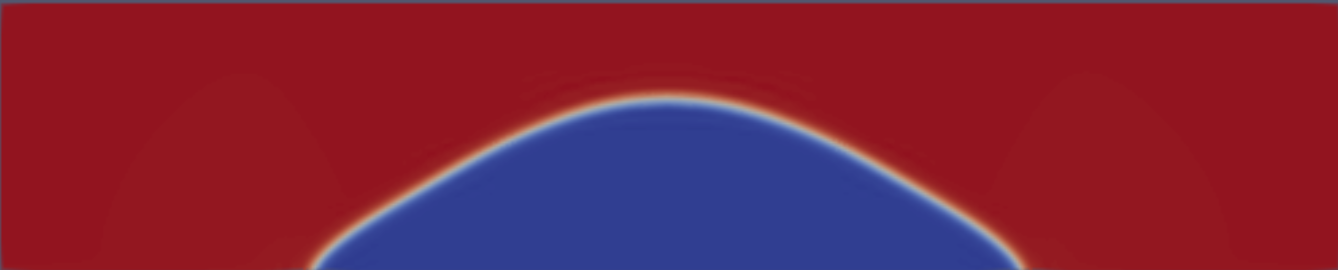}}
\subfigure[$t=0.25$]{
\includegraphics[width=0.24\textwidth]{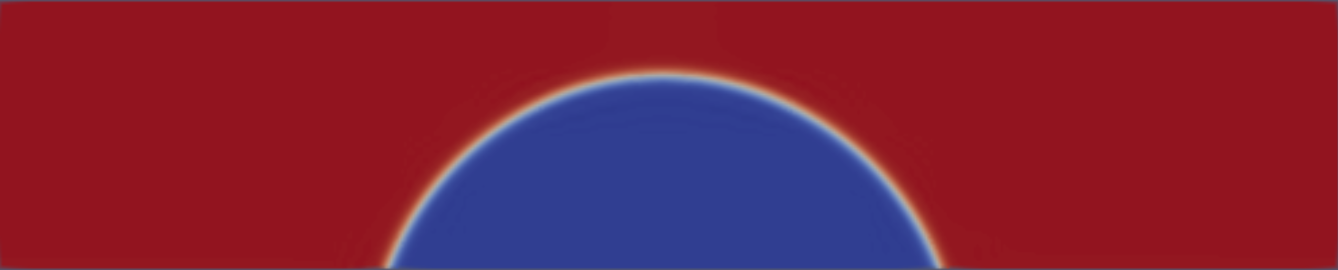}}
}
\caption{\small Drop coalescence: a temporal sequence of snapshots of the phase field distribution, showing the motion and coalescence of the two dielectric fluid drops.
}
\label{fg_19}
\end{figure}

We study the motion and coalescence of two dielectric fluid drops in this subsection. Figure~\ref{fg_17}(a) sketches  the flow configuration and  problem setting. We again consider a rectangular domain, periodic in the horizontal direction, open at the top, and with a solid wall at the bottom.
Two electrodes are embedded in the middle of the bottom wall, whose imposed voltages have the same magnitude but with opposite signs. Two liquid drops of the same dielectric fluid, initially at rest in 
the electrode-free regions of the wall,
are pulled toward each other when the electrodes are turned on, and merge into a single drop. Our goal is to simulate this process with the proposed method.

The simulation settings and the boundary conditions are similar to those employed in Section~\ref{sec:transdrop}.
We employ a computational domain $(x,y)\in\Omega=[0,2L_0]\times [0,\frac25 L_0]$, where $L_0=1mm$, and 
the phase field model given by the equations~\eqref{eq_6}--\eqref{eq_a8}. 
Figure~\ref{fg_17}(b) shows the mesh of $80$ quadrilateral spectral elements (with an element order $12$) employed in the simulations. The two electrodes occupy the following regions on the wall:
$x/L_0\in[0.8,0.9]$ (with voltage $V_0$),
and $x/L_0\in[1,1.1]$ (with voltage $-V_0$), where $V_0=300$volt.
The boundary conditions~\eqref{eq_19}--\eqref{eq_21} are imposed on the bottom wall, with  a static contact angle $\theta_s=75^0$ (measured on the side of the fluid drop). The boundary conditions~\eqref{eq_a19}--\eqref{eq_a25} are imposed on the top boundary ($y/L_0=0.4$).
Periodic conditions are imposed on the horizontal boundaries ($x/L_0=0, 2.0$) for all the dynamic variables.
%
Both drops are assumed to be shaped like a circular cap initially, with radius $R_0$ and their centers located at $(X_1,Y_1)$ and $(X_2,Y_2)$, respectively, as given by
\begin{equation}\small
R_0/L_0=\frac{0.3}{\sin\theta_s}, \quad
X_1/L_0=0.6, \quad
Y_1=-R_0\cos\theta_s, \quad
Y_2/L_0=1.4, \quad
Y_2=-R_0\cos\theta_s.
\end{equation}
The initial phase field distribution is
\begin{equation}\small
\phi(x,y,t=0)=\tanh\frac{\sqrt{(x-X_1)^2+(y-Y_1)^2}-R_0}{\sqrt{2}\eta}+\tanh\frac{\sqrt{(x-X_2)^2+(y-Y_2)^2}-R_0}{\sqrt{2}\eta}-1,
\end{equation}
where $\eta$ is the characteristic interfacial thickness. The initial velocity is set to zero.


We employ the following physical and simulation parameters for this problem:
\begin{equation}\small
\left\{
\begin{split}
&
\text{surface tension:}\ \gamma=1.136\times 10^{-1}kg/s^2; \\
&
\text{densities:}\ \rho_1=\rho_2=129.7kg/m^3;
\ (\text{ambient fluid}\ \rho_1,\ \text{drop}\ \rho_2) \\
&
\text{dynamic viscosities:}\
\text{(ambient fluid)}\ \mu_1=12.048\times 10^{-4}kg/(m\cdot s), \quad
\text{(drop)}\ \mu_2=2\mu_1; \\
&
\text{permittivities:}\ 
\text{(ambient fluid)}\ \epsilon_1=\epsilon_0=8.854\times 10^{-12}F/m, \quad
\text{(drop)}\ \epsilon_2=8.1\epsilon_0; \\
&
\text{static contact angle:}\
\theta_s=75^0\ \text{(measured on the drop side)}; \\
&
\text{model parameters:}\
\eta = 0.007L_0, \quad
\lambda = \frac{3}{2\sqrt{2}}\gamma\eta, \quad
\gamma_1 = 5\times 10^{-6}\frac{L_0^2}{\mu_1}, \quad
\Delta t=1\times 10^{-6}\frac{L_0^2\mu_1}{\epsilon_0 V_0^2}.
\end{split}
\right.
\end{equation}
The physical variables and parameters in the system are normalized based on those constants given in Table~\ref{tab_1}, in which $L_0$ and $\epsilon_0$ are as given above
and we set $\mu_0=\mu_1$ and $V_d=V_0$.

Figure~\ref{fg_18} shows the distribution of the electric potential in the domain, signifying a non-uniform potential and thus a non-uniform electric field. The electric field is stronger near the electrodes, and is weaker in the region farther away from the electrodes. As a result, the net Korteweg-Helmholtz force ($-\frac12(\mbs E\cdot\mbs E)\nabla\epsilon$) acting on the dielectric drops has a direction pointing toward the electrodes, which causes the drops to move inward toward each other. 

Figure~\ref{fg_19} shows a temporal sequence of snapshots of the phase field distribution in the domain. It can be observed that the two drops move along the wall and  merge with each other to form a single drop, which approaches an equilibrium state resting on top of the electrodes.

\section{Concluding Remarks}
\label{sec:summary}


In the current paper we have developed a  method for modeling and simulating multiphase flows consisting of two immiscible incompressible dielectric fluids, and their interactions with  external electric fields in two and three dimensions. We have first presented a thermodynamically-consistent and reduction-consistent formulation based on the phase-field framework for modeling two-phase dielectric fluids. The model honors the mass and momentum conservations, and the second law of thermodynamics.
When only one fluid component is present, the two-phase formulation reduces exactly to that for the single-phase system. In particular, the presented model accommodates an equilibrium solution that is compatible with the requirement of zero velocity based on physics.
This property provides a simpler method for simulating two-phase dielectric systems at equilibrium, 
by solving only a much simplified system consisting of  the phase field equation and the electric potential equation.

We have further presented an efficient semi-implicit type algorithm, together with a spectral-element discretization for 2D and a hybrid Fourier-spectral/spectral-element discretization for 3D in space, for simulating this class of problems. This algorithm allows the computation of different dynamic variables (electric potential, phase field function, pressure, velocity) successively in an uncoupled  fashion. Upon discretization the algorithm  involves only coefficient matrices that are constant and time-independent in the resultant linear algebraic systems, even when the physical properties of the two dielectric fluids (e.g.~the permittivities, densities, viscosities)  are different. 
This property is crucial and enables us to employ the combined Fourier spectral and spectral-element discretization and fast Fourier transforms (FFT) for 3D simulations. 


We have tested the performance of the presented method using several two-phase dielectric problems at equilibrium or in dynamic evolution. The simulation results obtained using the current method have been compared with theoretical models and with experimental measurements. The numerical results signify that the method developed herein can capture the physics well, and that it provides an effective technique for simulating this class of problems.


\section*{Acknowledgments}

This work was partially supported by the US National Science Foundation (DMS-2012415).


\section*{Appendix A: Development of Phase Field Model for Two-Phase Dielectric Flows}

In this appendix we outline the derivation of the phase field model for two-phase dielectric fluids based on the conservation laws and thermodynamic principles. Much of the following development builds upon the works of~\cite{abels2012thermodynamically,dong2014efficient}.

\paragraph{Mass Conservation}

We consider a system of two immiscible incompressible   dielectric fluids, and let $\rho_1,\rho_2$ denote the constant densities of these two pure fluids (without mixing). Consider an arbitrary control volume $V_c$ of the mixture, with mass $M$. Let $M_1$ and $M_2$ denote the mass of these two fluids  within $V_c$. 
Then $\hat{\rho}_1=\frac{M_1}{V_c}$ and $\hat{\rho}_2=\frac{M_2}{V_c}$ denote the densities of the two phases within the mixture. Naturally, we can  introduce 
 the mixture density $\rho$,
\begin{equation}\small
    \rho=\frac{M}{V_c}=\frac{M_1+M_2}{V_c}=\hat{\rho}_1+\hat{\rho}_2\label{bulk_density}
\end{equation}
Let $V_1$ and $V_2$ denote the volume occupied by each pure fluid component with mass $M_1$ and $M_2$. We  assume that when forming the mixture there is no volume loss or volume addition, i.e.
\begin{equation}\label{eq_95}\small
    V_c=V_1+V_2
\end{equation}
We introduce the volume fraction of each fluid by,
$ 
    \phi_i=\frac{V_i}{V_c}
    = \frac{M_i/\rho_i}{M_i/\hat{\rho}_i}
    =\frac{\hat{\rho}_i}{\rho_i}\ (i=1,2).
$ 
Note that $0\leq\phi_i\leq1$. Then equation~\eqref{eq_95} becomes
\begin{equation}\label{eq_97}\small
\phi_1 + \phi_2 = 1.
\end{equation}
We define the phase field variable by,
$ 
    \phi=\phi_1-\phi_2.
$ 

The mass conservation for each phase in the mixture is given by,
\begin{equation}\small
    \frac{\partial \hat{\rho}_i}{\partial t}+\nabla\cdot\hat{\textbf{J}}_i=0, \quad i=1,2 \label{mass_conservation},
\end{equation}
where $\hat{\mbs J}_i$ is the mass flux of phase $i$.
We introduce the velocity $\textbf{u}_i$ of each phase in the mixture  by
$ 
\hat{\textbf{J}}_i=\hat{\rho}_i\textbf{u}_i\ (i=1,2).
$ 
Then equation~\eqref{mass_conservation} is transformed into,
\begin{equation}\small
\dfrac{\partial \hat{\rho}_i}{\partial t}+\nabla\cdot(\hat{\rho}_i\textbf{u}_i)=0,\quad i=1,2.
\end{equation}
We define the bulk mixture velocity $\textbf{u}$ by the volume average of the velocities of the two phases,
\begin{equation}\small\label{eq_102}
\textbf{u}=\phi_1\textbf{u}_1+\phi_2\textbf{u}_2.
\end{equation}
Then it follows that
\begin{equation}\small\label{eq_103}
    \nabla\cdot\textbf{u}=\nabla \cdot\left(\dfrac{\hat{\rho}_1}{\rho_1}\textbf{u}_1+\dfrac{\hat{\rho}_2}{\rho_2}\textbf{u}_2\right)=\nabla\cdot\left(\dfrac{\hat{\textbf{J}}_1}{\rho_1}+\dfrac{\hat{\textbf{J}}_2}{\rho_2}\right)=-\dfrac{\partial}{\partial t}\left(\dfrac{\hat{\rho}_1}{\rho_1}+\dfrac{\hat{\rho}_2}{\rho_2}\right)
    =-\frac{\partial}{\partial t}(\phi_1+\phi_2)
    =-\dfrac{\partial 1}{\partial t}=0,
\end{equation}
where equation~\eqref{eq_97} has been used.

Equation~\eqref{eq_103} indicates the bulk mixture velocity as defined above is divergence free (see also~\cite{ding2007diffuse,abels2012thermodynamically}). One can also use the mass fraction to define the bulk velocity (see e.g.~\cite{LowengrubT1998}). However, in that case the bulk velocity will not be divergence free. In the current work we employ the volume-averaged velocity as the bulk mixture velocity, as given by~\eqref{eq_102}.

Finally, the mass conservation in terms of the bulk density $\rho$ is, by adding equation \eqref{mass_conservation} for $i=1,2$,
\begin{equation}\small
    \dfrac{\partial\rho}{\partial t}+\textbf{u}\cdot\nabla\rho=-\nabla\cdot\Tilde{\textbf{J}}\label{mass_conservation_bulk},
\end{equation}
where 
$ 
    \Tilde{\textbf{J}}= (\hat{\textbf{J}}_1-\hat{\rho}_1\textbf{u})+(\hat{\textbf{J}}_2-\hat{\rho}_2\textbf{u}).
$ 
$\Tilde{\textbf{J}}$ denotes the total difference of the mass flux of different phases with respect to the bulk.
It will be determined by an constitutive relation based on the energy inequality.
Note that equation (\ref{bulk_density}) implies $\rho=\frac{\rho_1+\rho_2}{2}+\frac{\rho_1-\rho_2}{2}\phi$. So equation~\eqref{mass_conservation_bulk} can be transformed into,
\begin{equation}\small\label{eq_106}
    \dfrac{\partial\phi}{\partial t}+\textbf{u}\cdot\nabla\phi=-\dfrac{2}{\rho_1-\rho_2}\nabla\cdot\Tilde{\textbf{J}}.
\end{equation}

\paragraph{Momentum Conservation}

Following~\cite{gurtin1996two,dong2014efficient}, we assume that the inertia and the kinetic energy due to the relative motion of each fluid phase with respect to  the bulk motion are negligible.
The conservation of momentum for each fluid phase is represented by,
\begin{equation}\small
\frac{\partial (\hat{\rho}_i\textbf{u}_i)}{\partial t}+\nabla \cdot(\hat{\rho}_i\textbf{u}_i\textbf{u}_i)=\nabla \cdot \textbf{T}_i+\bm\pi_i,
\quad i=1,2,
\end{equation}
where $\textbf{T}_i$ is the stress tensor of the phase $i$, and $\bm\pi_i$ ($i=1,2$) represents the interaction body force, with $\bm\pi_1+\bm\pi_2=0$.

We rewrite the above equation into,
\begin{equation}\small
\begin{aligned}
\frac{\partial (\hat{\rho}_i\textbf{u})}{\partial t}+\nabla \cdot(\hat{\rho}_i\textbf{u}\textbf{u})&+\frac{\partial (\hat{\rho}_i(\textbf{u}_i-\textbf{u}))}{\partial t}+\nabla \cdot(\hat{\rho}_i(\textbf{u}_i-\textbf{u})(\textbf{u}_i-\textbf{u}))\\
&+\nabla \cdot(\hat{\rho}_i (\textbf{u}_i-\textbf{u})\textbf{u})+\nabla \cdot (\hat{\rho}_i \textbf{u}(\textbf{u}-\textbf{u}_i))=\nabla \cdot \textbf{T}_i+\bm\pi_i.
\end{aligned}
\end{equation}
We omit the third and the fourth terms on the left hand side (LHS)
based on the assumption that the inertia and the kinetic energy of the differential motion relative to the bulk are negligible. 
We move the term $\nabla \cdot (\rho_i \textbf{u}(\textbf{u}-\textbf{u}_i))$ to the right hand side (RHS) and incorporate it into the $\nabla\mbs T_i$ term to get,
\begin{equation}\small\label{eq_109}
\frac{\partial (\hat{\rho}_i\textbf{u})}{\partial t}+\nabla \cdot(\hat{\rho}_i\textbf{u}\textbf{u})+\nabla \cdot(\textbf{J}_i\textbf{u})=\nabla \cdot \tilde{\textbf{T}}_i+\bm\pi_i, 
\quad i=1,2,
\end{equation}
where $\textbf{J}_i=\hat{\textbf{J}}_i-\hat{\rho}_i\textbf{u}$ and
$\tilde {\textbf{T}}_i=\textbf{T}_i-\hat{\rho}_i\textbf{u}(\textbf{u}-\textbf{u}_i)$.
Then we sum  up equation~\eqref{eq_109} for all the phases,
\begin{equation}\small
\frac{\partial (\rho \textbf{u})}{\partial t}+\nabla \cdot(\rho \textbf{u}\textbf{u})+\nabla \cdot(\Tilde{\textbf{J}}\textbf{u})=\nabla \cdot \textbf{T},
\end{equation}
where $\mbs T$ is a stress tensor with $\textbf{T}=\tilde{\textbf{T}}_1+\tilde{\textbf{T}}_2$, and we have used $\bm\pi_1+\bm\pi_2=0$. In light of~\eqref{mass_conservation_bulk}, this equation can be transformed into,
\begin{equation}\small
\rho\left(\frac{\partial  \textbf{u}}{\partial t}+\textbf{u}\cdot \nabla \textbf{u}\right)+\Tilde{\textbf{J}} \cdot \nabla \textbf{u}=\nabla \cdot \textbf{T}.
\label{momentum_equation_bulk}
\end{equation}

We assume that the 
the stress tensor $\textbf{T}$ is symmetric, and re-write it into
\begin{equation}\small\label{eq_112}
\mbs T = \frac13(\text{tr}\mbs T)\mbs I + \mbs S = -p\mbs I + \mbs S,
\end{equation}
where $\mbs I$ denotes the identity tensor, $\mbs S$ is a trace-free symmetric tensor, and $p=-\frac13\text{tr}\mbs T$ will be called the pressure.
Then equation~\eqref{momentum_equation_bulk} becomes,
\begin{equation}\small
\rho\left(\frac{\partial  \textbf{u}}{\partial t}+\textbf{u}\cdot \nabla \textbf{u}\right)+\Tilde{\textbf{J}} \cdot \nabla \textbf{u}= -\nabla p+\nabla \cdot \textbf{S}.
\label{eq_113}
\end{equation}
The tensor $\mbs S$ will be determined from a constitutive relation based on the energy inequality.

\paragraph{Quasi-Static Maxwell Equations}

We focus on a system of dielectrc fluids, which are non-conductive and contain no free electric charge. The characteristic velocity in the system is negligible compared with the speed of light. On the other hand, we would like to take into account the fluid motion and the momentum transport. So this is an electro quasi-static system~\cite{de2006electrodynamics}.

The quasi-static Maxwell equations are given by, 
\begin{subequations}\label{eq_a112}
\begin{align}\small
    &\nabla\cdot\textbf{D}=0, \label{Maxwell}\\
    &\nabla\times\textbf{E}=\textbf{0},
    \label{eq_112b}
\\
    &\dfrac{\partial \textbf{D}}{\partial t}=\nabla\times\textbf{H},
    \label{eq_112c}
\end{align}
\end{subequations}
where $\textbf{E}$ is the electric field, $\textbf{D}$ is electric displacement field ($\mbs D=\epsilon\mbs E$, with $\epsilon$ denoting the material permittivity), and $\textbf{H}$ is the magnetizing field. 
Equation \eqref{Maxwell} indicates that there is no free charge in the system. Equation~\eqref{eq_112b} allows us to introduce the electric potential $V(\mbs x)$ by 
\begin{equation}\small\label{eq_115}
\textbf{E}=\nabla V. 
\end{equation} 
Note that the equations \eqref{Maxwell} and \eqref{eq_112b} alone are sufficient  to determine the electric field. Equation \eqref{eq_112c} will not be solved in numerical simulations. But this equation  plays an important role in deriving the energy balance relation.  The magnetic field $\mbs H$ is weak based on the quasi-static assumption.

\paragraph{Energy Inequality and Constitutive Relations}

Let us now determine the forms of $\tilde{\mbs J}$ and $\mbs S$ involved in the mass/momentum balance equations based on the second law of thermodynamics.
We define the total energy of the system by,
\begin{equation}\small
E(t)=\int_{\Omega} \left[\dfrac{1}{2}\rho \textbf{u}\cdot \textbf{u}+F(\phi,\nabla \phi)+\dfrac{1}{2}\textbf{D}\cdot\textbf{E}\right]dV+\int_{\partial \Omega_{s}}\Theta(\phi) dS,
\end{equation}
where $\Omega$ is an arbitrary domain, $\partial\Omega_s$ denotes the wall boundary, $\frac12\mbs D\cdot\mbs E$ is the electric energy density,  $F(\phi,\nabla\phi)$ is the phase-field free energy density function  (see e.g.~\eqref{double_well}), and $\Theta(\phi)$ denotes a wall energy density to account for the contact angle effect. We assume that mixture permittivity is a function of the phase field function, $\epsilon = \epsilon(\phi)$.
By using equations \eqref{mass_conservation_bulk},
\eqref{momentum_equation_bulk} and \eqref{eq_a112}, we can derive
\begin{align}\small\label{eq_116}
 \frac{dE}{dt}=&\int_{\Omega} \dfrac{\rho_1-\rho_2}{2}\nabla\left[\dfrac{\partial F}{\partial \phi}-\nabla\cdot\dfrac{\partial F}{\partial\nabla\phi}-\dfrac{\epsilon'}{2}\textbf{E}\cdot\mbs E\right]\cdot\Tilde{\textbf{J}}-
\int_{\Omega}\left(\textbf{T}+\dfrac{\partial F}{\partial\nabla\phi}\otimes\nabla\phi-\boldsymbol{\sigma}_M\right):\nabla\textbf{u}\notag \\
 & -\int_{\partial \Omega}\left[\dfrac{\rho_1-\rho_2}{2}\left(\dfrac{\partial F}{\partial \phi}-\nabla\cdot\dfrac{\partial F}{\partial\nabla\phi}-\dfrac{\epsilon'}{2}\textbf{E}\cdot\mbs E\right)\Tilde{\textbf{J}}\cdot\textbf{n}-\dfrac{1}{2}(\textbf{u}\cdot\textbf{u})\Tilde{\textbf{J}}\cdot\textbf{n}\right]\notag\\
 &+\int_{\partial \Omega}\left[\left( \mbs T - \bm\sigma_M
+ \frac{\partial F}{\partial\nabla\phi}\otimes\nabla\phi
\right)\cdot\textbf{n}-F\textbf{n}-\dfrac{1}{2}(\textbf{u}\cdot\textbf{u})\textbf{n}\right]\cdot\textbf{u}\notag \\
 &+\int_{\partial \Omega}\left(\dfrac{\partial F}{\partial \nabla\phi}\cdot\textbf{n}\right)\dfrac{\partial \phi}{\partial t}+\int_{\partial \Omega_{s}}\Theta'(\phi)\dfrac{\partial \phi}{\partial t}
 -\int_{\partial \Omega}(\textbf{E}\times\textbf{H})\cdot\textbf{n}.
\end{align}
where $\epsilon'=\frac{d\epsilon}{d\phi}$, and $\boldsymbol{\sigma}_M$ is the Maxwell stress tensor~\cite{landau2013electrodynamics}, given by
\begin{equation}\small
\boldsymbol{\sigma}_M=\textbf{D}\otimes\textbf{E}-\dfrac{\epsilon}{2}(\textbf{E}\cdot\textbf{E})\textbf{I}.
\end{equation}

The second law of thermodynamics dictates that in the absence of  external forces (including surface forces acting on the boundary)
the system should be dissipative.
This means that the contributions of the volume integral terms involved in the above equation to $dE/dt$ should always be non-positive, while the contributions of the surface integral terms can be controlled if appropriate boundary conditions are imposed.
We would like to choose the constitutive relations about $\tilde{\mbs J}$ and $\mbs T$ such that the requirements of  the second law of thermodynamics are satisfied.

To ensure the non-positivity of the first volume integral on the RHS of~\eqref{eq_116}, we choose the following constitutive relation,
\begin{equation}\small\label{eq_118}
    \Tilde{\textbf{J}}=-\gamma_1\dfrac{\rho_1-\rho_2}{2}\nabla\left(\dfrac{\partial F}{\partial \phi}-\nabla\cdot\dfrac{\partial F}{\partial\nabla\phi}-\dfrac{\epsilon'}{2}\textbf{E}\cdot\mbs E\right),
\end{equation}
where $\gamma_1\geqslant 0$ is a non-negative constant or function.

Noting the symmetry of the tensors $\mbs T$ and $\bm{\sigma}_M$, the second volume integral on the RHS of~\eqref{eq_116} can be transformed into,
\begin{equation}\small
\begin{split}
&
-\int_{\Omega}\left( \mbs T - \bm\sigma_M
+ \frac{\partial F}{\partial\nabla\phi}\otimes\nabla\phi
\right):\nabla\mbs u \\
&=
-\int_{\Omega}\left[ \mbs T - \bm\sigma_M
+ \frac12\left( \frac{\partial F}{\partial\nabla\phi}\otimes\nabla\phi 
+ \nabla\phi\otimes\frac{\partial F}{\partial\nabla\phi}
\right)
\right]:\frac12\left(\nabla\mbs u
+ \nabla\mbs u^T\right) \\
&\quad
- \int_{\Omega} \frac12\left( \frac{\partial F}{\partial\nabla\phi}\otimes\nabla\phi 
- \nabla\phi\otimes\frac{\partial F}{\partial\nabla\phi}
\right) : 
\frac12\left(\nabla\mbs u
- \nabla\mbs u^T\right)\\
&=
-\int_{\Omega}\left[ \mbs S - \tilde{\bm\sigma}_M
+\tilde{\bm{\mathcal F}}_{\phi}
\right]:\frac12\left(\nabla\mbs u
+ \nabla\mbs u^T\right) \\
&\quad
- \int_{\Omega} \frac12\left( \frac{\partial F}{\partial\nabla\phi}\otimes\nabla\phi 
- \nabla\phi\otimes\frac{\partial F}{\partial\nabla\phi}
\right) : 
\frac12\left(\nabla\mbs u
- \nabla\mbs u^T\right)
\end{split}
\end{equation}
where we have used equations~\eqref{eq_103}, \eqref{eq_112} and~\eqref{eq_a112}, and
\begin{equation}\small\label{eq_a120}
\left\{
\begin{split}
&
\tilde{\bm\sigma}_M = \bm\sigma_M - \frac13\left[\text{tr}\bm\sigma_M\right]\mbs I,
\\
&
\tilde{\bm{\mathcal F}}_{\phi} = \bm{\mathcal F}_{\phi} - \frac13\left[\text{tr}\bm{\mathcal F}_{\phi}\right] \mbs I, \quad
\bm{\mathcal F}_{\phi} = 
\frac12\left( \frac{\partial F}{\partial\nabla\phi}\otimes\nabla\phi 
+ \nabla\phi\otimes\frac{\partial F}{\partial\nabla\phi}
\right).
\end{split}
\right.
\end{equation}
We choose the following constitutive relation to ensure its non-positivity,
\begin{align}\small
&
\mbs S - \tilde{\bm\sigma}_M
+ \tilde{\bm{\mathcal F}}_{\phi}
 = \mu(\nabla\mbs u + \nabla\mbs u^T), \label{eq_120} \\
&
\frac{\partial F}{\partial\nabla\phi}\otimes\nabla\phi 
- \nabla\phi\otimes\frac{\partial F}{\partial\nabla\phi} = 0. \label{eq_121}
\end{align}
where $\mu\geqslant 0$ is a non-negative constant or function. 
Equation~\eqref{eq_121} is a condition that the free energy density function $F(\phi,\nabla\phi)$ must satisfy. Equation~\eqref{eq_120} provides the tensor $\mbs S$,
\begin{equation}\small\label{eq_122}
\mbs S = \mu(\nabla\mbs u + \nabla\mbs u^T) 
+ \left(\bm\sigma_M - \frac13[\text{tr}\bm\sigma_M] \mbs I\right)
-\left(\bm{\mathcal F}_{\phi} - \frac13[\text{tr}\bm{\mathcal F}_{\phi}]\mbs I
\right).
\end{equation}

In light of equations~\eqref{eq_118}, \eqref{eq_122}, \eqref{eq_106}, \eqref{eq_113}, \eqref{eq_103}, \eqref{Maxwell} and~\eqref{eq_115},
we can write down the system of governing equations as follows,
\begin{subequations}\label{eq_124}
\begin{align}\small
    &\dfrac{\partial \phi}{\partial t}+\textbf{u}\cdot\nabla\phi=\nabla \cdot\left[ \gamma_1\nabla\left(\dfrac{\partial F}{\partial \phi}-\nabla\cdot\dfrac{\partial F}{\partial\nabla\phi}-\dfrac{\epsilon'}{2}\mbs E\cdot\mbs E\right)\right],\\
    &\rho\left(\dfrac{\partial \mbs u}{\partial t}+\mbs u\cdot\nabla\mbs u\right)+\Tilde{\mbs J}\cdot\nabla\mbs u=-\nabla\cdot \left( \dfrac{\partial F}{\partial \nabla \phi}\otimes \nabla\phi \right)
    -\frac{\nabla\epsilon}{2}(\mbs E\cdot\mbs E)
    +\nabla\cdot\left[\mu\left(\nabla\mbs u + \nabla\mbs u^T\right)\right]-\nabla {\mathcal P}, \\
    &\nabla\cdot\textbf{u}=0,\\
    &\nabla\cdot(\epsilon\nabla V)=0.
\end{align}
\end{subequations}
In the above equations, 
\begin{equation}\small
\mathcal P = p + \frac13\text{tr}\bm\sigma_M
-\frac13\text{tr}\bm{\mathcal F}_{\phi},
\end{equation}
and we have used
\begin{equation}\small
\nabla\cdot\boldsymbol{\sigma}_M
=\nabla\cdot\left[\mbs D\otimes\mbs E-\dfrac{\epsilon}{2}(\mbs E\cdot \mbs E) \mbs I\right]
=-\frac{\nabla\epsilon}{2}(\mbs E\cdot\mbs E).
\end{equation}
This set of equations constitutes the phase field model that describes the motion of a system of two immiscible incompressible dielectric fluids. This model is thermodynamically consistent.

In this model the form for the mixing energy density function $F(\phi,\nabla\phi)$ is still to be chosen, and it must satisfy the condition~\eqref{eq_121}. 
If we choose $F(\phi,\nabla\phi)$ based on equation~\eqref{double_well}, which satisfies the condition~\eqref{eq_121}, then the system~\eqref{eq_124} will be reduced to the system consisting of equations~\eqref{eq_6}--\eqref{eq_9}.

\bibliographystyle{plain}
\bibliography{obc,interface,multiphase,nphase,contact_line,mypub,engstab}

\end{document}